\documentclass[11pt]{amsart}
\usepackage{stmaryrd}
\newtheorem{theorem}{Theorem}[section]
\newtheorem{proposition}[theorem]{Proposition}
\newtheorem{lemma}[theorem]{Lemma}
\newtheorem{corollary}[theorem]{Corollary}
\theoremstyle{definition}
\newtheorem{definition}[theorem]{Definition}

\newcommand{\tracefreeRic}{\overset{\text{\rm o}}{\text{\rm Ric}}{}}
\newcommand{\tracefreeH}{\overset{\text{\rm o}}{H}{}}

\begin{document}

\title[Ricci flow on manifolds with positive isotropic curvature]{Ricci flow with surgery on manifolds with positive isotropic curvature}
\author{Simon Brendle}
\address{Department of Mathematics \\ Columbia University \\ New York, NY 10027}
\begin{abstract}
We study the Ricci flow for initial metrics with positive isotropic curvature (strictly PIC for short). 

In the first part of this paper, we prove new curvature pinching estimates which ensure that blow-up limits are uniformly PIC in all dimensions. Moreover, in dimension $n \geq 12$, we show that blow-up limits are weakly PIC2. This can be viewed as a higher-dimensional version of the fundamental Hamilton-Ivey pinching estimate in dimension $3$.

In the second part, we develop a theory of ancient solutions which have bounded curvature; are $\kappa$-noncollapsed; are weakly PIC2; and are uniformly PIC. This is an extension of Perelman's work; the additional ingredients needed in the higher dimensional setting are the differential Harnack inequality for solutions to the Ricci flow satisfying the PIC2 condition, and a rigidity result due to Brendle-Huisken-Sinestrari for ancient solutions that are uniformly PIC1. 

In the third part of this paper, we prove a Canonical Neighborhood Theorem for the Ricci flow with initial data with positive isotropic curvature, which holds in dimension $n \geq 12$. This relies on the curvature pinching estimates together with the structure theory for ancient solutions. This allows us to adapt Perelman's surgery procedure to this situation. As a corollary, we obtain a topological classification of all compact manifolds with positive isotropic curvature of dimension $n \geq 12$ which do not contain non-trivial incompressible $(n-1)$-dimensional space forms.
\end{abstract}
\maketitle

\tableofcontents 

\section{Introduction}

Our goal in this paper is to study the formation of singularities under the Ricci flow for initial metrics with positive isotropic curvature. Positive isotropic curvature is a natural curvature condition which makes sense in dimension $n \geq 4$. It was introduced in the work of Micallef and Moore \cite{Micallef-Moore} in their study of minimal two-spheres in Riemannian manifolds. Variants of this condition play a central role in the proof of the Differentiable Sphere Theorem \cite{Brendle-Schoen}. We first recall the relevant definitions:

\begin{definition}
\label{pic.and.its.variants}
(i) We denote by PIC the set of all algebraic curvature tensors which have nonnegative isotropic curvature in the sense that $R(\varphi,\bar{\varphi}) \geq 0$ for all complex two-forms of the form $\varphi = (e_1+ie_2) \wedge (e_3+ie_4)$, where $\{e_1,e_2,e_3,e_4\}$ is an orthonormal four-frame. \\ 
(ii) We denote by PIC1 the set of all algebraic curvature tensors satisfying $R(\varphi,\bar{\varphi}) \geq 0$ for all complex two-forms of the form $\varphi = (e_1+ie_2) \wedge (e_3+i\lambda e_4)$, where $\{e_1,e_2,e_3,e_4\}$ is an orthonormal four-frame and $\lambda \in [0,1]$. \\ 
(iii) We denote by PIC2 the set of all algebraic curvature tensors satisfying $R(\varphi,\bar{\varphi}) \geq 0$ for all complex two-forms of the form $\varphi = (e_1+i\mu e_2) \wedge (e_3+i\lambda e_4)$, where $\{e_1,e_2,e_3,e_4\}$ is an orthonormal four-frame and $\lambda,\mu \in [0,1]$. 
\end{definition}

Note that $\text{\rm PIC2} \subset \text{\rm PIC1} \subset \text{\rm PIC}$. The curvature tensor of a Riemannian manifold $M$ lies in the PIC1 cone if and only if the curvature tensor of $M \times \mathbb{R}$ lies in the PIC cone. Similarly, the curvature tensor of $M$ lies in the PIC2 cone if and only if the curvature tensor of $M \times \mathbb{R}^2$ lies in the PIC cone.

The significance of the curvature conditions above stems from the fact that they are all preserved by the Ricci flow. For an initial metric that is weakly PIC2, the subsequent solution of the Ricci flow satisfies a differential Harnack inequality (cf. \cite{Hamilton3}, \cite{Brendle2}). For an initial metric that is strictly PIC1, it was shown in \cite{Brendle1} that the Ricci flow will converge to a metric of constant curvature after rescaling (see \cite{Bohm-Wilking},\cite{Brendle-Schoen},\cite{Huisken},\cite{Margerin1},\cite{Margerin2},\cite{Margerin3},\cite{Nishikawa} for earlier work on the subject). For an initial metric that is strictly PIC, it has been conjectured that the Ricci flow should only form so-called neck-pinch singularities. For $n=4$, this was proved in a fundamental paper by Hamilton \cite{Hamilton5} (see also \cite{Chen-Zhu},\cite{Chen-Tang-Zhu}). Our goal in this paper is to confirm the conjecture for $n \geq 12$.

A key step in our analysis is a new curvature pinching estimate in higher dimensions. By work of Hamilton \cite{Hamilton1},\cite{Hamilton2}, the curvature tensor satisfies the evolution equation 
\[D_t R = \Delta R + Q(R).\] 
Here, $D_t$ denotes the covariant time derivative, and $Q(R)$ is a quadratic expression in the curvature tensor. More precisely, $Q(R) = R^2 + R^\#$, where $R^2$ and $R^\#$ are defined by 
\[(R^2)_{ijkl} = \sum_{p,q=1}^n R_{ijpq} R_{klpq}\] 
and 
\[(R^\#)_{ijkl} = 2 \sum_{p,q=1}^n (R_{ipkq} R_{jplq} - R_{iplq} R_{jpkq}).\] 
Note that the definitions of $R^2$, $R^\#$, and $Q(R)$ make sense for any algebraic curvature tensor $R$. In fact, the definitions even make sense if $R$ does not satisfy the first Bianchi identity. It is sometimes convenient to consider curvature-type tensors that do not satisfy the first Bianchi identity (see Section \ref{first.family.of.cones} below); however, unless stated otherwise, we will assume that the first Bianchi identity is satisfied.

In order to prove pinching estimates for the Ricci flow, we need to analyze the Hamilton ODE $\frac{d}{dt} R = Q(R)$ on the space of algebraic curvature tensors. Our first main result is a pinching estimate for the Hamilton ODE: 

\begin{theorem}
\label{main.theorem}
Assume that $n \geq 12$. Let $\mathcal{K}$ be a compact set of algebraic curvature tensors in dimension $n$ which is contained in the interior of the PIC cone, and let $T>0$ be given. Then there exist a small positive real number $\theta$, a large positive real number $N$, an increasing concave function $f>0$ satisfying $\lim_{s \to \infty} \frac{f(s)}{s} = 0$, and a continuous family of closed, convex, $O(n)$-invariant sets $\{\mathcal{F}_t: t \in [0,T]\}$ such that the family $\{\mathcal{F}_t: t \in [0,T]\}$ is invariant under the Hamilton ODE $\frac{d}{dt} R = Q(R)$; $\mathcal{K} \subset \mathcal{F}_0$; and 
\begin{align*} 
\mathcal{F}_t &\subset \{R: R - \theta \, \text{\rm scal} \, \text{\rm id} \owedge \text{\rm id} \in \text{\rm PIC}\} \\ 
&\cap \{R: \text{\rm Ric}_{11}+\text{\rm Ric}_{22} - \theta \, \text{\rm scal} + N \geq 0\} \\ 
&\cap \{R: R + f(\text{\rm scal}) \, \text{\rm id} \owedge \text{\rm id} \in \text{\rm PIC2}\}
\end{align*} 
for all $t \in [0,T]$.
\end{theorem}

Here, $\owedge$ denotes the Kulkarni-Nomizu product. More precisely, if $A$ and $B$ are symmetric bilinear forms, then $(A \owedge B)_{ijkl} = A_{ik} B_{jl} - A_{il} B_{jk} - A_{jk} B_{il} + A_{jl} B_{ik}$.

Via Hamilton's PDE-ODE principle (cf. \cite{Chow-Lu}, Theorem 3, or \cite{Chow-et-al}, Theorem 10.16), Theorem \ref{main.theorem} gives curvature pinching estimates for solutions to the Ricci flow starting from initial metrics with positive isotropic curvature: 

\begin{corollary}
\label{pinching.estimates.for.ricci.flow}
Let $(M,g_0)$ be a compact manifold of dimension $n \geq 12$ with positive isotropic curvature, and let $g(t)$ denote the solution to the Ricci flow with initial metric $g_0$. Then there exist a small positive real number $\theta$, a large positive real number $N$, and an increasing concave function $f$ satisfying $\lim_{s \to \infty} \frac{f(s)}{s} = 0$ such that the curvature tensor of $(M,g(t))$ satisfies $R - \theta \, \text{\rm scal} \, \text{\rm id} \owedge \text{\rm id} \in \text{\rm PIC}$, $\text{\rm Ric}_{11}+\text{\rm Ric}_{22} - \theta \, \text{\rm scal} + N \geq 0$, and $R + f(\text{\rm scal}) \, \text{\rm id} \owedge \text{\rm id} \in \text{\rm PIC2}$ for all $t \geq 0$.
\end{corollary}

The estimate $R + f(\text{\rm scal}) \, \text{\rm id} \owedge \text{\rm id} \in \text{\rm PIC2}$ ensures that blow-up limits are weakly PIC2. This can be viewed as a higher dimensional version of the fundamental Hamilton-Ivey pinching estimate in dimension $3$ (cf. \cite{Hamilton4}, \cite{Ivey}). The estimate $R - \theta \, \text{\rm scal} \, \text{\rm id} \owedge \text{\rm id} \in \text{\rm PIC}$ ensures that blow-up limits are uniformly PIC. Unlike the curvature pinching estimates in \cite{Hamilton5}, this estimate is not sharp on the cylinder.

The proof of Theorem \ref{main.theorem} will occupy Sections \ref{preliminary.pinching.estimate} -- \ref{proof.of.main.theorem}. In Section \ref{preliminary.pinching.estimate}, we construct a family of closed, convex, $O(n)$-invariant sets $\{\mathcal{G}_t^{(0)}: t \in [0,T]\}$ with the property that $\{\mathcal{G}_t^{(0)}: t \in [0,T]\}$ is invariant under the Hamilton ODE; $\mathcal{K} \subset \mathcal{G}_0^{(0)}$; and 
\begin{align*} 
\mathcal{G}_t^{(0)} 
&\subset \{R: R - \theta \, \text{\rm scal} \, \text{\rm id} \owedge \text{\rm id} \in \text{\rm PIC}\} \\ 
&\cap \{R: \text{\rm Ric}_{11}+\text{\rm Ric}_{22} - \theta \, \text{\rm scal} + N \geq 0\}. 
\end{align*}
This construction works in dimension $n \geq 5$. In Sections \ref{first.family.of.cones} and \ref{second.family.of.cones}, we construct two families of invariants cones $\mathcal{C}(b)$, $0 < b \leq b_{\text{\rm max}}$, and $\tilde{\mathcal{C}}(b)$, $0 < b \leq \tilde{b}_{\text{\rm max}}$. The family $\mathcal{C}(b)$ deforms the cone $\{R \in \text{\rm PIC}: \text{\rm Ric}_{11}+\text{\rm Ric}_{22} \geq 0\}$ inward. The family $\tilde{\mathcal{C}}(b)$ deforms the cone $\mathcal{C}(b_{\text{\rm max}}) \cap \text{\rm PIC1}$ outward. In dimension $n \geq 12$, we are able to join the two families of cones together (see Proposition \ref{pinching.cones.can.be.joined.together} below). This allows us to construct a family of sets which pinches toward PIC1. Combining this with ideas in \cite{Brendle1}, we are able to construct a family of sets which pinches toward PIC2. This is discussed in Section \ref{proof.of.main.theorem}.

In Section \ref{analysis.of.ancient.solutions}, we study ancient solutions to the Ricci flow which have bounded curvature; are $\kappa$-noncollapsed; are weakly PIC2; and are uniformly PIC. In particular, such ancient solutions satisfy a Harnack inequality and a longrange curvature estimate. Moreover, we show that such an ancient solution either splits locally as a product, or it is strictly PIC2. If the solution locally splits as a product, results in \cite{Brendle-Huisken-Sinestrari} imply that the solution is isometric to a family of shrinking cylinders $S^{n-1} \times \mathbb{R}$ or a quotient of $S^{n-1} \times \mathbb{R}$ by standard isometries. If the solution is strictly PIC2, it is either compact (in which case it is diffeomorphic to a quotient of $S^n$ by standard isometries) or noncompact (in which case it is diffeomorphic to $\mathbb{R}^n$ and has the structure of a tube with a cap attached).

In Section \ref{canonical.neighborhood.theorem}, we establish an analogue of Perelman's Canonical Neighborhood Theorem. To explain this, let $(M,g_0)$ be a compact manifold of dimension $n \geq 12$ with positive isotropic curvature, and let $g(t)$ denote the solution to the Ricci flow with initial metric $g_0$. The Canonical Neighborhood Theorem asserts that the high curvature regions in $(M,g(t))$ are modeled on ancient $\kappa$-solutions. Combining this with the results on the structure of ancient $\kappa$-solutions established in Section \ref{analysis.of.ancient.solutions}, we conclude that every point where the curvature is sufficiently large either lies on a neck, or on a cap adjacent to a neck, or on a quotient neck. Note that we may encounter singularities modeled on quotients of $S^{n-1} \times \mathbb{R}$. As in Hamilton's work \cite{Hamilton5}, these quotient necks can be ruled out if we assume that $M$ does not contain non-trivial incompressible $(n-1)$-dimensional space-forms. This makes it possible to extend the flow beyond singularities by a surgery procedure as in Perelman's work (cf. \cite{Perelman1},\cite{Perelman2},\cite{Perelman3}). Moreover, the surgically modified flow must become extinct in finite time. This is discussed in Sections \ref{first.singular.time} -- \ref{global.existence}. One simplification compared to Perelman's work is that we have an upper bound for the extinction time in terms of the infimum of the scalar curvature of the initial metric. This allows us to choose the surgery parameters independent of time $t$. 

As a corollary, we obtain a topological classification of all compact manifolds of dimension $n \geq 12$ which admit metrics of positive isotropic curvature and do not contain non-trivial incompressible $(n-1)$-dimensional space forms: 

\begin{theorem}
\label{topology}
Let $(M,g_0)$ be a compact manifold of dimension $n \geq 12$ with positive isotropic curvature. If $M$ does not contain any non-trivial incompressible $(n-1)$-dimensional space forms, then $M$ is diffeomorphic to a connected sum of finitely many spaces, each of which is a quotient of $S^n$ or $S^{n-1} \times \mathbb{R}$ by standard isometries.
\end{theorem}

Conversely, it follows from work of Micallef and Wang \cite{Micallef-Wang} that every manifold which is diffeomorphic to a connected sum of quotients of $S^n$ and $S^{n-1} \times \mathbb{R}$ admits a metric with positive isotropic curvature. Earlier results on the topology of manifolds with positive isotropic curvature (which rely on minimal surface techniques) are discussed in \cite{Fraser},\cite{Gromov},\cite{Micallef-Moore}. Theorem \ref{topology} is reminiscent of the topological classification of three-manifolds which admit metrics with positive scalar curvature (which is a consequence of Perelman's work).

Finally, we mention an interesting connection between Ricci flow on manifolds with positive isotropic curvature and mean curvature flow for two-convex hypersurfaces. If $M$ is a two-convex hypersurface in $\mathbb{R}^{n+1}$ ($n \geq 4$), then the induced metric on $M$ has positive isotropic curvature (see Lemma \ref{wedge.product}). If we evolve a two-convex hypersurface by mean curvature flow, then results of Huisken and Sinestrari \cite{Huisken-Sinestrari1},\cite{Huisken-Sinestrari2} imply that every blow-up limit is weakly convex, and that the flow only forms neck-pinch singularities. \\

\textbf{Acknowledgements.} The author would like to thank Hong Huang and Florian Johne for comments on an earlier version of this paper. He is especially grateful to an anonymous referee for many insightful remarks. This project was supported by the National Science Foundation and by the Simons Foundation. 

\section{A preliminary pinching estimate}

\label{preliminary.pinching.estimate}

Throughout this section, we assume that $n \geq 5$. Let $\mathcal{K}$ be a compact set of algebraic curvature tensors in dimension $n$ which is contained in the interior of the PIC cone. Our goal in this section is to construct a set of inequalities which are satisfied on the set $K$ and which are preserved under the Hamilton ODE. We first explain the intuition. Our starting point is the observation that, if $R \in \text{\rm PIC}$ and the sum of the two smallest eigenvalues of the Ricci tensor is negative, then the sum of the two smallest eigenvalues of the Ricci tensor is increasing under the Hamilton ODE. In particular, the set $\{R \in \text{\rm PIC}: \text{\rm Ric}_{11}+\text{\rm Ric}_{22} \geq -2\}$ is preserved under the Hamilton ODE.

Our goal is to deform the set $\{R \in \text{\rm PIC}: \text{\rm Ric}_{11}+\text{\rm Ric}_{22} \geq -2\}$ inward in such a way that the deformed sets are still preserved under the Hamilton ODE. To that end, we use the maps $\ell_{a,b}$ introduced in \cite{Bohm-Wilking}. Following \cite{Bohm-Wilking}, we define 
\[\ell_{a,b}(S) = S + b \, \text{\rm Ric}(S) \owedge \text{\rm id} + \frac{1}{n} \, (a-b) \, \text{\rm scal}(S) \, \text{\rm id} \owedge \text{\rm id}\] 
for every algebraic curvature tensor $S$. Under the map $\ell_{a,b}$, the scalar curvature changes by a factor of $1+2(n-1)a$, the tracefree Ricci tensor changes by a factor of $1+(n-2)b$, and the Weyl tensor is unchanged. It is shown in \cite{Bohm-Wilking} that 
\[\ell_{a,b}^{-1}(Q(\ell_{a,b}(S)) = Q(S) + D_{a,b}(S),\] 
where $D_{a,b}(S)$ is defined by 
\begin{align*} 
D_{a,b}(S) &= (2b+(n-2)b^2-2a) \, \tracefreeRic(S) \owedge \tracefreeRic(S) \\ 
&+ 2a \, \text{\rm Ric}(S) \owedge \text{\rm Ric}(S) + 2b^2 \, \tracefreeRic(S)^2 \owedge \text{\rm id} \\ 
&+ \frac{nb^2(1-2b)-2(a-b)(1-2b+nb^2)}{n(1+2(n-1)a)} \, |\tracefreeRic(S)|^2 \, \text{\rm id} \owedge \text{\rm id}.
\end{align*} 
The Ricci tensor of $D_{a,b}(S)$ is given by 
\begin{align*} 
\text{\rm Ric}(D_{a,b}(S)) 
&= -4b \, \text{\rm Ric}(S)^2 + \frac{4}{n} \, (2b+(n-2)a) \, \text{\rm scal}(S) \, \text{\rm Ric}(S) \\ 
&+ 2 \, \frac{n^2b^2 - 2(n-1)(a-b)(1-2b)}{n(1+2(n-1)a)} \, |\tracefreeRic(S)|^2 \, \text{\rm id} \\ 
&+ \frac{4}{n^2} \, (a-b) \, \text{\rm scal}(S)^2 \, \text{\rm id}.
\end{align*}
Hence, if $R$ evolves by the Hamilton ODE $\frac{d}{dt} R = Q(R)$, then $S = \ell_{a,b}^{-1}(R)$ evolves by the ODE $\frac{d}{dt} S = Q(S)+D_{a,b}(S)$. As in \cite{Bohm-Wilking}, it is convenient to consider combinations of small positive numbers $a,b$ such that $2a = 2b+(n-2)b^2$. In this case, the evolution equations for $S$ and $\text{\rm Ric}(S)$ take the form 
\[\frac{d}{dt} S = Q(S) + 2a \, \text{\rm Ric}(S) \owedge \text{\rm Ric}(S) + \text{\rm positive terms}\] 
and 
\begin{align*} 
\frac{d}{dt} \text{\rm Ric}(S) 
&= 2 \, S * \text{\rm Ric}(S) - 4b \, \text{\rm Ric}(S)^2 + \frac{4}{n} \, (2b+(n-2)a) \, \text{\rm scal}(S) \, \text{\rm Ric}(S) \\ 
&+ \text{\rm positive terms}, 
\end{align*}
where we have used the notation $(S * H)_{ik} := \sum_{p,q=1}^n S_{ipkq} H_{pq}$. In order to show that $S$ remains in the PIC cone, we need a lower bound for the sum of the two smallest eigenvalues of $\text{\rm Ric}(S)$. However, when we try to show that a lower bound for the sum of the smallest eigenvalues of $\text{\rm Ric}(S)$ is preserved, we encounter a problem, in that the evolution equation for $\text{\rm Ric}(S)_{11}+\text{\rm Ric}(S)_{22}$ contains a term of the form $-4b \, ((\text{\rm Ric}(S)^2)_{11}+(\text{\rm Ric}(S)^2)_{22})$, which has an unfavorable sign. To overcome this problem, we impose an additional inequality, which allows us to control the difference $\text{\rm Ric}(S)_{22}-\text{\rm Ric}(S)_{11}$ in terms of the sum $\sum_{p=3}^n (R_{1p1p}+R_{2p2p})$. This is the crucial ingredient needed to preserve the lower bound for $\text{\rm Ric}(S)_{11}+\text{\rm Ric}(S)_{22}$. The price to pay is that we need to verify that this additional inequality is itself preserved.

After these preparations, we now state the main result of this section:

\begin{theorem}
\label{rough.pinching}
Let $\mathcal{K}$ be a compact set of algebraic curvature tensors in dimension $n$ which is contained in the interior of the PIC cone, and let $T>0$ be given. Then there exist a small positive real number $\theta$, a large positive real number $N$ and a continuous family of closed, convex, $O(n)$-invariant sets $\{\mathcal{G}_t^{(0)}: t \in [0,T]\}$, such that the family $\{\mathcal{G}_t^{(0)}: t \in [0,T]\}$ is invariant under the Hamilton ODE $\frac{d}{dt} R = Q(R)$; $\mathcal{K} \subset \mathcal{G}_0^{(0)}$; and 
\begin{align*}
\mathcal{G}_t^{(0)} 
&\subset \{R: R - \theta \, \text{\rm scal} \, \text{\rm id} \owedge \text{\rm id} \in \text{\rm PIC}\} \\ 
&\cap \{R: \text{\rm Ric}_{11}+\text{\rm Ric}_{22} - \theta \, \text{\rm scal} + N \geq 0\}
\end{align*} 
for all $t \in [0,T]$.
\end{theorem}

In the remainder of this section, we give the proof of Theorem \ref{rough.pinching}. Without loss of generality, we may assume that $\text{\rm Ric}_{11}+\text{\rm Ric}_{22} \geq -1$ for all $R \in \mathcal{K}$. We first give the definition of the sets $\mathcal{G}_t^{(0)}$. The definition depends on a parameter $\delta>0$, which we choose small enough.

\begin{definition}
For $\delta>0$ small, we denote by $\mathcal{G}_t^{(0)}$ the set of all algebraic curvature tensors $R$ satisfying the following conditions: 
\begin{itemize}
\item[(i)] If $0 \leq b \leq \delta e^{-8t}$, $2a = 2b+(n-2)b^2$, and $S = \ell_{a,b}^{-1}(R)$, then $S \in \text{\rm PIC}$.
\item[(ii)] If $0 \leq b \leq \delta e^{-8t}$, $2a=2b+(n-2)b^2$, and $S = \ell_{a,b}^{-1}(R)$, then $\text{\rm Ric}(S)_{11}+\text{\rm Ric}(S)_{22} + 2b^{\frac{5}{4}} \, \text{\rm scal}(S) \geq -2$.
\item[(iii)] $R - 4\delta \, \text{\rm id} \owedge \text{\rm id} \in \text{\rm PIC}$.
\item[(iv)] For every orthonormal frame $\{e_1,\hdots,e_n\}$, the inequality  
\begin{align*} 
&R_{1313}+R_{1414}+R_{2323}+R_{2424}-2R_{1234} \\ 
&\geq \delta^{\frac{1}{4}} \, \text{\rm scal}(R)^{-3} \, \bigg ( \sum_{p,q=1}^n (R_{13pq}-R_{24pq})^2 + \sum_{p,q=1}^n (R_{14pq}+R_{23pq})^2 \bigg )^2
\end{align*} 
holds.
\end{itemize}
\end{definition}

It is easy to see that $\mathcal{K} \subset \mathcal{G}_0^{(0)}$ if $\delta>0$ is sufficiently small. 

\begin{lemma}
The set $\mathcal{G}_t^{(0)}$ is convex for each $t$.
\end{lemma}

\textbf{Proof.} 
It is clear that the inequalities (i),(ii),(iii) define a convex set. The inequality (iv) can be rewritten as 
\begin{align*} 
&\text{\rm scal}^{\frac{3}{4}} \, (R_{1313}+R_{1414}+R_{2323}+R_{2424}-2R_{1234})^{\frac{1}{4}} \\ 
&\geq \delta^{\frac{1}{16}} \, \bigg ( \sum_{p,q=1}^n (R_{13pq}-R_{24pq})^2 + \sum_{p,q=1}^n (R_{14pq}+R_{23pq})^2 \bigg )^{\frac{1}{2}}.
\end{align*} 
For each orthonormal frame $\{e_1,\hdots,e_n\}$, the function 
\[R \mapsto \text{\rm scal}^{\frac{3}{4}} \, (R_{1313}+R_{1414}+R_{2323}+R_{2424}-2R_{1234})^{\frac{1}{4}}\] 
is concave, while the function 
\[R \mapsto \bigg ( \sum_{p,q=1}^n (R_{13pq}-R_{24pq})^2 + \sum_{p,q=1}^n (R_{14pq}+R_{23pq})^2 \bigg )^{\frac{1}{2}}\] 
is convex. This shows that the inequality (iv) defines a convex set. \\

We now verify that the family $\{\mathcal{G}_t^{(0)}: t \in [0,T]\}$ is invariant under the Hamilton ODE $\frac{d}{dt} R = Q(R)$ if $\delta>0$ is sufficiently small. Let $R(t)$ be a solution of the ODE $\frac{d}{dt} R(t) = Q(R(t))$. Moreover, suppose that $R(t_0)$ lies in the interior of $\mathcal{G}_{t_0}^{(0)}$ for some time $t_0 \in [0,T]$. We will show that $R(t) \in \mathcal{G}_t^{(0)}$ for all $t \in [t_0,T]$, provided that $\delta>0$ is sufficiently small. To prove this, we argue by contradiction. Suppose that $R(t) \notin \mathcal{G}_t^{(0)}$ for some $t \in [t_0,T]$, and let 
\[\hat{t} := \inf \{t \in [t_0,T]: R(t) \notin \mathcal{G}_t^{(0)}\}.\] 
Clearly, $\hat{t} \in (t_0,T]$, $R(t) \in \mathcal{G}_t^{(0)}$ for all $t \in [0,\hat{t}]$, and $\hat{R} := R(\hat{t}) \in \partial \mathcal{G}_{\hat{t}}^{(0)}$. Our strategy is to show that, at time $\hat{t}$, the conditions (i) -- (iv) hold with strict inequalities. We begin with several technical lemmata:

\begin{lemma}
\label{technical.0}
For every orthonormal frame $\{e_1,\hdots,e_n\}$, the inequality  
\begin{align*} 
&\hat{R}_{1313}+\hat{R}_{1414}+\hat{R}_{2323}+\hat{R}_{2424} \\ 
&\geq \delta^{\frac{1}{4}} \, \text{\rm scal}(\hat{R})^{-3} \, \bigg ( \sum_{p,q=1}^n (R_{13pq}^2+R_{14pq}^2+R_{23pq}^2+R_{24pq}^2) \bigg )^2
\end{align*} 
holds.
\end{lemma}

\textbf{Proof.} 
The condition (iv) gives a lower bound for $\hat{R}_{1313}+\hat{R}_{1414}+\hat{R}_{2323}+\hat{R}_{2424}-2\hat{R}_{1234}$ and $\hat{R}_{1313}+\hat{R}_{1414}+\hat{R}_{2323}+\hat{R}_{2424}+2\hat{R}_{1234}$. If we add these inequalities, the assertion follows. \\

\begin{lemma}
\label{technical.1}
Suppose that $0 < \hat{b} \leq \delta e^{-8\hat{t}}$, $2\hat{a}=2\hat{b}+(n-2)\hat{b}^2$, and $\hat{S} := \ell_{\hat{a},\hat{b}}^{-1}(\hat{R})$. If we define 
\begin{align*} 
U &:= \frac{1}{1+(n-2)\hat{b}} \, \text{\rm Ric}(\hat{S}) \owedge \text{\rm id} \\ 
&- \frac{1}{n} \, \Big ( \frac{1}{1+(n-2)\hat{b}} - \frac{1+(n-2)\hat{b}}{1+2(n-1)\hat{a}} \Big ) \, \text{\rm scal}(\hat{S}) \, \text{\rm id} \owedge \text{\rm id}, 
\end{align*} 
then $U \in T_{\hat{S}} \text{\rm PIC}$. 
\end{lemma} 

\textbf{Proof.} 
Since condition (i) still holds at time $\hat{t}$, we know that $\ell_{b+(n-2)b^2/2,b}^{-1}(\hat{R}) \in \text{\rm PIC}$ for all $0 < b \leq \hat{b}$. This implies 
\[U = -\frac{d}{db} \ell_{b+(n-2)b^2/2,b}^{-1}(\hat{R}) \Big |_{b=\hat{b}} \in T_{\hat{S}} \text{\rm PIC},\] 
as claimed. \\

\begin{lemma} 
\label{technical.2}
Suppose that $0 < \hat{b} \leq \delta e^{-8\hat{t}}$, $2\hat{a}=2\hat{b}+(n-2)\hat{b}^2$, $\hat{S} := \ell_{\hat{a},\hat{b}}^{-1}(\hat{R})$, 
and 
\begin{align*} 
U &:= \frac{1}{1+(n-2)\hat{b}} \, \text{\rm Ric}(\hat{S}) \owedge \text{\rm id} \\ 
&- \frac{1}{n} \, \Big ( \frac{1}{1+(n-2)\hat{b}} - \frac{1+(n-2)\hat{b}}{1+2(n-1)\hat{a}} \Big ) \, \text{\rm scal}(\hat{S}) \, \text{\rm id} \owedge \text{\rm id}. 
\end{align*} 
Then 
\[8\hat{b} \, U - 4\hat{a} \, \text{\rm Ric}(\hat{S}) \owedge \text{\rm id} - 4\hat{a} \, \text{\rm id} \owedge \text{\rm id} + \hat{b}^{\frac{17}{8}} \, \text{\rm scal}(\hat{S})^2 \, \text{\rm id} \owedge \text{\rm id} \in T_{\hat{S}} \text{\rm PIC}\] 
if $\delta>0$ is sufficiently small.
\end{lemma}

\textbf{Proof.} 
Suppose that $\{e_1,e_2,e_3,e_4\}$ is an orthonormal four-frame satisfying 
\[\hat{S}_{1313}+\hat{S}_{1414}+\hat{S}_{2323}+\hat{S}_{2424}-2\hat{S}_{1234} = 0.\] 
We need to show that 
\begin{align*} 
&\hat{b} \, (U_{1313}+U_{1414}+U_{2323}+U_{2424}-2U_{1234}) \\ 
&- \hat{a} \, (\text{\rm Ric}(\hat{S})_{11}+\text{\rm Ric}(\hat{S})_{22}+\text{\rm Ric}(\hat{S})_{33}+\text{\rm Ric}(\hat{S})_{44}) - 4\hat{a} + \hat{b}^{\frac{17}{8}} \, \text{\rm scal}(\hat{S})^2 \geq 0 
\end{align*} 
for this particular orthonormal four-frame $\{e_1,e_2,e_3,e_4\}$. 

\textit{Case 1:} Suppose first that $\text{\rm Ric}(\hat{S})_{11}+\text{\rm Ric}(\hat{S})_{22}+\text{\rm Ric}(\hat{S})_{33}+\text{\rm Ric}(\hat{S})_{44} \leq 8$. Using condition (iii), we obtain $\hat{R} - 4\hat{b} \, \text{\rm id} \owedge \text{\rm id} \in \text{\rm PIC}$. Therefore, 
\begin{align*} 
32\hat{b} 
&\leq \hat{R}_{1313}+\hat{R}_{1414}+\hat{R}_{2323}+\hat{R}_{2424}-2\hat{R}_{1234} \\ 
&= \hat{S}_{1313}+\hat{S}_{1414}+\hat{S}_{2323}+\hat{S}_{2424}-2\hat{S}_{1234} \\ 
&+ 2\hat{b} \, (\text{\rm Ric}(\hat{S})_{11}+\text{\rm Ric}(\hat{S})_{22}+\text{\rm Ric}(\hat{S})_{33}+\text{\rm Ric}(\hat{S})_{44}) + \frac{4(n-2)\hat{b}^2}{n} \, \text{\rm scal}(\hat{S}) \\ 
&\leq 16\hat{b} + \frac{4(n-2)\hat{b}^2}{n} \, \text{\rm scal}(\hat{S}),
\end{align*} 
hence $\hat{b} \, \text{\rm scal}(\hat{S}) \geq \frac{4n}{n-2}$. Since $U \in T_{\hat{S}} \text{\rm PIC}$ by Lemma \ref{technical.1}, we conclude that  
\begin{align*} 
&\hat{b} \, (U_{1313}+U_{1414}+U_{2323}+U_{2424}-2U_{1234}) \\ 
&- \hat{a} \, (\text{\rm Ric}(\hat{S})_{11}+\text{\rm Ric}(\hat{S})_{22}+\text{\rm Ric}(\hat{S})_{33}+\text{\rm Ric}(\hat{S})_{44}) - 4\hat{a} + \hat{b}^{\frac{17}{8}} \, \text{\rm scal}(\hat{S})^2 \\ 
&\geq -\hat{a} \, (\text{\rm Ric}(\hat{S})_{11}+\text{\rm Ric}(\hat{S})_{22}+\text{\rm Ric}(\hat{S})_{33}+\text{\rm Ric}(\hat{S})_{44}) - 4\hat{a} + \hat{b}^{\frac{17}{8}} \, \text{\rm scal}(\hat{S})^2 \\ 
&\geq -12\hat{a} + \frac{16n^2}{(n-2)^2} \, \hat{b}^{\frac{1}{8}} \\ 
&> 0
\end{align*} 
if $\hat{b}>0$ is sufficiently small.

\textit{Case 2:} Suppose finally that $\text{\rm Ric}(\hat{S})_{11}+\text{\rm Ric}(\hat{S})_{22}+\text{\rm Ric}(\hat{S})_{33}+\text{\rm Ric}(\hat{S})_{44} \geq 8$. By definition of $U$, we obtain 
\begin{align*} 
&\hat{b} \, (U_{1313}+U_{1414}+U_{2323}+U_{2424}-2U_{1234}) \\ 
&- \hat{a} \, (\text{\rm Ric}(\hat{S})_{11}+\text{\rm Ric}(\hat{S})_{22}+\text{\rm Ric}(\hat{S})_{33}+\text{\rm Ric}(\hat{S})_{44}) - 4\hat{a} + \hat{b}^{\frac{17}{8}} \, \text{\rm scal}(\hat{S})^2 \\ 
&= \Big ( \frac{2\hat{b}}{1+(n-2)\hat{b}} - \hat{a} \Big ) \, (\text{\rm Ric}(\hat{S})_{11}+\text{\rm Ric}(\hat{S})_{22}+\text{\rm Ric}(\hat{S})_{33}+\text{\rm Ric}(\hat{S})_{44}) \\ 
&- \frac{8\hat{b}}{n} \, \Big ( \frac{1}{1+(n-2)\hat{b}} - \frac{1+(n-2)\hat{b}}{1+2(n-1)\hat{a}} \Big ) \, \text{\rm scal}(\hat{S}) - 4\hat{a} + \hat{b}^{\frac{17}{8}} \, \text{\rm scal}(\hat{S})^2 \\ 
&\geq 8 \, \Big ( \frac{2\hat{b}}{1+(n-2)\hat{b}} - \hat{a} \Big ) - 4\hat{a} + \hat{b}^{\frac{17}{8}} \, \text{\rm scal}(\hat{S})^2 \\ 
&- \frac{8\hat{b}}{n} \, \Big ( \frac{1}{1+(n-2)\hat{b}} - \frac{1+(n-2)\hat{b}}{1+2(n-1)\hat{a}} \Big ) \, \text{\rm scal}(\hat{S}) \\ 
&\geq \hat{b} + \hat{b}^{\frac{17}{8}} \, \text{\rm scal}(\hat{S})^2 - C(n) \, \hat{b}^2 \, \text{\rm scal}(\hat{S}) \\ 
&> 0 
\end{align*} 
if $\hat{b}>0$ is sufficiently small. \\

\begin{lemma} 
\label{i.preserved}
Let $\hat{a},\hat{b}$ be real numbers such that $0 \leq \hat{b} \leq \delta e^{-8\hat{t}}$, $2\hat{a} = 2\hat{b}+(n-2)\hat{b}^2$, and let $\hat{S} := \ell_{a,b}^{-1}(\hat{R})$. Then $\hat{S}$ lies in the interior of the PIC cone.
\end{lemma} 

\textbf{Proof.} 
In view of condition (iii), $\hat{R}$ lies in the interior of the PIC cone. Therefore, the assertion is true for $\hat{b}=0$. Hence, it suffices to consider the case $0 < \hat{b} \leq \delta e^{-8\hat{t}}$. Let us define functions $a(t),b(t)$ by $b(t) = e^{8(\hat{t}-t)} \hat{b}$ and $a(t) = 2b(t)+(n-2)b(t)^2$. Moreover, let $S(t) := \ell_{a(t),b(t)}^{-1}(R(t))$. Since condition (i) holds up to time $\hat{t}$, we know that $S(t) \in \text{\rm PIC}$ for all $t \in [0,\hat{t}]$. The evolution of $S(t)$ is given by 
\[\frac{d}{dt} S(t) \Big |_{t=\hat{t}} = Q(\hat{S}) + D_{\hat{a},\hat{b}}(\hat{S}) - b'(\hat{t}) \, U = Q(\hat{S}) + D_{\hat{a},\hat{b}}(\hat{S}) + 8\hat{b} \, U,\] 
where 
\begin{align*} 
U &= \frac{1}{1+(n-2)\hat{b}} \, \text{\rm Ric}(\hat{S}) \owedge \text{\rm id} \\ 
&- \frac{1}{n} \, \Big ( \frac{1}{1+(n-2)\hat{b}} - \frac{1+(n-2)\hat{b}}{1+2(n-1)\hat{a}} \Big ) \, \text{\rm scal}(\hat{S}) \, \text{\rm id} \owedge \text{\rm id}. 
\end{align*} 
We claim that $D_{\hat{a},\hat{b}}(\hat{S}) + 8\hat{b} \, U$ lies in the interior of the tangent cone $T_{\hat{S}} \text{\rm PIC}$. To verify this, we distinguish two cases: 

\textit{Case 1:} Suppose first that $\text{\rm Ric}(\hat{S})$ is strictly two-positive. By Lemma \ref{wedge.product}, $\text{\rm Ric}(\hat{S}) \owedge \text{\rm Ric}(\hat{S})$ lies in the interior of the PIC cone. Since $2\hat{a} = 2\hat{b}+(n-2)\hat{b}^2$, this implies that $D_{\hat{a},\hat{b}}(\hat{S})$ lies in the interior of the PIC cone. Since $U \in T_{\hat{S}} \text{\rm PIC}$ by Lemma \ref{technical.1}, it follows that $D_{\hat{a},\hat{b}}(\hat{S}) + 8\hat{b} \, U$ lies in the interior of the tangent cone $T_{\hat{S}} \text{\rm PIC}$, as claimed. 

\textit{Case 2:} Suppose next that $\text{\rm Ric}(\hat{S})$ is not strictly two-positive. In this case, $\text{\rm scal}(\hat{S}) \leq C(n) \, |\tracefreeRic(\hat{S})|$. Since $2\hat{a}=2\hat{b}+(n-2)\hat{b}^2$, it follows that 
\[D_{\hat{a},\hat{b}}(\hat{S}) - 2\hat{a} \, \text{\rm Ric}(\hat{S}) \owedge \text{\rm Ric}(\hat{S}) - 2\hat{b}^{\frac{17}{8}} \, \text{\rm scal}(\hat{S})^2 \, \text{\rm id} \owedge \text{\rm id} \in \text{\rm PIC}\] 
if $\hat{b}>0$ is sufficiently small. The condition (ii) implies that $\text{\rm Ric}(\hat{S}) + (1+\hat{b}^{\frac{5}{4}} \, \text{\rm scal}(\hat{S})) \, \text{\rm id}$ is weakly two-positive. Using Lemma \ref{wedge.product}, we obtain 
\begin{align*} 
&2\hat{a} \, \text{\rm Ric}(\hat{S}) \owedge \text{\rm Ric}(\hat{S}) + 4\hat{a}(1+\hat{b}^{\frac{5}{4}} \, \text{\rm scal}(\hat{S})) \, \text{\rm Ric}(\hat{S}) \owedge \text{\rm id} + 4\hat{a} (1+\hat{b}^{\frac{5}{2}} \, \text{\rm scal}(\hat{S})^2) \, \text{\rm id} \owedge \text{\rm id} \\ 
&= 2\hat{a} \, [\text{\rm Ric}(\hat{S}) + (1+\hat{b}^{\frac{5}{4}} \, \text{\rm scal}(\hat{S})) \, \text{\rm id}] \owedge [\text{\rm Ric}(\hat{S}) + (1+\hat{b}^{\frac{5}{4}} \, \text{\rm scal}(\hat{S})) \, \text{\rm id}] \\ 
&+ 2\hat{a} (1-\hat{b}^{\frac{5}{4}} \, \text{\rm scal}(\hat{S}))^2 \, \text{\rm id} \owedge \text{\rm id} \in \text{\rm PIC}.
\end{align*} 
On the other hand, Lemma \ref{technical.2} gives 
\[8\hat{b} \, U - 4\hat{a} \, \text{\rm Ric}(\hat{S}) \owedge \text{\rm id} - 4\hat{a} \, \text{\rm id} \owedge \text{\rm id} + \hat{b}^{\frac{17}{8}} \, \text{\rm scal}(\hat{S})^2 \, \text{\rm id} \owedge \text{\rm id} \in T_{\hat{S}} \text{\rm PIC}.\] 
Moreover, if $\hat{b}>0$ is sufficiently small, then $\hat{a}\hat{b}^{\frac{5}{4}}$ is much smaller than $\hat{b}^{\frac{17}{8}}$. Since $|\text{\rm Ric}(\hat{S})| \leq C(n) \, \text{\rm scal}(\hat{S})$, it follows that  
\[\hat{b}^{\frac{17}{8}} \, \text{\rm scal}(S)^2 \, \text{\rm id} \owedge \text{\rm id} - 4\hat{a}\hat{b}^{\frac{5}{4}} \, \text{\rm scal}(\hat{S}) \, \text{\rm Ric}(\hat{S}) \owedge \text{\rm id} - 4\hat{a}\hat{b}^{\frac{5}{2}} \, \text{\rm scal}(\hat{S})^2 \, \text{\rm id} \owedge \text{\rm id}\] 
lies in the interior of the PIC cone. Adding all four formulae, we conclude that $D_{\hat{a},\hat{b}}(\hat{S}) + 8\hat{b} \, U$ lies in the interior of the tangent cone $T_{\hat{S}} \text{\rm PIC}$. 

To summarize, we have shown that $D_{\hat{a},\hat{b}}(\hat{S}) + 8\hat{b} \, U$ lies in the interior of the tangent cone $T_{\hat{S}} \text{\rm PIC}$. Since $Q(\hat{S}) \in T_{\hat{S}} \text{\rm PIC}$ by Proposition 7.5 in \cite{Brendle-book}, we conclude that $\frac{d}{dt} S(t) \big |_{t=\hat{t}}$ lies in the interior of the tangent cone $T_{\hat{S}} \text{\rm PIC}$. Hence, if $\hat{S}$ lies on the boundary of the PIC cone, then $S(t)$ lies outside the PIC cone if $t \in [0,\hat{t})$ is sufficiently close to $\hat{t}$. This is a contradiction. This completes the proof of Lemma \ref{i.preserved}. \\

\begin{lemma}
\label{ii.preserved}
Let $a,b$ be real numbers such that $0 \leq b \leq \delta e^{-8\hat{t}}$, $2a = 2b+(n-2)b^2$, and let $\hat{S} := \ell_{a,b}^{-1}(\hat{R})$. Then $\text{\rm Ric}(\hat{S})_{11}+\text{\rm Ric}(\hat{S})_{22} + 2b^{\frac{5}{4}} \, \text{\rm scal}(\hat{S}) > -2$ for every pair of orthonormal vectors $\{e_1,e_2\}$. 
\end{lemma}

\textbf{Proof.} 
Suppose that the assertion is false. Then 
\[\text{\rm Ric}(\hat{S})_{11}+\text{\rm Ric}(\hat{S})_{22} + 2b^{\frac{5}{4}} \, \text{\rm scal}(\hat{S}) = -2\] 
for some pair of orthonormal vectors $\{e_1,e_2\}$. We may assume that $\{e_1,e_2\}$ are eigenvectors of $\text{\rm Ric}(\hat{S})$, and $\text{\rm Ric}(\hat{S})_{22} \geq \text{\rm Ric}(\hat{S})_{11}$. Let us extend $\{e_1,e_2\}$ to an eigenbasis $\{e_1,\hdots,e_n\}$ of $\text{\rm Ric}(\hat{S})$. Clearly, $\text{\rm Ric}(\hat{S})_{pp} \geq \text{\rm Ric}(\hat{S})_{22}$ for $3 \leq p \leq n$. For each $t \in [0,\hat{t}]$, we define $S(t) := l_{a,b}^{-1}(R(t))$. Since the condition (ii) holds up to time $\hat{t}$, we know that  
\[\text{\rm Ric}(S(t))_{11}+\text{\rm Ric}(S(t))_{22} + 2b^{\frac{5}{4}} \, \text{\rm scal}(S(t)) \geq -2\] 
for all $t \in [0,\hat{t}]$. The evolution of $S(t)$ is given by $\frac{d}{dt} S(t) \big |_{t=\hat{t}} = Q(\hat{S}) + D_{a,b}(\hat{S})$. Using the formula for $\text{\rm Ric}(D_{a,b}(\hat{S}))$, we obtain 
\begin{align*} 
0 &\geq \frac{1}{2} \, \frac{d}{dt} \big ( \text{\rm Ric}(S(t))_{11}+\text{\rm Ric}(S(t))_{22} + 2b^{\frac{5}{4}} \, \text{\rm scal}(S(t)) \big ) \Big |_{t=\hat{t}} \\ 
&\geq \sum_{p=1}^n (\hat{S}_{1p1p}+\hat{S}_{2p2p}) \text{\rm Ric}(\hat{S})_{pp} \\ 
&+ \frac{3}{2} \, b^{\frac{5}{4}} \, |\text{\rm Ric}(\hat{S})|^2 - 2b \, ((\text{\rm Ric}(\hat{S})^2)_{11}+(\text{\rm Ric}(\hat{S})^2)_{22}) \\ 
&+ \frac{4}{n} \, (2b+(n-2)a) \, \text{\rm scal}(\hat{S}) \, (\text{\rm Ric}(\hat{S})_{11}+\text{\rm Ric}(\hat{S})_{22}) \\ 
&= \sum_{p=3}^n (\hat{S}_{1p1p}+\hat{S}_{2p2p}) (\text{\rm Ric}(\hat{S})_{pp} - \frac{1}{2} \, (\text{\rm Ric}(\hat{S})_{11}+\text{\rm Ric}(\hat{S})_{22})) \\ 
&+ \frac{3}{2} \, b^{\frac{5}{4}} \, |\text{\rm Ric}(\hat{S})|^2 + \frac{1}{2} \, (1-2b) \, (\text{\rm Ric}(\hat{S})_{11}+\text{\rm Ric}(\hat{S})_{22})^2 - b \, (\text{\rm Ric}(\hat{S})_{22}-\text{\rm Ric}(\hat{S})_{11})^2 \\ 
&+ \frac{4}{n} \, (2b+(n-2)a) \, \text{\rm scal}(\hat{S}) \, (\text{\rm Ric}(\hat{S})_{11}+\text{\rm Ric}(\hat{S})_{22}).
\end{align*} 
if $b \geq 0$ is sufficiently small. Note that $|\text{\rm Ric}(\hat{S})|^2 \geq \frac{1}{n} \, \text{\rm scal}(\hat{S})^2$. Hence, if $b \geq 0$ is sufficiently small, then 
\begin{align*} 
&\frac{1}{2} \, b^{\frac{5}{4}} \, |\text{\rm Ric}(\hat{S})|^2 + \frac{1}{2} \, (1-2b) \, (\text{\rm Ric}(\hat{S})_{11}+\text{\rm Ric}(\hat{S})_{22})^2 \\ 
&> -\frac{4}{n} \, (2b+(n-2)a) \, \text{\rm scal}(\hat{S}) \, (\text{\rm Ric}(\hat{S})_{11}+\text{\rm Ric}(\hat{S})_{22}) 
\end{align*}
by Young's inequality. This gives 
\begin{align*} 
0 &\geq \frac{1}{2} \, \frac{d}{dt} \big ( \text{\rm Ric}(S(t))_{11}+\text{\rm Ric}(S(t))_{22} + 2b^{\frac{5}{4}} \, \text{\rm scal}(S(t)) \big ) \Big |_{t=\hat{t}} \\ 
&> \sum_{p=3}^n (\hat{S}_{1p1p}+\hat{S}_{2p2p}) (\text{\rm Ric}(\hat{S})_{pp} - \frac{1}{2} \, (\text{\rm Ric}(\hat{S})_{11}+\text{\rm Ric}(\hat{S})_{22})) \\ 
&+ b^{\frac{5}{4}} \, |\text{\rm Ric}(\hat{S})|^2 - b \, (\text{\rm Ric}(\hat{S})_{22}-\text{\rm Ric}(\hat{S})_{11})^2. 
\end{align*} 
At this point, we distinguish two cases: 

\textit{Case 1:} Suppose first that $\hat{S}_{1p1p}+\hat{S}_{2p2p} \geq 0$ for all $p \in \{3,\hdots,n\}$. In this case, 
\begin{align*} 
&\sum_{p=3}^n (\hat{S}_{1p1p}+\hat{S}_{2p2p}) (\text{\rm Ric}(\hat{S})_{pp} - \frac{1}{2} \, (\text{\rm Ric}(\hat{S})_{11}+\text{\rm Ric}(\hat{S})_{22})) \\ 
&\geq \frac{1}{2} \sum_{p=3}^n (\hat{S}_{1p1p}+\hat{S}_{2p2p}) \, (\text{\rm Ric}(\hat{S})_{22}-\text{\rm Ric}(\hat{S})_{11}),
\end{align*}
hence 
\begin{align*} 
0 &\geq \frac{1}{2} \, \frac{d}{dt} \big ( \text{\rm Ric}(S(t))_{11}+\text{\rm Ric}(S(t))_{22} + 2b^{\frac{5}{4}} \, \text{\rm scal}(S(t)) \big ) \Big |_{t=\hat{t}} \\ 
&> \frac{1}{2} \sum_{p=3}^n (\hat{S}_{1p1p}+\hat{S}_{2p2p}) \, (\text{\rm Ric}(\hat{S})_{22}-\text{\rm Ric}(\hat{S})_{11}) \\ 
&+ b^{\frac{5}{4}} \, |\text{\rm Ric}(\hat{S})|^2 - b \, (\text{\rm Ric}(\hat{S})_{22}-\text{\rm Ric}(\hat{S})_{11})^2.  
\end{align*} 
If $b=0$ or $\text{\rm Ric}(\hat{S})_{22}-\text{\rm Ric}(\hat{S})_{11} \leq b^{\frac{1}{8}} \, |\text{\rm Ric}(\hat{S})|$ or $\sum_{p=3}^n (\hat{S}_{1p1p}+\hat{S}_{2p2p}) \geq 2b \, (\text{\rm Ric}(\hat{S})_{22}-\text{\rm Ric}(\hat{S})_{11})$, then the right hand side is nonnegative and we arrive at a contradiction. Therefore, we must have $b>0$ and $\text{\rm Ric}(\hat{S})_{22}-\text{\rm Ric}(\hat{S})_{11} \geq b^{\frac{1}{8}} \, |\text{\rm Ric}(\hat{S})|$ and $\sum_{p=3}^n (\hat{S}_{1p1p}+\hat{S}_{2p2p}) \leq 2b \, (\text{\rm Ric}(\hat{S})_{22}-\text{\rm Ric}(\hat{S})_{11})$. Since $|\hat{R}-\hat{S}| \leq C(n)b \, \text{\rm scal}(\hat{R})$, it follows that $\text{\rm Ric}(\hat{R})_{22}-\text{\rm Ric}(\hat{R})_{11} \geq c(n) \, b^{\frac{1}{8}} \, \text{\rm scal}(\hat{R})$ and $\sum_{p=3}^n (\hat{R}_{1p1p}+\hat{R}_{2p2p}) \leq C(n) b \, \text{\rm scal}(\hat{R})$. On the other hand, using Lemma \ref{technical.0}, we obtain 
\begin{align*} 
&\hat{R}_{1p1p}+\hat{R}_{2p2p}+\hat{R}_{1q1q}+\hat{R}_{2q2q} \\ 
&\geq c(n) \, b^{\frac{1}{4}} \, \text{\rm scal}(\hat{R})^{-3} \, (|\hat{R}_{1p1p}|+|\hat{R}_{2p2p}|+|\hat{R}_{1q1q}|+|\hat{R}_{2q2q}|)^4 
\end{align*}
for $3 \leq p < q \leq n$. Summation over $p$ and $q$ gives 
\begin{align*} 
\sum_{p=3}^n (\hat{R}_{1p1p}+\hat{R}_{2p2p}) 
&\geq c(n) \, b^{\frac{1}{4}} \, \text{\rm scal}(\hat{R})^{-3} \, \bigg ( \sum_{p=3}^n (|\hat{R}_{1p1p}|+|\hat{R}_{2p2p}|) \bigg )^4 \\ 
&\geq c(n) \, b^{\frac{1}{4}} \, \text{\rm scal}(\hat{R})^{-3} \, (\text{\rm Ric}(\hat{R})_{22}-\text{\rm Ric}(\hat{R})_{11})^4 \\ 
&\geq c(n) \, b^{\frac{3}{4}} \, \text{\rm scal}(\hat{R}). 
\end{align*} 
This contradicts the inequality $\sum_{p=3}^n (\hat{R}_{1p1p}+\hat{R}_{2p2p}) \leq C(n) b \, \text{\rm scal}(\hat{R})$ if $b$ is sufficiently small.

\textit{Case 2:} Suppose next that $\hat{S}_{1m1m}+\hat{S}_{2m2m} < 0$ for some $m \in \{3,\hdots,n\}$. Since $\hat{S} \in \text{\rm PIC}$, Lemma \ref{estimate.for.largest.ricci.eigenvalue} implies that each eigenvalue of $\text{\rm Ric}(\hat{S})$ is bounded from above by the sum of all the other eigenvalues. Therefore, 
\[\sum_{p \in \{1,\hdots,n\} \setminus \{m\}} \text{\rm Ric}(\hat{S})_{pp} \geq \text{\rm Ric}(\hat{S})_{mm}.\] 
Since $\text{\rm Ric}(\hat{S})_{11}+\text{\rm Ric}(\hat{S})_{22} \leq 0$, it follows that 
\begin{align*} 
&\sum_{p \in \{3,\hdots,n\} \setminus \{m\}}  (\text{\rm Ric}(\hat{S})_{pp} - \frac{1}{2} \, (\text{\rm Ric}(\hat{S})_{11}+\text{\rm Ric}(\hat{S})_{22})) \\ 
&\geq \text{\rm Ric}(\hat{S})_{mm} - \frac{1}{2} \, (\text{\rm Ric}(\hat{S})_{11}+\text{\rm Ric}(\hat{S})_{22}). 
\end{align*} 
This implies 
\begin{align*} 
&\sum_{p=3}^n (\hat{S}_{1p1p}+\hat{S}_{2p2p}) (\text{\rm Ric}(\hat{S})_{pp} - \frac{1}{2} \, (\text{\rm Ric}(\hat{S})_{11}+\text{\rm Ric}(\hat{S})_{22})) \\ 
&\geq \sum_{p \in \{3,\hdots,n\} \setminus \{m\}} (\hat{S}_{1m1m}+\hat{S}_{2m2m}+\hat{S}_{1p1p}+\hat{S}_{2p2p}) \\ 
&\hspace{20mm} \cdot (\text{\rm Ric}(\hat{S})_{pp} - \frac{1}{2} \, (\text{\rm Ric}(\hat{S})_{11}+\text{\rm Ric}(\hat{S})_{22})) \\ 
&\geq \frac{1}{2} \sum_{p \in \{3,\hdots,n\} \setminus \{m\}} (\hat{S}_{1m1m}+\hat{S}_{2m2m}+\hat{S}_{1p1p}+\hat{S}_{2p2p}) \, (\text{\rm Ric}(\hat{S})_{22}-\text{\rm Ric}(\hat{S})_{11}), 
\end{align*} 
hence 
\begin{align*} 
0 &\geq \frac{1}{2} \, \frac{d}{dt} \big ( \text{\rm Ric}(S(t))_{11}+\text{\rm Ric}(S(t))_{22} + 2b^{\frac{5}{4}} \, \text{\rm scal}(S(t)) \big ) \Big |_{t=\hat{t}} \\ 
&> \frac{1}{2} \sum_{p \in \{3,\hdots,n\} \setminus \{m\}} (\hat{S}_{1m1m}+\hat{S}_{2m2m}+\hat{S}_{1p1p}+\hat{S}_{2p2p}) \, (\text{\rm Ric}(\hat{S})_{22}-\text{\rm Ric}(\hat{S})_{11}) \\ 
&+ b^{\frac{5}{4}} \, |\text{\rm Ric}(\hat{S})|^2 - b \, (\text{\rm Ric}(\hat{S})_{22}-\text{\rm Ric}(\hat{S})_{11})^2.  
\end{align*} 
If $b=0$ or $\text{\rm Ric}(\hat{S})_{22}-\text{\rm Ric}(\hat{S})_{11} \leq b^{\frac{1}{8}} \, |\text{\rm Ric}(\hat{S})|$ or $\sum_{p \in \{3,\hdots,n\} \setminus \{m\}} (\hat{S}_{1m1m}+\hat{S}_{2m2m}+\hat{S}_{1p1p}+\hat{S}_{2p2p}) \geq 2b \, (\text{\rm Ric}(\hat{S})_{22}-\text{\rm Ric}(\hat{S})_{11})$, then the right hand side is nonnegative, and we arrive at a contradiction. Therefore, we must have $b>0$ and $\text{\rm Ric}(\hat{S})_{22}-\text{\rm Ric}(\hat{S})_{11} \geq b^{\frac{1}{8}} \, |\text{\rm Ric}(\hat{S})|$ and $\sum_{p \in \{3,\hdots,n\} \setminus \{m\}} (\hat{S}_{1m1m}+\hat{S}_{2m2m}+\hat{S}_{1p1p}+\hat{S}_{2p2p}) \leq 2b \, (\text{\rm Ric}(\hat{S})_{22}-\text{\rm Ric}(\hat{S})_{11})$. Since $|\hat{R}-\hat{S}| \leq C(n) b \, \text{\rm scal}(\hat{R})$, it follows that $\text{\rm Ric}(\hat{R})_{22}-\text{\rm Ric}(\hat{R})_{11} \geq c(n) \, b^{\frac{1}{8}} \, \text{\rm scal}(\hat{R})$ and $\sum_{p \in \{3,\hdots,n\} \setminus \{m\}} (\hat{R}_{1m1m}+\hat{R}_{2m2m}+\hat{R}_{1p1p}+\hat{R}_{2p2p}) \leq C(n) b \, \text{\rm scal}(\hat{R})$. On the other hand, using Lemma \ref{technical.0}, we obtain 
\begin{align*} 
&\hat{R}_{1m1m}+\hat{R}_{2m2m}+\hat{R}_{1p1p}+\hat{R}_{2p2p} \\ 
&\geq c(n) \, b^{\frac{1}{4}} \, \text{\rm scal}(\hat{R})^{-3} \, (|\hat{R}_{1m1m}|+|\hat{R}_{2m2m}|+|\hat{R}_{1p1p}|+|\hat{R}_{2p2p}|)^4 
\end{align*}
for $p \in \{3,\hdots,n\} \setminus \{m\}$. Summation over $p \in \{3,\hdots,n\} \setminus \{m\}$ gives 
\begin{align*} 
&\sum_{p \in \{3,\hdots,n\} \setminus \{m\}} (\hat{R}_{1m1m}+\hat{R}_{2m2m}+\hat{R}_{1p1p}+\hat{R}_{2p2p}) \\ 
&\geq c(n) \, b^{\frac{1}{4}} \, \text{\rm scal}(\hat{R})^{-3} \, \bigg ( \sum_{p=3}^n (|\hat{R}_{1p1p}|+|\hat{R}_{2p2p}|) \bigg )^4 \\ 
&\geq c(n) \, b^{\frac{1}{4}} \, \text{\rm scal}(\hat{R})^{-3} \, (\text{\rm Ric}(\hat{R})_{22}-\text{\rm Ric}(\hat{R})_{11})^4 \\ 
&\geq c(n) \, b^{\frac{3}{4}} \, \text{\rm scal}(\hat{R}). 
\end{align*} 
This contradicts the inequality $\sum_{p \in \{3,\hdots,n\} \setminus \{m\}} (\hat{R}_{1m1m}+\hat{R}_{2m2m}+\hat{R}_{1p1p}+\hat{R}_{2p2p}) \leq C(n) b \, \text{\rm scal}(\hat{R})$ if $b$ is sufficiently small. \\

\begin{lemma}
\label{iii.preserved}
The tensor $\hat{R} - 4\delta \, \text{\rm id} \owedge \text{\rm id}$ lies in the interior of the PIC cone.
\end{lemma}

\textbf{Proof.} 
By assumption, $R(t_0) - 4\delta \, \text{\rm id} \owedge \text{\rm id}$ lies in the interior of the PIC cone. Using Proposition \ref{minimum.isotropic.curvature}, we conclude that $R(\hat{t}) - 4\delta \, \text{\rm id} \owedge \text{\rm id}$ lies in the interior of the PIC cone. \\

\begin{lemma} 
\label{iv.preserved}
If $\delta>0$ is sufficiently small, then 
\begin{align*} 
&\hat{R}_{1313}+\hat{R}_{1414}+\hat{R}_{2323}+\hat{R}_{2424}-2\hat{R}_{1234} \\ 
&\geq \delta^{\frac{1}{4}} \, \text{\rm scal}^{-3} \, \bigg ( \sum_{p,q=1}^n (\hat{R}_{13pq}-\hat{R}_{24pq})^2 + \sum_{p,q=1}^n (\hat{R}_{14pq}+\hat{R}_{23pq})^2 \bigg )^2 
\end{align*} 
for every orthonormal frame $\{e_1,\hdots,e_n\}$.
\end{lemma}

\textbf{Proof.} 
Suppose that
\begin{align*} 
&\hat{R}_{1313}+\hat{R}_{1414}+\hat{R}_{2323}+\hat{R}_{2424}-2\hat{R}_{1234} \\ 
&= \delta^{\frac{1}{4}} \, \text{\rm scal}(\hat{R})^{-3} \, \bigg ( \sum_{p,q=1}^n (\hat{R}_{13pq}-\hat{R}_{24pq})^2 + \sum_{p,q=1}^n (\hat{R}_{14pq}+\hat{R}_{23pq})^2 \bigg )^2 
\end{align*} 
for some orthonormal frame $\{e_1,\hdots,e_n\}$. In view of condition (iii), the left hand side of the equation is strictly positive. We now argue as in Section 7.3 of \cite{Brendle-book}. Lemma 7.8, Lemma 7.10, and Lemma 7.12 in \cite{Brendle-book} hold with error terms which are bounded by $C(n) \, \delta^{\frac{1}{4}} \, \text{\rm scal}(\hat{R})^{-1} \, \big ( \sum_{p,q=1}^n (\hat{R}_{13pq}-\hat{R}_{24pq})^2 + \sum_{p,q=1}^n (\hat{R}_{14pq}+\hat{R}_{23pq})^2 \big )$. Consequently, Lemma 7.9, Lemma 7.11, and Lemma 7.13 in \cite{Brendle-book} hold with error terms which are  bounded by $C(n) \, \delta^{\frac{1}{4}} \, \big ( \sum_{p,q=1}^n (\hat{R}_{13pq}-\hat{R}_{24pq})^2 + \sum_{p,q=1}^n (\hat{R}_{14pq}+\hat{R}_{23pq})^2 \big )$. Putting these facts together, we obtain an inequality of the form 
\begin{align*} 
&(\hat{R}^\#)_{1313}+(\hat{R}^\#)_{1414}+(\hat{R}^\#)_{2323}+(\hat{R}^\#)_{2424}+2(\hat{R}^\#)_{1342}+2(\hat{R}^\#)_{1423} \\ 
&\geq -C(n) \, \delta^{\frac{1}{4}} \, \bigg ( \sum_{p,q=1}^n (\hat{R}_{13pq}-\hat{R}_{24pq})^2 + \sum_{p,q=1}^n (\hat{R}_{14pq}+\hat{R}_{23pq})^2 \bigg ). 
\end{align*} 
This implies 
\begin{align*} 
&\frac{d}{dt} (R(t)_{1313}+R(t)_{1414}+R(t)_{2323}+R(t)_{2424}-2R(t)_{1234}) \Big |_{t=\hat{t}} \\ 
&= (\hat{R}^2)_{1313}+(\hat{R}^2)_{1414}+(\hat{R}^2)_{2323}+(\hat{R}^2)_{2424}+2(\hat{R}^2)_{1342}+2(\hat{R}^2)_{1423} \\ 
&+ (\hat{R}^\#)_{1313}+(\hat{R}^\#)_{1414}+(\hat{R}^\#)_{2323}+(\hat{R}^\#)_{2424}+2(\hat{R}^\#)_{1342}+2(\hat{R}^\#)_{1423} \\ 
&\geq (1-C(n) \, \delta^{\frac{1}{4}}) \, \bigg ( \sum_{p,q=1}^n (\hat{R}_{13pq}-\hat{R}_{24pq})^2 + \sum_{p,q=1}^n (\hat{R}_{14pq}+\hat{R}_{23pq})^2 \bigg ). 
\end{align*}
On the other hand, since the norm of $\hat{R}$ is controlled by $\text{\rm scal}(\hat{R})$, we obtain 
\begin{align*} 
&\frac{d}{dt} \bigg [ \text{\rm scal}(R(t))^{-3} \, \bigg ( \sum_{p,q=1}^n (R(t)_{13pq}-R(t)_{24pq})^2 + \sum_{p,q=1}^n (R(t)_{14pq}+R(t)_{23pq})^2 \bigg )^2 \bigg ] \bigg |_{t=\hat{t}} \\ 
&\leq C(n) \, \bigg ( \sum_{p,q=1}^n (\hat{R}_{13pq}-\hat{R}_{24pq})^2 + \sum_{p,q=1}^n (\hat{R}_{14pq}+\hat{R}_{23pq})^2 \bigg ). 
\end{align*} 
Thus, we conclude that 
\begin{align*} 
&\frac{d}{dt} \bigg [ R(t)_{1313}+R(t)_{1414}+R(t)_{2323}+R(t)_{2424}-2R(t)_{1234} \\ 
&- \delta^{\frac{1}{4}} \, \text{\rm scal}(R(t))^{-3} \, \bigg ( \sum_{p,q=1}^n (R(t)_{13pq}-R(t)_{24pq})^2 + \sum_{p,q=1}^n (R(t)_{14pq}+R(t)_{23pq})^2 \bigg )^2 \bigg ] \bigg |_{t=\hat{t}} \\ 
&\geq (1-C(n) \, \delta^{\frac{1}{4}}) \, \bigg ( \sum_{p,q=1}^n (\hat{R}_{13pq}-\hat{R}_{24pq})^2 + \sum_{p,q=1}^n (\hat{R}_{14pq}+\hat{R}_{23pq})^2 \bigg ).
\end{align*} 
If we choose $\delta>0$ sufficiently small, then the right hand side is strictly positive. Hence, if $t \in [0,\hat{t})$ is sufficiently close to $\hat{t}$, then 
\begin{align*} 
&R(t)_{1313}+R(t)_{1414}+R(t)_{2323}+R(t)_{2424}-2R(t)_{1234} \\
&< \delta^{\frac{1}{4}} \, \text{\rm scal}(R(t))^{-3} \, \bigg ( \sum_{p,q=1}^n (R(t)_{13pq}-R(t)_{24pq})^2 + \sum_{p,q=1}^n (R(t)_{14pq}+R(t)_{23pq})^2 \bigg ).
\end{align*} 
This contradicts the fact that condition (iv) holds up to time $\hat{t}$. \\ 

Combining Lemma \ref{i.preserved}, Lemma \ref{ii.preserved}, Lemma \ref{iii.preserved}, and Lemma \ref{iv.preserved}, we conclude that $\hat{R}$ lies in the interior of the set $\mathcal{G}_{\hat{t}}^{(0)}$. This contradicts the definition of $\hat{t}$. Therefore, the family of sets $\{\mathcal{G}_t^{(0)}: t \in [0,T]\}$ is invariant under the Hamilton ODE. 

Finally, we verify that, for a curvature tensor in $\mathcal{G}_t^{(0)}$, the isotropic curvature is uniformly bounded from below by a small multiple of the scalar curvature:

\begin{lemma}
\label{uniformly.pic}
We can find a small positive constant $\theta$ (depending only on $\mathcal{K}$ and $T$) such that 
\begin{align*} 
\mathcal{G}_t^{(0)} 
&\subset \{R: R - \theta \, \text{\rm scal} \, \text{\rm id} \owedge \text{\rm id} \in \text{\rm PIC}\} \\ 
&\cap \{R: \text{\rm Ric}_{11}+\text{\rm Ric}_{22} - \theta \, \text{\rm scal} + 4 \geq 0\}. 
\end{align*}
\end{lemma}

\textbf{Proof.} 
Suppose that $R \in \mathcal{G}_t^{(0)}$. By definition, we may write $R = \ell_{a,b}(S)$, where $b=\delta e^{-8t}$, $2a=2b+(n-2)b^2$, $S \in \text{\rm PIC}$, and $\text{\rm Ric}(S)_{11}+\text{\rm Ric}(S)_{22} + 2b^{\frac{5}{4}} \, \text{\rm scal}(S) \geq -2$. Since $S \in \text{\rm PIC}$, Lemma \ref{Ric.four.positive} implies $\text{\rm Ric}(S) \owedge \text{\rm id} \in \text{\rm PIC}$. Consequently, 
\[R - \frac{(n-2)b^2}{2n} \, \text{\rm scal}(S) \, \text{\rm id} \owedge \text{\rm id} = S + b \, \text{\rm Ric}(S) \owedge \text{\rm id} \in \text{\rm PIC}.\] 
Moreover, 
\begin{align*} 
\text{\rm Ric}_{11}+\text{\rm Ric}_{22} 
&= (1+(n-2)b) \, (\text{\rm Ric}(S)_{11}+\text{\rm Ric}(S)_{22}) \\ 
&+ \frac{2}{n} \, (nb + (n-1)(n-2)b^2) \, \text{\rm scal}(S) \\ 
&\geq \Big [ \frac{2}{n} \, (nb + (n-1)(n-2)b^2) - 2(1+(n-2)b) \, b^{\frac{5}{4}} \Big ] \, \text{\rm scal}(S) \\ 
&- 2(1+(n-2)b). 
\end{align*}
From this, the assertion follows.

\section{A one-parameter family of invariant cones $\mathcal{C}(b)$}

\label{first.family.of.cones}

Theorem \ref{rough.pinching} gives uniform lower bounds for the isotropic curvature and for the sum of the two smallest eigenvalues of the Ricci tensor. However, the pinching constant $\theta$ in Theorem \ref{rough.pinching} is extremely small. In this section, we prove estimates of a similar nature, but with a more efficient pinching constant. To that end, we construct a one-parameter family of cones $\mathcal{C}(b)$ which are invariant under the Hamilton ODE $\frac{d}{dt} R = Q(R)$. Like in Section \ref{preliminary.pinching.estimate}, the basic strategy is to consider the cone $\{R \in \text{\rm PIC}: \text{\rm Ric}_{11}+\text{\rm Ric}_{22} \geq 0\}$ (which is invariant under the Hamilton ODE), and deform it inward using the maps $\ell_{a,b}$ introduced in \cite{Bohm-Wilking}. However, these deformed sets by themselves are not preserved under the Hamilton ODE, so we need to impose additional conditions. These conditions are formulated in terms of a tensor $T$. This tensor $T$ is required to satisfy all algebraic properties of a curvature tensor, except for the first Bianchi identity. We say that $T \geq 0$ if $T(\varphi,\bar{\varphi}) \geq 0$ for every complex two-form $\varphi$. As in Definition \ref{pic.and.its.variants}, we say that $T \in \text{\rm PIC}$ if $T(\varphi,\bar{\varphi}) \geq 0$ for every complex two-form of the form $\varphi = (e_1+ie_2) \wedge (e_3+ie_4)$, where $\{e_1,e_2,e_3,e_4\}$ is an orthonormal four-frame. In other words, $T \in \text{\rm PIC}$ if and only if $T_{1313}+T_{1414}+T_{2323}+T_{2424}+2T_{1342}+2T_{1423} \geq 0$ for every orthonormal four-frame $\{e_1,e_2,e_3,e_4\}$.

Throughout this section, we assume that $n \geq 12$. Let $b_{\text{\rm max}} = \frac{1}{4n}$. To each $0 < b \leq b_{\text{\rm max}}$, we associate real numbers $a$, $\gamma$, and $\omega$ by 
\begin{align*} 
a &= \frac{(2+(n-2)b)^2}{2(2+(n-3)b)} \, b, \\  
\gamma &= \frac{b}{2+(n-3)b}, \\ 
\omega^2 &= \frac{27(2+(n-2)b)}{8} \, \frac{(1-4(n-2)b^2)^4(1+(n-2)b)^2}{n^2b^2(2+(n-3)b)^2}.
\end{align*}

\begin{definition}
For $0 < b \leq b_{\text{\rm max}}$ as above, let $\mathcal{E}(b)$ denote the set of all curvature tensors $S$ for which there exists a tensor $T$ satisfying the following conditions: 
\begin{itemize} 
\item[(i)] $T$ satisfies all the algebraic properties of a curvature tensor, except for the first Bianchi identity. Moreover, $T \geq 0$.
\item[(ii)] $S-T \in \text{\rm PIC}$.
\item[(iii)] $\text{\rm Ric}(S)_{11}+\text{\rm Ric}(S)_{22} + \frac{2\gamma}{n} \, \text{\rm scal}(S) \geq 0$. 
\item[(iv)] For every orthonormal frame $\{e_1,\hdots,e_n\}$, the inequality
\begin{align*} 
&\text{\rm Ric}(S)_{22}-\text{\rm Ric}(S)_{11} \\ 
&\leq \omega^{\frac{1}{2}} \, \text{\rm scal}(S)^{\frac{1}{2}} \, \Big ( \sum_{p=3}^n (T_{1p1p}+T_{2p2p}) \Big )^{\frac{1}{2}} + 2(n-2)b \sum_{p=3}^n (T_{1p1p}+T_{2p2p}) 
\end{align*} 
holds.
\end{itemize}
Moreover, we define $\mathcal{C}(b) = \ell_{a,b}(\mathcal{E}(b))$.
\end{definition} 

Clearly, $\mathcal{C}(b)$ is convex for each $b$. The following is the main result of this section:

\begin{theorem}
\label{Cb.is.preserved}
For each $0 < b \leq b_{\text{\rm max}}$, the cone $\mathcal{C}(b)$ is transversally invariant under the Hamilton ODE $\frac{d}{dt} R = Q(R)$.
\end{theorem}

In the remainder of this section, we give the proof of Theorem \ref{Cb.is.preserved}. We begin with an elementary lemma: 

\begin{lemma}
\label{parameters.1}
Let $n \geq 12$, let $0 < b \leq b_{\text{\rm max}}$, and let $a$ and $\omega$ be defined as above. Then 
\[2+(n-8)b-2(n+2)(n-2)b^2-n(n-2)^2b^3 > 0\] 
and 
\[\Big ( \frac{n-2}{n-3}+\frac{2(n-2)}{n} \, (2b+(n-2)a)+8(n-3)b^2 \Big ) \Big ( \frac{n-2}{n-3}+2(2b+(n-2)a) \Big ) \leq \frac{\omega}{2}.\] 
\end{lemma}

\textbf{Proof.} 
Note that $2-2(n+2)(n-2)b^2 \geq 0$ and $(n-8)b-n(n-2)^2b^3 \geq 0$. If we add these inequalities, the first statement follows.

To prove the second statement, we observe that $\frac{2+(n-2)b}{(2+(n-3)b)^2} \geq \frac{1}{2+nb} \geq \frac{4}{9}$. This implies 
\begin{align*} 
\omega 
&= (1-4(n-2)b^2)^2 \, \Big ( \frac{1}{nb} + \frac{n-2}{n} \Big ) \, \sqrt{\frac{27}{8} \, \frac{2+(n-2)b}{(2+(n-3)b)^2}} \\ 
&\geq \Big ( 1-\frac{n-2}{4n^2} \Big )^2 \, \frac{5n-2}{n} \, \sqrt{\frac{3}{2}}. 
\end{align*}
On the other hand, we have $a = \frac{(2+(n-2)b)^2}{2(2+(n-3)b)} \, b \leq \frac{2+nb}{2} \, b \leq \frac{9}{8} \, b$. This implies $2b+(n-2)a \leq na \leq \frac{9n}{8} \, b \leq \frac{9}{32}$. This gives 
\[\frac{n-2}{n-3}+\frac{2(n-2)}{n} \, (2b+(n-2)a)+8(n-3)b^2 \leq \frac{n-2}{n-3} + \frac{9(n-2)}{16n} + \frac{n-3}{2n^2}\] 
and  
\[\frac{n-2}{n-3}+2(2b+(n-2)a) \leq \frac{n-2}{n-3}+\frac{9}{16}.\] 
It is elementary to verify that 
\begin{align*} 
&\Big ( \frac{n-2}{n-3} + \frac{9(n-2)}{16n} + \frac{n-3}{2n^2} \Big ) \, \Big ( \frac{n-2}{n-3}+\frac{9}{16} \Big ) \\ 
&\leq \Big ( 1-\frac{n-2}{4n^2} \Big )^2 \, \frac{5n-2}{n} \, \sqrt{\frac{3}{8}} 
\end{align*} 
for $n \geq 12$. This proves the assertion. \\

\begin{lemma}
\label{quadratic.inequality}
Let $0 < b \leq b_{\text{\rm max}}$, and let $a$ be defined as above. Then 
\[(2b+(n-2)b^2-2a) xy + 2a (x+2)(y+2) + b^2 (x^2+y^2) \geq 0\] 
whenever $x,y \geq -\frac{2(2+(n-2)b)}{2+(n-3)b}$.
\end{lemma}

\textbf{Proof.} 
Let $D = \{(x,y) \in \mathbb{R}^2: x,y \geq -\frac{2(2+(n-2)b)}{2+(n-3)b}\}$, and define $\psi: D \to \mathbb{R}$ by $\psi(x,y) = (2b+(n-2)b^2-2a) xy + 2a (x+2)(y+2) + b^2 (x^2+y^2)$. Clearly, $\psi(x,y) \to \infty$ as $(x,y) \in D$ approaches infinity. Hence, there exists a point in $D$ where $\psi$ attains its minimum. The Hessian of $\psi$ is given by 
\[\begin{bmatrix} 2b^2 & 2b+(n-2)b^2 \\ 2b+(n-2)b^2 & 2b^2 \end{bmatrix}.\]
Since the Hessian of $\psi$ has two eigenvalues of opposite signs, the function $\psi$ attains its minimum on the boundary of $D$. On the other hand, a straightforward calculation gives 
\[\psi \Big ( -\frac{2(2+(n-2)b)}{2+(n-3)b},y \Big ) = b^2 y^2 \geq 0.\] 
Thus, $\psi \geq 0$ on $\partial D$. This implies $\psi \geq 0$ on $D$. \\

\begin{lemma}
\label{Dab.term.1}
Suppose that $S \in \mathcal{E}(b)$. Then $(2b+(n-2)b^2-2a) \, \tracefreeRic(S) \owedge \tracefreeRic(S) + 2a \, \text{\rm Ric}(S) \owedge \text{\rm Ric}(S) + 2b^2 \, \tracefreeRic(S)^2 \owedge \text{\rm id} \in \text{\rm PIC}$. 
\end{lemma} 

\textbf{Proof.} 
Let 
\[U := (2b+(n-2)b^2-2a) \, \tracefreeRic(S) \owedge \tracefreeRic(S) + 2a \, \text{\rm Ric}(S) \owedge \text{\rm Ric}(S) + 2b^2 \, \tracefreeRic(S)^2 \owedge \text{\rm id}.\] 
Let $\zeta,\eta \in \mathbb{C}^n$ be linearly independent vectors satisfying $g(\zeta,\zeta)=g(\zeta,\eta)=g(\eta,\eta)=0$. We claim that $U(\zeta,\eta,\bar{\zeta},\bar{\eta}) \geq 0$. We can find vectors $z,w \in \text{\rm span}\{\zeta,\eta\}$ such that $g(z,\bar{z})=g(w,\bar{w})=2$, $g(z,\bar{w})=0$, and $\text{\rm Ric}(S)(z,\bar{w}) = 0$. The identities $g(\zeta,\zeta)=g(\zeta,\eta)=g(\eta,\eta)=0$ give $g(z,z)=g(z,w)=g(w,w)=0$. Consequently, we may write $z=e_1+ie_2$ and $w=e_3+ie_4$ for some orthonormal four-frame $\{e_1,e_2,e_3,e_4\} \subset \mathbb{R}^n$. Using the identity $\text{\rm Ric}(S)(z,\bar{w}) = 0$, we obtain 
\begin{align*} 
&U(z,w,\bar{z},\bar{w}) \\ 
&= 2(2b+(n-2)b^2-2a) \, (\tracefreeRic(S)_{11}+\tracefreeRic(S)_{22}) (\tracefreeRic(S)_{33}+\tracefreeRic_{44}) \\ 
&+ 4a \, (\text{\rm Ric}(S)_{11}+\text{\rm Ric}(S)_{22}) (\text{\rm Ric}(S)_{33}+\text{\rm Ric}(S)_{44}) \\ 
&+ 4b^2 \, ((\tracefreeRic(S)^2)_{11}+(\tracefreeRic(S)^2)_{22}+(\tracefreeRic(S)^2)_{33}+(\tracefreeRic(S)^2)_{44}) \\ 
&\geq 2(2b+(n-2)b^2-2a) \, (\tracefreeRic(S)_{11}+\tracefreeRic(S)_{22}) (\tracefreeRic(S)_{33}+\tracefreeRic_{44}) \\ 
&+ 4a \, (\text{\rm Ric}(S)_{11}+\text{\rm Ric}(S)_{22}) (\text{\rm Ric}(S)_{33}+\text{\rm Ric}(S)_{44}) \\ 
&+ 2b^2 \, ((\tracefreeRic(S)_{11}+\tracefreeRic(S)_{22})^2+(\tracefreeRic(S)_{33}+\tracefreeRic(S)_{44})^2) \\ 
&= \frac{2 \, \text{\rm scal}(S)^2}{n^2} \, [(2b+(n-2)b^2-2a) xy + 2a (x+2)(y+2) + b^2 (x^2+y^2)], 
\end{align*}
where $x$ and $y$ are defined by $x := \frac{n}{\text{\rm scal}(S)} \, (\tracefreeRic(S)_{11}+\tracefreeRic(S)_{22})$ and $y := \frac{n}{\text{\rm scal}(S)} \, (\tracefreeRic(S)_{33}+\tracefreeRic(S)_{44})$. The condition (iii) implies $x,y \geq -2\gamma-2 = -\frac{2(2+(n-2)b)}{2+(n-3)b}$. Using Lemma \ref{quadratic.inequality}, we obtain $U(z,w,\bar{z},\bar{w}) \geq 0$. Since $\text{\rm span}\{\zeta,\eta\} = \text{\rm span}\{z,w\}$, it follows that $U(\zeta,\eta,\bar{\zeta},\bar{\eta}) \geq 0$. Thus, $U \in \text{\rm PIC}$, as claimed. \\

\begin{lemma}
\label{Dab.term.2}
Suppose that $S \in \mathcal{E}(b) \setminus \{0\}$. Then $D_{a,b}(S)$ lies in the interior of the PIC cone.
\end{lemma}

\textbf{Proof.} 
If $|\tracefreeRic(S)|^2=0$, the assertion is trivial. If $|\tracefreeRic(S)|^2 > 0$, the assertion follows from Lemma \ref{Dab.term.1}, taking into account the fact that 
\begin{align*} 
&nb^2(1-2b)-2(a-b)(1-2b+nb^2) \\ 
&= \frac{b^2}{2+(n-3)b} \, \big ( 2+(n-8)b-2(n+2)(n-2)b^2-n(n-2)^2b^3 \big ) > 0 
\end{align*} 
for $0 < b \leq b_{\text{\rm max}}$. \\

After these preparations, we now give the proof of Theorem \ref{Cb.is.preserved}. It suffices to show that the cone $\mathcal{E}(b)$ is transversally invariant under the ODE $\frac{d}{dt} S = Q(S) + D_{a,b}(S)$ for each $0 < b \leq b_{\text{\rm max}}$. By Lemma \ref{Dab.term.2}, we can find a small positive number $\varepsilon$ (depending on $b$) such that $D_{a,b}(S) - 2\varepsilon \, \text{\rm scal}(S)^2 \, \text{\rm id} \owedge \text{\rm id} \in \text{\rm PIC}$ for all $S \in \mathcal{E}(b)$. We evolve $T$ by the ODE $\frac{d}{dt} T = S^2 + \varepsilon \, \text{\rm scal}(S)^2 \, \text{\rm id} \owedge \text{\rm id}$. 

\begin{lemma} 
\label{positivity.of.T.preserved}
Suppose that $S \in \mathcal{E}(b) \setminus \{0\}$, $T \geq 0$, and $S-T \in \text{\rm PIC}$. Then $\frac{d}{dt} T > 0$. 
\end{lemma}

\textbf{Proof.} 
This follows immediately from the evolution equation for $T$. \\

\begin{lemma} 
\label{S-T.in.PIC.preserved}
Suppose that $S \in \mathcal{E}(b) \setminus \{0\}$, $T \geq 0$, and $S-T \in \text{\rm PIC}$. Then $\frac{d}{dt} (S-T)$ lies in the interior of the tangent cone $T_{S-T} \text{\rm PIC}$.
\end{lemma}

\textbf{Proof.} 
We compute $\frac{d}{dt} (S-T) = S^\# + D_{a,b}(S) - \varepsilon \, \text{\rm scal}(S)^2 \, \text{\rm id} \owedge \text{\rm id}$. Applying Proposition \ref{sharp} with $U := S-T$, we conclude that $S^\# \in T_{S-T} \text{\rm PIC}$. Moreover, by our choice of $\varepsilon$, $D_{a,b}(S) - \varepsilon \, \text{\rm scal}(S)^2 \, \text{\rm id} \owedge \text{\rm id}$ lies in the interior of the PIC cone. Putting these facts together, the assertion follows. \\

\begin{lemma} 
\label{Ricci.bound.preserved}
Suppose that $S \in \mathcal{E}(b) \setminus \{0\}$ and $\text{\rm Ric}(S)_{11}+\text{\rm Ric}(S)_{22} + \frac{2\gamma}{n} \, \text{\rm scal}(S) = 0$ for some pair of orthonormal vectors $\{e_1,e_2\}$. Then $\frac{d}{dt} (\text{\rm Ric}(S)_{11}+\text{\rm Ric}(S)_{22} + \frac{2\gamma}{n} \, \text{\rm scal}(S)) > 0$. 
\end{lemma} 

\textbf{Proof.} 
We may assume that $\{e_1,e_2\}$ are eigenvectors of $\text{\rm Ric}(S)$, and $\text{\rm Ric}(S)_{22} \geq \text{\rm Ric}(S)_{11}$. Let us extend $\{e_1,e_2\}$ to an eigenbasis $\{e_1,\hdots,e_n\}$ of $\text{\rm Ric}(S)$. Clearly, $\text{\rm Ric}(S)_{pp} \geq \text{\rm Ric}(S)_{22}$ for $3 \leq p \leq n$. We compute 
\begin{align*} 
&\frac{1}{2} \frac{d}{dt} \big ( \text{\rm Ric}(S)_{11}+\text{\rm Ric}(S)_{22} + \frac{2\gamma}{n} \, \text{\rm scal}(S) \big ) \\ 
&= \sum_{p=1}^n (S_{1p1p}+S_{2p2p}) \text{\rm Ric}(S)_{pp} \\ 
&+ \frac{4(1+\gamma)}{n^2} \, (a-b) \, \text{\rm scal}(S)^2 + \frac{2\gamma}{n} \, (1-2b) \, |\text{\rm Ric}(S)|^2 - 2b \, (\text{\rm Ric}(S)_{11}^2+\text{\rm Ric}(S)_{22}^2) \\ 
&+ 2 (1+\gamma) \, \frac{n^2b^2-2(n-1)(a-b)(1-2b)}{n(1+2(n-1)a)} \, |\tracefreeRic(S)|^2 \\ 
&= \sum_{p=3}^n (S_{1p1p}+S_{2p2p}) (\text{\rm Ric}(S)_{pp} - \frac{1}{2} \, (\text{\rm Ric}(S)_{11}+\text{\rm Ric}(S)_{22})) \\ 
&+ \frac{4(1+\gamma)}{n^2} \, (a-b) \, \text{\rm scal}(S)^2 + \frac{2\gamma(1+\gamma)}{n^2} \, (1-2b) \, \text{\rm scal}(S)^2  - b \, (\text{\rm Ric}(S)_{22}-\text{\rm Ric}(S)_{11})^2 \\ 
&+ \frac{2\gamma}{n} \, (1-2b) \, |\tracefreeRic(S)|^2 + 2 (1+\gamma) \, \frac{n^2b^2-2(n-1)(a-b)(1-2b)}{n(1+2(n-1)a)} \, |\tracefreeRic(S)|^2. 
\end{align*} 
Using the identity 
\[2(a-b) + \gamma \, (1-2b) = \frac{b(1+(n-2)b)^2}{2+(n-3)b}\] 
and the inequality 
\begin{align*} 
&n^2b^2-2(n-1)(a-b)(1-2b) \\ 
&= \frac{2b^2}{2+(n-3)b} \, \big ( (2n-1)+ (3n-2)(n-2)b + (n-1)(n-2)^2b^2 \big ) \geq 0, 
\end{align*} 
we obtain 
\begin{align*} 
&\frac{1}{2} \frac{d}{dt} \big ( \text{\rm Ric}(S)_{11}+\text{\rm Ric}(S)_{22} + \frac{2\gamma}{n} \, \text{\rm scal}(S) \big ) \\ 
&\geq \sum_{p=3}^n (S_{1p1p}+S_{2p2p}) (\text{\rm Ric}(S)_{pp} - \frac{1}{2} \, (\text{\rm Ric}(S)_{11}+\text{\rm Ric}(S)_{22})) \\ 
&+ \frac{2(1+\gamma)}{n^2} \, \frac{b(1+(n-2)b)^2}{2+(n-3)b} \, \text{\rm scal}(S)^2 - b \, (\text{\rm Ric}(S)_{22}-\text{\rm Ric}(S)_{11})^2. 
\end{align*} 
By Lemma \ref{estimate.for.largest.ricci.eigenvalue}, each eigenvalue of $\text{\rm Ric}(S)$ is bounded from above by the sum of all the other eigenvalues. Therefore, 
\[\sum_{p \in \{1,\hdots,n\} \setminus \{m\}} \text{\rm Ric}(S)_{pp} \geq \text{\rm Ric}(S)_{mm}.\] 
Since $\text{\rm Ric}(S)_{11}+\text{\rm Ric}(S)_{22} \leq 0$, we deduce that 
\begin{align*} 
&\sum_{p \in \{3,\hdots,n\} \setminus \{m\}}  (\text{\rm Ric}(S)_{pp} - \frac{1}{2} \, (\text{\rm Ric}(S)_{11}+\text{\rm Ric}(S)_{22})) \\ 
&\geq \text{\rm Ric}(S)_{mm} - \frac{1}{2} \, (\text{\rm Ric}(S)_{11}+\text{\rm Ric}(S)_{22}) 
\end{align*} 
for each $m \in \{3,\hdots,n\}$. Since $S-T \in \text{\rm PIC}$, it follows that 
\[\sum_{p=3}^n (S_{1p1p}+S_{2p2p}-T_{1p1p}-T_{2p2p}) \, (\text{\rm Ric}(S)_{pp} - \frac{1}{2} \, (\text{\rm Ric}(S)_{11}+\text{\rm Ric}(S)_{22})) \geq 0.\] 
Since $T \geq 0$, we know that $T_{1p1p}+T_{2p2p} \geq 0$ for $3 \leq p \leq n$. This gives 
\begin{align*} 
&\sum_{p=3}^n (S_{1p1p}+S_{2p2p}) (\text{\rm Ric}(S)_{pp} - \frac{1}{2} \, (\text{\rm Ric}(S)_{11}+\text{\rm Ric}(S)_{22})) \\ 
&\geq \sum_{p=3}^n (T_{1p1p}+T_{2p2p}) (\text{\rm Ric}(S)_{pp} - \frac{1}{2} \, (\text{\rm Ric}(S)_{11}+\text{\rm Ric}(S)_{22})) \\ 
&\geq \frac{1}{2} \sum_{p=3}^n (T_{1p1p}+T_{2p2p}) \, (\text{\rm Ric}(S)_{22}-\text{\rm Ric}(S)_{11}). 
\end{align*}
Putting these facts together, we obtain
\begin{align*} 
&\frac{1}{2} \frac{d}{dt} \big ( \text{\rm Ric}(S)_{11}+\text{\rm Ric}(S)_{22} + \frac{2\gamma}{n} \, \text{\rm scal}(S) \big ) \\ 
&\geq \frac{1}{2} \sum_{p=3}^n (T_{1p1p}+T_{2p2p}) \, (\text{\rm Ric}(S)_{22}-\text{\rm Ric}(S)_{11}) \\ 
&+ \frac{2(1+\gamma)}{n^2} \, \frac{b(1+(n-2)b)^2}{2+(n-3)b} \, \text{\rm scal}(S)^2 - b \, (\text{\rm Ric}(S)_{22}-\text{\rm Ric}(S)_{11})^2. 
\end{align*} 
At this point, we distinguish two cases: 

\textit{Case 1:} Suppose that $\sum_{p=3}^n (T_{1p1p}+T_{2p2p}) \geq 2b \, (\text{\rm Ric}(S)_{22}-\text{\rm Ric}(S)_{11})$. In this case, we clearly have 
\[\frac{1}{2} \frac{d}{dt} \big ( \text{\rm Ric}(S)_{11}+\text{\rm Ric}(S)_{22} + \frac{2\gamma}{n} \, \text{\rm scal}(S) \big ) > 0.\] 

\textit{Case 2:} Suppose next that $\sum_{p=3}^n (T_{1p1p}+T_{2p2p}) < 2b \, (\text{\rm Ric}(S)_{22}-\text{\rm Ric}(S)_{11})$. In this case, the condition (iv) gives 
\begin{align*} 
&\text{\rm Ric}(S)_{22}-\text{\rm Ric}(S)_{11} \\ 
&\leq \omega^{\frac{1}{2}} \, \text{\rm scal}(S)^{\frac{1}{2}} \, \Big ( \sum_{p=3}^n (T_{1p1p}+T_{2p2p}) \Big )^{\frac{1}{2}} + 2(n-2)b \sum_{p=3}^n (T_{1p1p}+T_{2p2p}) \\ 
&< \omega^{\frac{1}{2}} \, \text{\rm scal}(S)^{\frac{1}{2}} \, \Big ( \sum_{p=3}^n (T_{1p1p}+T_{2p2p}) \Big )^{\frac{1}{2}} + 4(n-2)b^2 \, (\text{\rm Ric}(S)_{22}-\text{\rm Ric}(S)_{11}). 
\end{align*} 
Consequently, 
\[\sum_{p=3}^n (T_{1p1p}+T_{2p2p}) > \frac{(1-4(n-2)b^2)^2}{\omega} \, \frac{(\text{\rm Ric}(S)_{22}-\text{\rm Ric}(S)_{11})^2}{\text{\rm scal}(S)}.\] 
Using the elementary inequality $2x+y \geq 3 (x^2y)^{\frac{1}{3}}$ for $x,y \geq 0$, we obtain  
\begin{align*} 
&\frac{(1-4(n-2)b^2)^2}{2\omega} \, \frac{(\text{\rm Ric}(S)_{22}-\text{\rm Ric}(S)_{11})}{\text{\rm scal}(S)} \\ 
&+ \frac{2(1+\gamma)}{n^2} \, \frac{b(1+(n-2)b)^2}{2+(n-3)b} \, \frac{\text{\rm scal}(S)^2}{(\text{\rm Ric}(S)_{22}-\text{\rm Ric}(S)_{11})^2} \\ 
&\geq 3 \, \Big ( \frac{(1-4(n-2)b^2)^4}{8\omega^2} \, \frac{1+\gamma}{n^2} \, \frac{b(1+(n-2)b)^2}{2+(n-3)b} \Big )^{\frac{1}{3}} \\ 
&= b, 
\end{align*}
where in the last step we have used the definitions of $\gamma$ and $\omega$. Thus, we conclude that 
\begin{align*} 
&\frac{1}{2} \frac{d}{dt} \big ( \text{\rm Ric}(S)_{11}+\text{\rm Ric}(S)_{22} + \frac{2\gamma}{n} \, \text{\rm scal}(S) \big ) \\ 
&> \frac{(1-4(n-2)b^2)^2}{2\omega} \, \frac{(\text{\rm Ric}(S)_{22}-\text{\rm Ric}(S)_{11})^3}{\text{\rm scal}(S)} \\ 
&+ \frac{2(1+\gamma)}{n^2} \, \frac{b(1+(n-2)b)^2}{2+(n-3)b} \, \text{\rm scal}(S)^2  - b \, (\text{\rm Ric}(S)_{22} - \text{\rm Ric}(S)_{11})^2 \\ 
&\geq 0. 
\end{align*} 
This proves the assertion. \\

\begin{lemma}
Suppose that $S \in \mathcal{E}(b) \setminus \{0\}$, $T \geq 0$, and $S-T \in \text{\rm PIC}$. Moreover, suppose that 
\begin{align*} 
&\text{\rm Ric}(S)_{22}-\text{\rm Ric}(S)_{11} \\ 
&= \omega^{\frac{1}{2}} \, \text{\rm scal}(S)^{\frac{1}{2}} \, \Big ( \sum_{p=3}^n (T_{1p1p}+T_{2p2p}) \Big )^{\frac{1}{2}} + 2(n-2)b \sum_{p=3}^n (T_{1p1p}+T_{2p2p}). 
\end{align*} 
Then 
\begin{align*} 
&\frac{d}{dt}(\text{\rm Ric}(S)_{22}-\text{\rm Ric}(S)_{11}) \\ 
&< \frac{d}{dt} \bigg [ \omega^{\frac{1}{2}} \, \text{\rm scal}(S)^{\frac{1}{2}} \, \Big ( \sum_{p=3}^n (T_{1p1p}+T_{2p2p}) \Big )^{\frac{1}{2}} + 2(n-2)b \sum_{p=3}^n (T_{1p1p}+T_{2p2p}) \bigg ]. 
\end{align*} 
\end{lemma}

\textbf{Proof.} 
If $\sum_{p=3}^n (T_{1p1p}+T_{2p2p}) = 0$, then the right hand side is infinite, and the inequality is trivially true. In the following, we assume that $\sum_{p=3}^n (T_{1p1p}+T_{2p2p}) > 0$. Clearly, $\text{\rm Ric}(S)_{22} - \text{\rm Ric}(S)_{11} \geq 0$ and $\text{\rm Ric}(S)_{12}=0$. Without loss of generality, we may assume that $\text{\rm Ric}(S)_{pq}=0$ for $3 \leq p<q \leq n$. We may write 
\[\frac{d}{dt} (\text{\rm Ric}(S)_{22}-\text{\rm Ric}(S)_{11}) = J_1 + J_2 + J_3 + J_4 + J_5 + J_6,\] 
where 
\begin{align*} 
J_1 &:= 2 \sum_{p=3}^n (S_{2p2p}-S_{1p1p}) \, \text{\rm Ric}(S)_{pp}, \\ 
J_2 &:= -4 \sum_{p=3}^n S_{121p} \, \text{\rm Ric}(S)_{2p}, \\ 
J_3 &:= 4 \sum_{p=3}^n S_{212p} \, \text{\rm Ric}(S)_{1p}, \\ 
J_4 &:= -2 S_{1212} \, (\text{\rm Ric}(S)_{22}-\text{\rm Ric}(S)_{11}), \\ 
J_5 &:= 4b \, ((\text{\rm Ric}(S)^2)_{11} - (\text{\rm Ric}(S)^2)_{22}), \\
J_6 &:= \frac{4}{n} \, (2b+(n-2)a) \, \text{\rm scal}(S) \, (\text{\rm Ric}(S)_{22}-\text{\rm Ric}(S)_{11}).
\end{align*}
The terms $J_1,J_2,J_3$ can be estimated by 
\begin{align*} 
J_1 &\leq \tau \sum_{p=3}^n (S_{2p2p}-S_{1p1p})^2 + \frac{1}{\tau} \sum_{p=3}^n (\text{\rm Ric}(S)_{pp})^2, \\ 
J_2 &\leq \frac{2(n-2)\tau}{n-3} \sum_{p=3}^n (S_{121p})^2 + \frac{2(n-3)}{(n-2)\tau} \sum_{p=3}^n (\text{\rm Ric}(S)_{2p})^2, \\ 
J_3 &\leq \frac{2(n-2)\tau}{n-3} \sum_{p=3}^n (S_{212p})^2 + \frac{2(n-3)}{(n-2)\tau} \sum_{p=3}^n (\text{\rm Ric}(S)_{1p})^2, 
\end{align*} 
where $\tau$ is an arbitrary positive number. To estimate $J_4$, we observe that 
\begin{align*} 
0 
&\leq \sum_{3 \leq p,q \leq n, \, p \neq q} (S_{1212}+S_{1p1p}+S_{2q2q}+S_{pqpq}) \\ 
&= (n-2)(n-3) S_{1212} + (n-4) \sum_{p=3}^n (S_{1p1p}+S_{2p2p}) + \sum_{p=3}^n \text{\rm Ric}(S)_{pp} 
\end{align*} 
since $S \in \text{\rm PIC}$. This gives 
\begin{align*} 
J_4 
&\leq \frac{2(n-4)}{(n-2)(n-3)} \sum_{p=3}^n (S_{1p1p}+S_{2p2p}) \, (\text{\rm Ric}(S)_{22}-\text{\rm Ric}(S)_{11}) \\ 
&+ \frac{2}{(n-2)(n-3)} \, (\text{\rm Ric}(S)_{22}-\text{\rm Ric}(S)_{11}) \, \Big ( \sum_{p=3}^n \text{\rm Ric}(S)_{pp} \Big ) \\ 
&= \frac{2(n-4)}{(n-2)(n-3)} \sum_{p=3}^n (S_{1p1p}+S_{2p2p}) \, (\text{\rm Ric}(S)_{22}-\text{\rm Ric}(S)_{11}) \\ 
&+ \frac{2}{(n-2)(n-3)}  \, \Big ( \sum_{p=3}^n (S_{2p2p}-S_{1p1p}) \Big ) \, \Big ( \sum_{p=3}^n \text{\rm Ric}(S)_{pp} \Big ) \\ 
\end{align*}
hence 
\begin{align*}
J_4 
&\leq \frac{(n-2)\tau}{n-3} \sum_{p=3}^n (S_{1p1p}+S_{2p2p})^2 + \frac{(n-4)^2}{(n-2)^2(n-3)\tau} \, (\text{\rm Ric}(S)_{22}-\text{\rm Ric}(S)_{11})^2 \\ 
&+ \frac{\tau}{n-3} \sum_{p=3}^n (S_{2p2p}-S_{1p1p})^2 + \frac{1}{(n-3)\tau} \sum_{p=3}^n (\text{\rm Ric}(S)_{pp})^2 
\end{align*} 
for each $\tau>0$. Adding these inequalities gives 
\begin{align*} 
&J_1+J_2+J_3+J_4 \\ 
&\leq \frac{2(n-2)\tau}{n-3} \sum_{p=3}^n ((S_{1p1p})^2+(S_{2p2p})^2+(S_{121p})^2+(S_{212p})^2) \\ 
&+ \frac{n-2}{(n-3)\tau} \sum_{p=3}^n (\text{\rm Ric}(S)_{pp})^2 + \frac{2(n-3)}{(n-2)\tau} \sum_{p=3}^n ((\text{\rm Ric}(S)_{1p})^2 +  (\text{\rm Ric}(S)_{2p})^2) \\ 
&+ \frac{(n-4)^2}{(n-2)^2(n-3)\tau} \, (\text{\rm Ric}(S)_{22}-\text{\rm Ric}(S)_{11})^2 \\ 
&\leq \frac{(n-2)\tau}{n-3} \sum_{p=3}^n ((S^2)_{1p1p}+(S^2)_{2p2p}) \\ 
&+ \frac{n-2}{(n-3)\tau} \, |\text{\rm Ric}(S)|^2 - \frac{n-2}{2(n-3)\tau} \, (\text{\rm Ric}(S)_{11}+\text{\rm Ric}(S)_{22})^2 
\end{align*} 
for each $\tau>0$. Moreover, 
\begin{align*} 
J_5 
&\leq -4b \, (\text{\rm Ric}(S)_{22}-\text{\rm Ric}(S)_{11}) \, (\text{\rm Ric}(S)_{11}+\text{\rm Ric}(S)_{22}) + 4b \sum_{p=3}^n (\text{\rm Ric}(S)_{1p})^2 \\ 
&= -4b \sum_{p=3}^n (S_{2p2p}-S_{1p1p}) \, (\text{\rm Ric}(S)_{11}+\text{\rm Ric}(S)_{22}) + 4b \sum_{p=3}^n \Big ( S_{12p2} + \sum_{q=3}^n S_{1qpq} \Big )^2 \\ 
&\leq 8(n-3)b^2\tau \sum_{p=3}^n (S_{2p2p}-S_{1p1p})^2 + \frac{n-2}{2(n-3)\tau} \, (\text{\rm Ric}(S)_{11}+\text{\rm Ric}(S)_{22})^2 \\ 
&+ 4(n-2)b \sum_{p=3}^n \Big ( (S_{12p2})^2 + \sum_{q=3}^n (S_{1qpq})^2 \Big ) \\ 
&\leq (8(n-3)b^2\tau+2(n-2)b) \sum_{p=3}^n ((S^2)_{1p1p}+(S^2)_{2p2p}) \\ 
&+ \frac{n-2}{2(n-3)\tau} \, (\text{\rm Ric}(S)_{11}+\text{\rm Ric}(S)_{22})^2 
\end{align*}
for each $\tau>0$. Finally, 
\begin{align*} 
J_6 
&= \frac{4}{n} \, (2b+(n-2)a) \, \Big ( \sum_{p=3}^n (S_{2p2p}-S_{1p1p}) \Big ) \, \Big ( \sum_{p=1}^n \text{\rm Ric}(S)_{pp} \Big ) \\
&\leq \frac{2(n-2)\tau}{n} \, (2b+(n-2)a) \sum_{p=3}^n (S_{2p2p}-S_{1p1p})^2 \\ 
&+ \frac{2(2b+(n-2)a)}{\tau} \, |\text{\rm Ric}(S)|^2 \\ 
&\leq \frac{2(n-2)\tau}{n} \, (2b+(n-2)a) \sum_{p=3}^n ((S^2)_{1p1p}+(S^2)_{2p2p}) \\ 
&+ \frac{2(2b+(n-2)a)}{\tau} \, |\text{\rm Ric}(S)|^2 
\end{align*} 
for each $\tau>0$. Putting these facts together, we obtain 
\begin{align*} 
&\frac{d}{dt} (\text{\rm Ric}(S)_{22}-\text{\rm Ric}(S)_{11}) \\ 
&\leq \bigg [ \Big ( \frac{n-2}{n-3}+\frac{2(n-2)}{n} \, (2b+(n-2)a)+8(n-3)b^2 \Big ) \tau + 2(n-2)b \bigg ] \\ 
&\hspace{10mm} \cdot \sum_{p=3}^n ((S^2)_{1p1p}+(S^2)_{2p2p}) \\ 
&+ \Big ( \frac{n-2}{n-3}+2(2b+(n-2)a) \Big ) \, \frac{1}{\tau} \, |\text{\rm Ric}(S)|^2 
\end{align*} 
for each $\tau>0$. If we put $\tau = 2\sigma \, \big ( \frac{n-2}{n-3}+2(2b+(n-2)a) \big )$ and use the inequality 
\[\Big ( \frac{n-2}{n-3}+\frac{2(n-2)}{n} \, (2b+(n-2)a)+8(n-3)b^2 \Big ) \Big ( \frac{n-2}{n-3}+2(2b+(n-2)a) \Big ) \leq \frac{\omega}{2}\] 
(cf. Lemma \ref{parameters.1}), we conclude that 
\begin{align*} 
&\frac{d}{dt} (\text{\rm Ric}(S)_{22}-\text{\rm Ric}(S)_{11}) \\ 
&\leq (\omega\sigma+2(n-2)b) \sum_{p=3}^n ((S^2)_{1p1p}+(S^2)_{2p2p}) + \frac{1}{2\sigma} \, |\text{\rm Ric}(S)|^2
\end{align*} 
for each $\sigma>0$. 

On the other hand, Lemma \ref{Dab.term.2} implies $\text{\rm scal}(D_{a,b}(S)) > 0$, hence $\frac{d}{dt} \text{\rm scal}(S) > 2 \, |\text{\rm Ric}(S)|^2$. Using the inequality $\frac{d}{dt} \sum_{p=3}^n (T_{1p1p}+T_{2p2p}) > \sum_{p=3}^n ((S^2)_{1p1p}+(S^2)_{2p2p})$, we obtain 
\begin{align*}
&\frac{d}{dt} \bigg [ \omega^{\frac{1}{2}} \, \text{\rm scal}(S)^{\frac{1}{2}} \, \Big ( \sum_{p=3}^n (T_{1p1p}+T_{2p2p}) \Big )^{\frac{1}{2}} + 2(n-2)b \sum_{p=3}^n (T_{1p1p}+T_{2p2p}) \bigg ] \\ 
&> (\omega\sigma+2(n-2)b) \sum_{p=3}^n ((S^2)_{1p1p}+(S^2)_{2p2p}) + \frac{1}{2\sigma} \, |\text{\rm Ric}(S)|^2, 
\end{align*} 
where $\sigma := \frac{1}{2} \, \omega^{-\frac{1}{2}} \, \text{\rm scal}(S)^{\frac{1}{2}} \, \big ( \sum_{p=3}^n (T_{1p1p}+T_{2p2p}) \big )^{-\frac{1}{2}}$. Putting these facts together, the assertion follows.

\section{A one-parameter family of invariant cones $\tilde{\mathcal{C}}(b)$ which pinch toward PIC1}

\label{second.family.of.cones}

In this section, we construct another one-parameter family of cones $\tilde{\mathcal{C}}(b)$ which are invariant under the Hamilton ODE. These cones pinch toward the PIC1 cone and can be joined continuously to the family $\mathcal{C}(b)$ constructed in Section \ref{first.family.of.cones}. More precisely, we start from the cone $\mathcal{C}(b_{\text{\rm max}}) \cap \text{\rm PIC1}$ (which is invariant under the Hamilton ODE), and deform it outward in such a way that the deformed cones are all preserved under the ODE. To that end, we consider curvature conditions which interpolate between the PIC condition and the PIC1 condition. To show that these conditions are preserved, we will use Proposition \ref{interpolation.between.pic.and.pic1} in a crucial way. We also make use of the quantitative estimates in Section \ref{first.family.of.cones}. 

Throughout this section, we assume that $n \geq 12$. We define 
\[\tilde{b}_{\text{\rm max}} = \begin{cases} \frac{1}{3(n-6)^2} & \text{\rm $n \geq 13$} \\ \frac{1}{115} & \text{\rm $n=12$.} \end{cases}\]

\begin{definition} 
Assume that $0 < b \leq \tilde{b}_{\text{\rm max}}$, and let $2a=2b+(n-2)b^2$. We denote by $\tilde{\mathcal{E}}(b)$ the set of all algebraic curvature tensors $S$ such that $\ell_{a,b}(S) \in \mathcal{C}(b_{\text{\rm max}})$ and 
\begin{align*} 
Z 
&:= S_{1313}+\lambda^2 S_{1414}+S_{2323}+\lambda^2 S_{2424}-2\lambda S_{1234} \\ 
&+ \sqrt{2a} \, (1-\lambda^2) \, (\text{\rm Ric}(S)_{11}+\text{\rm Ric}(S)_{22}) \geq 0 
\end{align*}
for every orthonormal four-frame $\{e_1,e_2,e_3,e_4\}$ and every $\lambda \in [0,1]$. Moreover, we define $\tilde{\mathcal{C}}(b) = \ell_{a,b}(\tilde{\mathcal{E}}(b))$. 
\end{definition} 

Clearly, $\tilde{\mathcal{C}}(b)$ is convex for each $b$. In the first step, we verify that the family of cones $\mathcal{C}(b)$, $0 < b \leq b_{\text{\rm max}}$, can be joined with the family of cones $\tilde{\mathcal{C}}(b)$, $0 < b \leq \tilde{b}_{\text{\rm max}}$, constructed in Section \ref{first.family.of.cones}. Recall that $b_{\text{\rm max}} = \frac{1}{4n}$. For abbreviation, we define 
\[a_{\text{\rm max}} = \frac{(2+(n-2)b_{\text{\rm max}})^2}{2(2+(n-3)b_{\text{\rm max}})} \, b_{\text{\rm max}}\] 
and 
\[\gamma_{\text{\rm max}} = \frac{b_{\text{\rm max}}}{2+(n-3)b_{\text{\rm max}}}.\]

\begin{lemma}
\label{parameters.2}
Let $\tilde{b}_{\text{\rm max}}$ be defined as above, and let $\tilde{a}_{\text{\rm max}}$ be defined by $2\tilde{a}_{\text{\rm max}} = 2\tilde{b}_{\text{\rm max}} + (n-2)\tilde{b}_{\text{\rm max}}^2$. Then 
\[\frac{1+(n-2)b_{\text{\rm max}}}{1+(n-2)\tilde{b}_{\text{\rm max}}} \, \sqrt{2\tilde{a}_{\text{\rm max}}} \geq \frac{n^2-5n+4}{n^2-7n+14} \, \frac{1}{n-4},\] 
\[\frac{a_{\text{\rm max}}-\tilde{a}_{\text{\rm max}}}{1+2(n-1)\tilde{a}_{\text{\rm max}}} - \frac{b_{\text{\rm max}}-\tilde{b}_{\text{\rm max}}}{1+(n-2)\tilde{b}_{\text{\rm max}}} \geq 0,\] 
and
\begin{align*}
&2 \, \Big ( \frac{a_{\text{\rm max}}-\tilde{a}_{\text{\rm max}}}{1+2(n-1)\tilde{a}_{\text{\rm max}}} - (1+\gamma_{\text{\rm max}}) \, \frac{b_{\text{\rm max}}-\tilde{b}_{\text{\rm max}}}{1+(n-2)\tilde{b}_{\text{\rm max}}} \Big ) \\ 
&+ \Big ( \frac{2(n-1)(a_{\text{\rm max}}-\tilde{a}_{\text{\rm max}})}{1+2(n-1)\tilde{a}_{\text{\rm max}}} - \frac{(n-2)(b_{\text{\rm max}}-\tilde{b}_{\text{\rm max}})}{1+(n-2)\tilde{b}_{\text{\rm max}}} \Big ) \, \sqrt{2\tilde{a}_{\text{\rm max}}} \\ 
&\geq \frac{1+(n-2)b_{\text{\rm max}}}{1+(n-2)\tilde{b}_{\text{\rm max}}} \, \frac{n\sqrt{2\tilde{a}_{\text{\rm max}}}}{n^2-5n+4}. 
\end{align*}
\end{lemma}

\textbf{Proof.} 
To verify the first statement, we observe that $\tilde{a}_{\text{\rm max}} \geq (1+(n-2)\tilde{b}_{\text{\rm max}})^2 \, (1-\frac{3(n-2)}{4} \, \tilde{b}_{\text{\rm max}})^2 \, \tilde{b}_{\text{\rm max}}$. This implies 
\begin{align*} 
\frac{1+(n-2)b_{\text{\rm max}}}{1+(n-2)\tilde{b}_{\text{\rm max}}} \, \sqrt{2\tilde{a}_{\text{\rm max}}} 
&\geq (1+(n-2)b_{\text{\rm max}}) \, \Big ( 1-\frac{3(n-2)}{4} \, \tilde{b}_{\text{\rm max}} \Big ) \, \sqrt{2\tilde{b}_{\text{\rm max}}} \\ 
&\geq \Big ( 1+\frac{n-2}{4n} \Big ) \, \Big ( 1-\frac{n-2}{4(n-6)^2} \Big ) \, \frac{1}{n-6} \, \sqrt{\frac{2}{3}} \\ 
&\geq \frac{n^2-5n+4}{n^2-7n+14} \, \frac{1}{n-4} 
\end{align*} 
for $n \geq 13$. Moreover, it is straightforward to verify that the first inequality holds for $n=12$.

To prove the second statement, we observe that $1+2(n-1)\tilde{a}_{\text{\rm max}} \leq (1+n\tilde{b}_{\text{\rm max}}) (1+(n-2)\tilde{b}_{\text{\rm max}})$. This implies 
\begin{align*} 
&\frac{a_{\text{\rm max}}-\tilde{a}_{\text{\rm max}}}{1+2(n-1)\tilde{a}_{\text{\rm max}}} - \frac{b_{\text{\rm max}}-\tilde{b}_{\text{\rm max}}}{1+(n-2)\tilde{b}_{\text{\rm max}}} \\ 
&\geq \frac{(a_{\text{\rm max}}-\tilde{a}_{\text{\rm max}}) - (1+n\tilde{b}_{\text{\rm max}}) (b_{\text{\rm max}}-\tilde{b}_{\text{\rm max}})}{(1+n\tilde{b}_{\text{\rm max}})(1+(n-2)\tilde{b}_{\text{\rm max}})} \\ 
&\geq \frac{a_{\text{\rm max}} - (1+n\tilde{b}_{\text{\rm max}}) b_{\text{\rm max}}}{(1+n\tilde{b}_{\text{\rm max}})(1+(n-2)\tilde{b}_{\text{\rm max}})} \\ 
&\geq 0
\end{align*} 
for $n \geq 12$. From this, the second statement follows. 

It remains to verify the third statement. Arguing as above, we obtain 
\begin{align*} 
&\frac{a_{\text{\rm max}}-\tilde{a}_{\text{\rm max}}}{1+2(n-1)\tilde{a}_{\text{\rm max}}} - (1+\gamma_{\text{\rm max}}) \, \frac{b_{\text{\rm max}}-\tilde{b}_{\text{\rm max}}}{1+(n-2)\tilde{b}_{\text{\rm max}}} \\ 
&\geq \frac{(a_{\text{\rm max}}-\tilde{a}_{\text{\rm max}}) - (1+n\tilde{b}_{\text{\rm max}}) (1+\gamma_{\text{\rm max}}) (b_{\text{\rm max}}-\tilde{b}_{\text{\rm max}})}{(1+n\tilde{b}_{\text{\rm max}})(1+(n-2)\tilde{b}_{\text{\rm max}})} \\ 
&\geq \frac{a_{\text{\rm max}} - (1+n\tilde{b}_{\text{\rm max}}) (1+\gamma_{\text{\rm max}}) b_{\text{\rm max}}}{(1+n\tilde{b}_{\text{\rm max}})(1+(n-2)\tilde{b}_{\text{\rm max}})}. 
\end{align*} 
Using the identity $a_{\text{\rm max}} = (1+\frac{n-2}{2} \, b_{\text{\rm max}}) (1+\gamma_{\text{\rm max}}) b_{\text{\rm max}}$ and the inequality $\frac{n-2}{2} \, b_{\text{\rm max}}-n\tilde{b}_{\text{\rm max}} \geq \frac{n-8}{8n} \, (1+n\tilde{b}_{\text{\rm max}})$ for $n \geq 20$, we obtain 
\begin{align*} 
&\frac{a_{\text{\rm max}}-\tilde{a}_{\text{\rm max}}}{1+2(n-1)\tilde{a}_{\text{\rm max}}} - (1+\gamma_{\text{\rm max}}) \, \frac{b_{\text{\rm max}}-\tilde{b}_{\text{\rm max}}}{1+(n-2)\tilde{b}_{\text{\rm max}}} \\ 
&\geq \frac{(\frac{n-2}{2} \, b_{\text{\rm max}}-n\tilde{b}_{\text{\rm max}}) (1+\gamma_{\text{\rm max}}) b_{\text{\rm max}}}{(1+n\tilde{b}_{\text{\rm max}})(1+(n-2)\tilde{b}_{\text{\rm max}})} \\ 
&\geq \frac{1}{1+(n-2)\tilde{b}_{\text{\rm max}}} \, \frac{n-8}{32n^2}
\end{align*} 
for $n \geq 20$. On the other hand, the inequality $2\tilde{a}_{\text{\rm max}} \leq 3\tilde{b}_{\text{\rm max}} = \frac{1}{(n-6)^2}$ gives 
\[\frac{1+(n-2) b_{\text{\rm max}}}{1+(n-2)\tilde{b}_{\text{\rm max}}} \, \frac{n\sqrt{2\tilde{a}_{\text{\rm max}}}}{n^2-5n+4} \leq \frac{1}{1+(n-2)\tilde{b}_{\text{\rm max}}} \, \frac{5}{4(n-6)} \, \frac{n}{n^2-5n+4}\] 
for $n \geq 13$. Since $\frac{n-8}{16n^2} \geq \frac{5}{4(n-6)} \, \frac{n}{n^2-5n+4}$ for $n \geq 36$, we conclude that the third statement holds for $n \geq 36$. It is straightforward to verify that the third inequality holds for $12 \leq n \leq 35$. \\

\begin{proposition}
\label{pinching.cones.can.be.joined.together}
Let $b_{\text{\rm max}}$ and $\tilde{b}_{\text{\rm max}}$ be defined as above. Then $\tilde{\mathcal{C}}(\tilde{b}_{\text{\rm max}}) = \mathcal{C}(b_{\text{\rm max}})$. 
\end{proposition}

\textbf{Proof.} 
The inclusion $\tilde{\mathcal{C}}(\tilde{b}_{\text{\rm max}}) \subset \mathcal{C}(b_{\text{\rm max}})$ follows immediately from the definition. We now prove the reverse inclusion $\mathcal{C}(b_{\text{\rm max}}) \subset \tilde{\mathcal{C}}(\tilde{b}_{\text{\rm max}})$. For abbreviation, we put $b = \tilde{b}_{\text{\rm max}}$ and $a=\tilde{a}_{\text{\rm max}}$. Moreover, suppose that $S$ is an algebraic curvature tensor such that $\ell_{a,b}(S) \in \mathcal{C}(b_{\text{\rm max}})$. We claim that $S \in \tilde{\mathcal{E}}(b)$. To verify this claim, we need to show that 
\begin{align*} 
Z 
&:= S_{1313}+\lambda^2 S_{1414}+S_{2323}+\lambda^2 S_{2424}-2\lambda S_{1234} \\ 
&+ \sqrt{2a} \, (1-\lambda^2) \, (\text{\rm Ric}(S)_{11}+\text{\rm Ric}(S)_{22}) \geq 0 
\end{align*}
for every orthonormal four-frame $\{e_1,e_2,e_3,e_4\}$ and every $\lambda \in [0,1]$. To prove this, we consider the algebraic curvature tensor $T = \ell_{a_{\text{\rm max}},b_{\text{\rm max}}}^{-1}(\ell_{a,b}(S)) \in \mathcal{E}(b_{\text{\rm max}})$. Clearly, $T \in \text{\rm PIC}$ and $\text{\rm Ric}(T)_{11}+\text{\rm Ric}(T)_{22} + \frac{2\gamma_{\text{\rm max}}}{n} \, \text{\rm scal}(T) \geq 0$.
 
Note that $S$ and $T$ have the same Weyl tensor. Moreover, $\tracefreeRic(S) = \frac{1+(n-2)b_{\text{\rm max}}}{1+(n-2)b} \, \tracefreeRic(T)$ and $\text{\rm scal}(S) = \frac{1+2(n-1)a_{\text{\rm max}}}{1+2(n-1)a} \, \text{\rm scal}(T)$. This gives 
\begin{align*} 
\text{\rm Ric}(S) 
&= \frac{1+(n-2)b_{\text{\rm max}}}{1+(n-2)b} \, \text{\rm Ric}(T) \\ 
&+ \frac{1}{n} \, \Big ( \frac{2(n-1)(a_{\text{\rm max}}-a)}{1+2(n-1)a} - \frac{(n-2)(b_{\text{\rm max}}-b)}{1+(n-2)b} \Big ) \, \text{\rm scal}(T) \, \text{\rm id} 
\end{align*} 
and 
\begin{align*} 
S &= T + \frac{b_{\text{\rm max}}-b}{1+(n-2)b} \, \text{\rm Ric}(T) \owedge \text{\rm id} \\ 
&+ \frac{1}{n} \, \Big ( \frac{a_{\text{\rm max}}-a}{1+2(n-1)a} - \frac{b_{\text{\rm max}}-b}{1+(n-2)b} \Big ) \, \text{\rm scal}(T) \, \text{\rm id} \owedge \text{\rm id}. 
\end{align*}
Consequently, 
\begin{align*} 
Z &= T_{1313}+\lambda^2 T_{1414}+T_{2323}+\lambda^2 T_{2424}-2\lambda T_{1234} \\ 
&+ \frac{b_{\text{\rm max}}-b}{1+(n-2)b} \, [(1+\lambda^2) \, \text{\rm Ric}(T)_{11}+(1+\lambda^2) \, \text{\rm Ric}(T)_{22} \\ 
&\hspace{30mm} +2 \, \text{\rm Ric}(T)_{33}+2\lambda^2 \, \text{\rm Ric}(T)_{44}] \\ 
&+ \frac{4}{n} \, \Big ( \frac{a_{\text{\rm max}}-a}{1+2(n-1)a} - \frac{b_{\text{\rm max}}-b}{1+(n-2)b} \Big ) \, (1+\lambda^2) \, \text{\rm scal}(T) \\ 
&+ \frac{1+(n-2)b_{\text{\rm max}}}{1+(n-2)b} \, \sqrt{2a} \, (1-\lambda^2) \, (\text{\rm Ric}(T)_{11}+\text{\rm Ric}(T)_{22}) \\ 
&+ \frac{2}{n} \, \Big ( \frac{2(n-1)(a_{\text{\rm max}}-a)}{1+2(n-1)a} - \frac{(n-2)(b_{\text{\rm max}}-b)}{1+(n-2)b} \Big ) \, \sqrt{2a} \, (1-\lambda^2) \, \text{\rm scal}(T).
\end{align*} 
Using Lemma \ref{pic.implies.pic1.with.error.term}, we obtain 
\begin{align*} 
Z &\geq -\frac{1-\lambda^2}{n-4} \, (\text{\rm Ric}(T)_{11}+\text{\rm Ric}(T)_{22}-2T_{1212}) \\ 
&+ \frac{b_{\text{\rm max}}-b}{1+(n-2)b} \, [(1+\lambda^2) \, \text{\rm Ric}(T)_{11}+(1+\lambda^2) \, \text{\rm Ric}(T)_{22} \\ 
&\hspace{30mm} +2 \, \text{\rm Ric}(T)_{33}+2\lambda^2 \, \text{\rm Ric}(T)_{44}] \\ 
&+ \frac{4}{n} \, \Big ( \frac{a_{\text{\rm max}}-a}{1+2(n-1)a} - \frac{b_{\text{\rm max}}-b}{1+(n-2)b} \Big ) \, (1+\lambda^2) \, \text{\rm scal}(T) \\ 
&+ \frac{1+(n-2)b_{\text{\rm max}}}{1+(n-2)b} \, \sqrt{2a} \, (1-\lambda^2) \, (\text{\rm Ric}(T)_{11}+\text{\rm Ric}(T)_{22}) \\ 
&+ \frac{2}{n} \, \Big ( \frac{2(n-1)(a_{\text{\rm max}}-a)}{1+2(n-1)a} - \frac{(n-2)(b_{\text{\rm max}}-b)}{1+(n-2)b} \Big ) \, \sqrt{2a} \, (1-\lambda^2) \, \text{\rm scal}(T) \\ 
&=: RHS. 
\end{align*} 
We claim that $RHS \geq 0$ for all $\lambda \in [0,1]$. Since $RHS$ is a linear function of $\lambda^2$, it suffices to examine the endpoints of the interval:

\textit{Case 1:} Suppose first that $\lambda=1$. In this case,  
\begin{align*} 
RHS &= 2 \, \frac{b_{\text{\rm max}}-b}{1+(n-2)b} \, [\text{\rm Ric}(T)_{11}+\text{\rm Ric}(T)_{22}+\text{\rm Ric}(T)_{33}+\text{\rm Ric}(T)_{44}] \\ 
&+ \frac{8}{n} \, \Big ( \frac{a_{\text{\rm max}}-a}{1+2(n-1)a} - \frac{b_{\text{\rm max}}-b}{1+(n-2)b} \Big ) \, \text{\rm scal}(T). 
\end{align*} 
By Lemma \ref{Ric.four.positive}, $\text{\rm Ric}(T)_{11}+\text{\rm Ric}(T)_{22}+\text{\rm Ric}(T)_{33}+\text{\rm Ric}(T)_{44} \geq 0$. Using Lemma \ref{parameters.2}, we conclude that $RHS \geq 0$.

\textit{Case 2:} Suppose next that $\lambda=0$. Recall that 
\[\frac{1}{n-4} \leq \frac{1+(n-2)b_{\text{\rm max}}}{1+(n-2)b} \, \sqrt{2a} \, \frac{n^2-7n+14}{n^2-5n+4}\] 
by Lemma \ref{parameters.2}. Using the inequalities 
\[\text{\rm Ric}(T)_{11}+\text{\rm Ric}(T)_{22}-2T_{1212} \geq 0\] 
and 
\begin{align*} 
&\text{\rm Ric}(T)_{11}+\text{\rm Ric}(T)_{22}+2 \, \text{\rm Ric}(T)_{33} \\ 
&= (\text{\rm Ric}(T)_{11}+\text{\rm Ric}(T)_{33}) + (\text{\rm Ric}(T)_{22}+\text{\rm Ric}(T)_{33}) \geq -\frac{4\gamma_{\text{\rm max}}}{n} \, \text{\rm scal}(T), 
\end{align*}
we obtain 
\begin{align*} 
RHS &\geq -\frac{1+(n-2)b_{\text{\rm max}}}{1+(n-2)b} \, \sqrt{2a} \, \frac{n^2-7n+14}{n^2-5n+4} \, (\text{\rm Ric}(T)_{11}+\text{\rm Ric}(T)_{22}-2T_{1212}) \\ 
&- \frac{4\gamma_{\text{\rm max}}}{n} \, \frac{b_{\text{\rm max}}-b}{1+(n-2)b} \, \text{\rm scal}(T) \\ 
&+ \frac{4}{n} \, \Big ( \frac{a_{\text{\rm max}}-a}{1+2(n-1)a} - \frac{b_{\text{\rm max}}-b}{1+(n-2)b} \Big ) \, \text{\rm scal}(T) \\ 
&+ \frac{1+(n-2)b_{\text{\rm max}}}{1+(n-2)b} \, \sqrt{2a} \, (\text{\rm Ric}(T)_{11}+\text{\rm Ric}(T)_{22}) \\ 
&+ \frac{2}{n} \, \Big ( \frac{2(n-1)(a_{\text{\rm max}}-a)}{1+2(n-1)a} - \frac{(n-2)(b_{\text{\rm max}}-b)}{1+(n-2)b} \Big ) \, \sqrt{2a} \, \text{\rm scal}(T). 
\end{align*}
Since $T \in \text{\rm PIC}$, we can estimate the term $T_{1212}$ in terms of $\text{\rm Ric}(T)_{11}+\text{\rm Ric}(T)_{22}$ and $\text{\rm scal}(T)$: 
\begin{align*} 
0 &\leq \sum_{3 \leq p,q \leq n, \, p \neq q} (T_{1212}+T_{1p1p}+T_{2q2q}+T_{pqpq}) \\ 
&= (n^2-7n+14) \, T_{1212} + (n-5) (\text{\rm Ric}(T)_{11}+\text{\rm Ric}(T)_{22}) + \text{\rm scal}(T). 
\end{align*} 
This finally gives 
\begin{align*} 
RHS &\geq -\frac{1+(n-2)b_{\text{\rm max}}}{1+(n-2)b} \, \frac{2\sqrt{2a}}{n^2-5n+4} \, \text{\rm scal}(T) - \frac{4\gamma_{\text{\rm max}}}{n} \, \frac{b_{\text{\rm max}}-b}{1+(n-2)b} \, \text{\rm scal}(T) \\ 
&+ \frac{4}{n} \, \Big ( \frac{a_{\text{\rm max}}-a}{1+2(n-1)a} - \frac{b_{\text{\rm max}}-b}{1+(n-2)b} \Big ) \, \text{\rm scal}(T) \\ 
&+ \frac{2}{n} \, \Big ( \frac{2(n-1)(a_{\text{\rm max}}-a)}{1+2(n-1)a} - \frac{(n-2)(b_{\text{\rm max}}-b)}{1+(n-2)b} \Big ) \, \sqrt{2a} \, \text{\rm scal}(T). 
\end{align*}
Using Lemma \ref{parameters.2}, we obtain $RHS \geq 0$. This completes the proof of Proposition \ref{pinching.cones.can.be.joined.together}. \\

In the remainder of this section, we show that the cone $\tilde{\mathcal{C}}(b)$ is preserved by the Hamilton ODE for $0 < b \leq \tilde{b}_{\text{\rm max}}$. We first state an elementary lemma:

\begin{lemma}
\label{parameters.3} 
Assume that $0 < b \leq \tilde{b}_{\text{\rm max}}$, and let $a$ be defined by $2a=2b+(n-2)b^2$. Then 
\[nb^2(1-2b)-2(a-b)(1-2b+nb^2) \geq 0,\] 
\[n^2b^2-2(n-1)(a-b)(1-2b) \geq 0.\] 
Moreover, if we put 
\[\zeta := \frac{1+2(n-1)a}{1+2(n-1)a_{\text{\rm max}}} \, \frac{1+(n-2)b_{\text{\rm max}}}{1+(n-2)b} \, (1+\gamma_{\text{\rm max}}),\] 
then $\zeta \leq 1$ and $1 + (n-2)(1-\zeta) \geq 2\zeta^2 \, \frac{n^2-2n+2}{(n-2)^2}$.
\end{lemma}
 
\textbf{Proof.}
To prove the first and second statement, we write 
\[nb^2(1-2b) - 2(a-b)(1-2b+nb^2) = b^2(2-4b-n(n-2)b^2) \geq 0\] 
and 
\[n^2b^2 - 2(n-1)(a-b)(1-2b) = b^2 (3n-2+2(n-1)(n-2)b) \geq 0.\] 
To verify the third and fourth statement, we observe that $\frac{1+2(n-1)a}{1+(n-2)b} \leq \frac{1+2(n-1)\tilde{a}_{\text{\rm max}}}{1+(n-2)\tilde{b}_{\text{\rm max}}}$, hence 
\[\zeta \leq \zeta_{\text{\rm max}} := \frac{1+2(n-1)\tilde{a}_{\text{\rm max}}}{1+2(n-1)a_{\text{\rm max}}} \, \frac{1+(n-2)b_{\text{\rm max}}}{1+(n-2)\tilde{b}_{\text{\rm max}}} \, (1+\gamma_{\text{\rm max}}).\] 
It is elementary to check that $\zeta_{\text{\rm max}} \leq 1$ and $1+(n-2)(1-\zeta_{\text{\rm max}}) \geq 2\zeta_{\text{\rm max}}^2 \, \frac{n^2-2n+2}{(n-2)^2}$. This implies $\zeta \leq 1$ and $1 + (n-2)(1-\zeta) \geq 2\zeta^2 \, \frac{n^2-2n+2}{(n-2)^2}$. \\

\begin{lemma}
\label{Ricci.bound.coming.from.Cbmax}
Suppose that $S \in \tilde{\mathcal{E}}(b)$, and 
\[\zeta := \frac{1+2(n-1)a}{1+2(n-1)a_{\text{\rm max}}} \, \frac{1+(n-2)b_{\text{\rm max}}}{1+(n-2)b} \, (1+\gamma_{\text{\rm max}}) \leq 1.\] 
Then $\text{\rm Ric}(S)_{11}+\text{\rm Ric}(S)_{22} \geq \frac{2(1-\zeta)}{n} \, \text{\rm scal}(S)$. 
\end{lemma}
 
\textbf{Proof.} 
Let $T = \ell_{a_{\text{\rm max}},b_{\text{\rm max}}}^{-1}(\ell_{a,b}(S)) \in \mathcal{E}(b_{\text{\rm max}})$. Recall that $\tracefreeRic(S) = \frac{1+(n-2)b_{\text{\rm max}}}{1+(n-2)b} \, \tracefreeRic(T)$ and $\text{\rm scal}(S) = \frac{1+2(n-1)a_{\text{\rm max}}}{1+2(n-1)a} \, \text{\rm scal}(T)$. Using the inequality $\tracefreeRic(T)_{11}+\tracefreeRic(T)_{22} + \frac{2(1+\gamma_{\text{\rm max}})}{n} \, \text{\rm scal}(T) \geq 0$, we obtain $\tracefreeRic(S)_{11}+\tracefreeRic(S)_{22} + \frac{2\zeta}{n} \, \text{\rm scal}(S) \geq 0$, which implies the claim. \\ 

After these preparations, we now prove the main result of this section: 
 
\begin{theorem}
\label{tildeCb.is.preserved}
For each $0 < b \leq \tilde{b}_{\text{\rm max}}$, the cone $\tilde{\mathcal{C}}(b)$ is transversally invariant under the Hamilton ODE $\frac{d}{dt} R = Q(R)$.
\end{theorem}

\textbf{Proof.} It suffices to show that $\tilde{\mathcal{E}}(b)$ is transversally invariant under the ODE $\frac{d}{dt} S = Q(S)+D_{a,b}(S)$. Suppose that $S \in \tilde{\mathcal{E}}(b) \setminus \{0\}$. Then 
\begin{align*} 
Z &:= S_{1313}+\lambda^2 S_{1414}+S_{2323}+\lambda^2 S_{2424}-2\lambda S_{1234} \\ 
&+ \sqrt{2a} \, (1-\lambda^2) \, (\text{\rm Ric}(S)_{11}+\text{\rm Ric}(S)_{22}) \geq 0 
\end{align*}
for every orthonormal four-frame $\{e_1,e_2,e_3,e_4\}$ and every $\lambda \in [0,1]$. Moreover, we assume that $Z=0$ for some orthonormal four-frame $\{e_1,e_2,e_3,e_4\}$ and some $\lambda \in [0,1]$. We will show that $\frac{d}{dt} Z > 0$ for this particular four-frame $\{e_1,e_2,e_3,e_4\}$ and this particular $\lambda \in [0,1]$. We distinguish two cases:

\textit{Case 1:} Suppose that $\lambda \in [0,1)$. Using the formula for $D_{a,b}(S)$, we obtain 
\begin{align*} 
\frac{d}{dt} S &= Q(S) + 2a \, \text{\rm Ric}(S) \owedge \text{\rm Ric}(S) + 2b^2 \, \tracefreeRic(S)^2 \owedge \text{\rm id} \\ 
&+ \frac{nb^2(1-2b)-2(a-b)(1-2b+nb^2)}{n(1+2(n-1)a)} \, |\tracefreeRic(S)|^2 \, \text{\rm id} \owedge \text{\rm id} 
\end{align*} 
and 
\begin{align*} 
\frac{d}{dt} \text{\rm Ric}(S) 
&= 2 \, S * \text{\rm Ric}(S) - 4b \, \tracefreeRic(S)^2 \\ 
&+ \frac{4(n-2)a}{n} \, \text{\rm scal}(S) \, \text{\rm Ric}(S) + \frac{4a}{n^2} \, \text{\rm scal}(S)^2 \, \text{\rm id} \\ 
&+ 2 \, \frac{n^2b^2-2(n-1)(a-b)(1-2b)}{n(1+2(n-1)a)} \, |\tracefreeRic(S)|^2 \, \text{\rm id},
\end{align*} 
where $(S * H)_{ik} := \sum_{p,q=1}^n S_{ipkq} H_{pq}$. This implies  
\begin{align*} 
\frac{d}{dt} Z 
&\geq Q(S)_{1313}+\lambda^2 Q(S)_{1414}+Q(S)_{2323}+\lambda^2 Q(S)_{2424}-2\lambda Q(S)_{1234} \\ 
&+ 2a \, (\text{\rm Ric}(S) \owedge \text{\rm Ric}(S))_{1313}+2a\lambda^2 \, (\text{\rm Ric}(S) \owedge \text{\rm Ric}(S))_{1414} \\ 
&+ 2a \, (\text{\rm Ric}(S) \owedge \text{\rm Ric}(S))_{2323}+2a\lambda^2 \, (\text{\rm Ric}(S) \owedge \text{\rm Ric}(S))_{2424} \\ 
&- 4a\lambda \, (\text{\rm Ric}(S) \owedge \text{\rm Ric}(S))_{1234} \\ 
&+ 2\sqrt{2a} \, (1-\lambda^2) \, ((S * \text{\rm Ric}(S))_{11}+(S * \text{\rm Ric}(S))_{22}) \\ 
&+ \frac{4(n-2)a}{n} \, \sqrt{2a} \, (1-\lambda^2) \, \text{\rm scal}(S) \, (\text{\rm Ric}(S)_{11}+\text{\rm Ric}(S)_{22}) \\ 
&- 4b \, \sqrt{2a} \, (1-\lambda^2) \, ((\tracefreeRic(S)^2)_{11}+(\tracefreeRic(S)^2)_{22}) \\ 
&+ \frac{8a}{n^2} \, \sqrt{2a} \, (1-\lambda^2) \, \text{\rm scal}(S)^2.
\end{align*} 
Applying Proposition \ref{interpolation.between.pic.and.pic1} with $H = \sqrt{2a} \, \text{\rm Ric}(S)$, we obtain 
\begin{align*} 
&Q(S)_{1313}+\lambda^2 Q(S)_{1414}+Q(S)_{2323}+\lambda^2 Q(S)_{2424}-2\lambda Q(S)_{1234} \\ 
&+ 2a \, (\text{\rm Ric}(S) \owedge \text{\rm Ric}(S))_{1313}+2a\lambda^2 \, (\text{\rm Ric}(S) \owedge \text{\rm Ric}(S))_{1414} \\ 
&+ 2a \, (\text{\rm Ric}(S) \owedge \text{\rm Ric}(S))_{2323}+2a\lambda^2 \, (\text{\rm Ric}(S) \owedge \text{\rm Ric}(S))_{2424} \\ 
&- 4a\lambda \, (\text{\rm Ric}(S) \owedge \text{\rm Ric}(S))_{1234} \\ 
&+ 2\sqrt{2a} \, (1-\lambda^2) \, ((S * \text{\rm Ric}(S))_{11}+(S * \text{\rm Ric}(S))_{22}) \geq 0. 
\end{align*} 
Since $S \in \tilde{\mathcal{E}}(b) \subset \text{\rm PIC}$, Lemma \ref{estimate.for.largest.ricci.eigenvalue} implies that the largest eigenvalue of $\text{\rm Ric}(S)$ is bounded from above by $\frac{1}{2} \, \text{\rm scal}(S)$. Moreover, by Lemma \ref{Ricci.bound.coming.from.Cbmax}, the sum of the two smallest eigenvalues of $\text{\rm Ric}(S)$ is bounded from below by $\frac{2(1-\zeta)}{n} \, \text{\rm scal}(S)$. Applying  Lemma \ref{algebraic.fact} with $H = \text{\rm Ric}(S)$ and $\rho = \frac{b}{a}$ gives 
\begin{align*} 
&\frac{(n-2)a}{n} \, \text{\rm scal}(S) \, (\text{\rm Ric}(S)_{11}+\text{\rm Ric}(S)_{22}) - b \, ((\tracefreeRic(S)^2)_{11}+(\tracefreeRic(S)^2)_{22}) \\ 
&\geq \frac{2}{n^2} \, \Big ( a(n-2)(1-\zeta) - 2b\zeta^2 \, \frac{n^2-2n+2}{(n-2)^2} \Big ) \, \text{\rm scal}(S)^2. 
\end{align*} 
Putting these facts together, we conclude that 
\[\frac{d}{dt} Z \geq \frac{8}{n^2} \, \Big ( a(1+(n-2)(1-\zeta)) - 2b\zeta^2 \, \frac{n^2-2n+2}{(n-2)^2} \Big ) \, \sqrt{2a} \, (1-\lambda^2) \, \text{\rm scal}(S)^2.\]
On the other hand, Lemma \ref{parameters.3} implies 
\[a(1+(n-2)(1-\zeta)) \geq 2a\zeta^2 \, \frac{n^2-2n+2}{(n-2)^2} > 2b\zeta^2 \, \frac{n^2-2n+2}{(n-2)^2}.\] 
Since $\lambda \in [0,1)$, it follows that $\frac{d}{dt} Z > 0$. 

\textit{Case 2:} Finally, suppose that $\lambda=1$. In this case, 
\begin{align*} 
\frac{d}{dt} Z 
&= Q(S)_{1313}+Q(S)_{1414}+Q(S)_{2323}+Q(S)_{2424}-2Q(S)_{1234} \\ 
&+ 2a \, (\text{\rm Ric}(S) \owedge \text{\rm Ric}(S))_{1313}+2a \, (\text{\rm Ric}(S) \owedge \text{\rm Ric}(S))_{1414} \\ 
&+ 2a \, (\text{\rm Ric}(S) \owedge \text{\rm Ric}(S))_{2323}+2a \, (\text{\rm Ric}(S) \owedge \text{\rm Ric}(S))_{2424} \\ 
&- 4a \, (\text{\rm Ric}(S) \owedge \text{\rm Ric}(S))_{1234}.  
\end{align*} 
Since $S \in \text{\rm PIC}$ and $Z = S_{1313}+S_{1414}+S_{2323}+S_{2424}-2S_{1234} = 0$, we obtain 
\[Q(S)_{1313}+Q(S)_{1414}+Q(S)_{2323}+Q(S)_{2424}-2Q(S)_{1234} \geq 0\] 
for this particular four-frame $\{e_1,e_2,e_3,e_4\}$. By Lemma \ref{Ricci.bound.coming.from.Cbmax}, the sum of the two smallest eigenvalues of $\text{\rm Ric}(S)$ is strictly positive. This implies 
\begin{align*} 
&(\text{\rm Ric}(S) \owedge \text{\rm Ric}(S))_{1313}+(\text{\rm Ric}(S) \owedge \text{\rm Ric}(S))_{1414} \\ 
&+ (\text{\rm Ric}(S) \owedge \text{\rm Ric}(S))_{2323}+(\text{\rm Ric}(S) \owedge \text{\rm Ric}(S))_{2424} \\ 
&- 2 \, (\text{\rm Ric}(S) \owedge \text{\rm Ric}(S))_{1234} > 0 
\end{align*} 
by Lemma \ref{wedge.product}. Putting these facts together, we conclude that $\frac{d}{dt} Z > 0$. This completes the proof of Theorem \ref{tildeCb.is.preserved}.

\section{Pinching towards PIC2 and proof of Theorem \ref{main.theorem}} 

\label{proof.of.main.theorem}

\begin{theorem}
\label{pinching.towards.pic1}
Assume that $n \geq 12$. Let $\mathcal{K}$ be a compact set of algebraic curvature tensors in dimension $n$ which is contained in the interior of the PIC cone, and let $T>0$ be given. Then there exist a small positive real number $\theta$, a large positive real number $N$, a smooth, increasing, and concave function $g>0$ satisfying $\lim_{s \to \infty} \frac{g(s)}{s} = 0$, and a continuous family of closed, convex, $O(n)$-invariant sets $\{\mathcal{G}_t: t \in [0,T]\}$ such that the family $\{\mathcal{G}_t: t \in [0,T]\}$ is invariant under the Hamilton ODE $\frac{d}{dt} R = Q(R)$; $\mathcal{K} \subset \mathcal{G}_0$; and 
\begin{align*} 
\mathcal{G}_t 
&\subset \{R: R - \theta \, \text{\rm scal} \, \text{\rm id} \owedge \text{\rm id} \in \text{\rm PIC}\} \\ 
&\cap \{R: \text{\rm Ric}_{11}+\text{\rm Ric}_{22} - \theta \, \text{\rm scal} + N \geq 0\} \\ 
&\cap \{R: R + g(\text{\rm scal}) \, \text{\rm id} \owedge \text{\rm id} \in \text{\rm PIC1}\}
\end{align*} 
for all $t \in [0,T]$.
\end{theorem}

\textbf{Proof.} 
By Proposition \ref{pinching.cones.can.be.joined.together}, we have $\tilde{\mathcal{C}}(\tilde{b}_{\text{\rm max}}) = \mathcal{C}(b_{\text{\rm max}})$. For each $b \in (0,b_{\text{\rm max}}+\tilde{b}_{\text{\rm max}})$, we define 
\[\hat{\mathcal{C}}(b) := \begin{cases} \mathcal{C}(b) & \text{\rm for $b \in (0,b_{\text{\rm max}}]$} \\ \tilde{\mathcal{C}}(b_{\text{\rm max}}+\tilde{b}_{\text{\rm max}}-b) & \text{\rm for $b \in (b_{\text{\rm max}},b_{\text{\rm max}}+\tilde{b}_{\text{\rm max}})$.} \end{cases}\]
Therefore, $\hat{\mathcal{C}}(b)$ is a family of closed, convex, $O(n)$-invariant cones which depends continuously on the parameter $b \in (0,b_{\text{\rm max}}+\tilde{b}_{\text{\rm max}})$. Moreover, for each $b \in (0,b_{\text{\rm max}}+\tilde{b}_{\text{\rm max}})$, the cone $\hat{\mathcal{C}}(b)$ is transversally invariant under the Hamilton ODE. 

By Theorem \ref{rough.pinching} and Corollary \ref{uniformly.pic}, there exists a family of closed, convex, $O(n)$-invariant sets $\{\mathcal{G}_t^{(0)}: t \in [0,T]\}$ such that the family $\{\mathcal{G}_t^{(0)}: t \in [0,T]\}$ is invariant under the Hamilton ODE $\frac{d}{dt} R = Q(R)$; $\mathcal{K} \subset \mathcal{G}_0^{(0)}$; and 
\begin{align*}
\mathcal{G}_t^{(0)} 
&\subset \{R: R - \theta \, \text{\rm scal} \, \text{\rm id} \owedge \text{\rm id} \in \text{\rm PIC}\} \\ 
&\cap \{R: \text{\rm Ric}_{11}+\text{\rm Ric}_{22} - \theta \, \text{\rm scal} + N \geq 0\}
\end{align*} 
for all $t \in [0,T]$. Consequently, we can find a small positive number $\beta_0$ (depending on $\theta$) such that 
\[\mathcal{G}_t^{(0)} \subset \{R: R + N \, \text{\rm id} \owedge \text{\rm id} \in \hat{\mathcal{C}}(\beta_0)\}\] 
for all $t \in [0,T]$.

Arguing as in Theorem 4.1 in \cite{Bohm-Wilking} (see also Proposition 16 in \cite{Brendle-Schoen}), we can construct an increasing sequence $\beta_0 < \beta_1 < \beta_2 < \hdots$ with $\lim_{j \to \infty} \beta_j = b_{\text{\rm max}}+\tilde{b}_{\text{\rm max}}$ and a sequence of families of sets $\{\mathcal{G}_t^{(j)}: t \in [0,T]\}$ with the following properties: 
\begin{itemize}
\item For $j=0$, the family of sets $\{\mathcal{G}_t^{(0)}: t \in [0,T]\}$ coincides with the one constructed in Theorem \ref{rough.pinching}. 
\item For each $j \geq 1$, the family of sets $\{\mathcal{G}_t^{(j)}: t \in [0,T]\}$ is defined by 
\[\mathcal{G}_t^{(j)} = \{R \in \mathcal{G}_t^{(j-1)}: R + 2^j N \, \text{\rm id} \owedge \text{\rm id} \in \hat{\mathcal{C}}(\beta_j)\}\] 
for $t \in [0,T]$. 
\item For each $j \geq 1$, the family $\{\mathcal{G}_t^{(j)}: t \in [0,T]\}$ is invariant under the Hamilton ODE. 
\item For each $j \geq 1$, we have $\mathcal{K} \subset \mathcal{G}_0^{(j)}$.
\end{itemize}
Since $\lim_{j \to \infty} \beta_j = b_{\text{\rm max}}+\tilde{b}_{\text{\rm max}}$, we can find a sequence of positive numbers $\varepsilon_j$ such that $\lim_{j \to \infty} \varepsilon_j = 0$ and $\hat{\mathcal{C}}(\beta_j) \subset \{R: R + \varepsilon_j \, \text{\rm scal} \, \text{\rm id} \owedge \text{\rm id} \in \text{\rm PIC1}\}$. This implies 
\[\mathcal{G}_t^{(j)} \subset \{R: R + \varepsilon_j \, \text{\rm scal} \, \text{\rm id} \owedge \text{\rm id} + 2^j N \, (1+2n(n-1)\varepsilon_j) \, \text{\rm id} \owedge \text{\rm id} \in \text{\rm PIC1}\}\] 
for all $t \in [0,T]$.

We now define $\mathcal{G}_t = \bigcap_{j \in \mathbb{N}} \mathcal{G}_t^{(j)}$ for each $t \in [0,T]$. Then $\{\mathcal{G}_t: t \in [0,T]\}$ is a family of closed, convex, $O(n)$-invariant sets satisfying $\mathcal{K} \subset \mathcal{G}_0$. Moreover, the family of sets $\{\mathcal{G}_t: t \in [0,T]\}$ is invariant under the Hamilton ODE. Furthermore,  
\[\mathcal{G}_t \subset \{R \in \mathcal{G}_t^{(0)}: R + h(\text{\rm scal}) \, \text{\rm id} \owedge \text{\rm id} \in \text{\rm PIC1}\}\] 
for all $t \in [0,T]$, where $h(s) := \inf_{j \in \mathbb{N}} (\varepsilon_j s + 2^j N \, (1+2n(n-1)\varepsilon_j))$. Clearly, $h$ is an increasing concave function satisfying $\lim_{s \to \infty} \frac{h(s)}{s} = 0$. Finally, we put $g(s) := \int_0^1 \chi(r) \, h(s+r) \, dr$, where $\chi$ is a nonnegative smooth function which is supported in $[0,1]$ and satisfies $\int_0^1 \chi(r) \, dr = 1$. Clearly, $g(s) \geq h(s)$. It is easy to see that $g$ has all the required properties. \\

\begin{proposition}
\label{positivity.for.lambda.mu.close.to.1}
There exists a small positive constant $\kappa>0$ with the following property. If $R \in \mathcal{G}_t$, then
\[R_{1313}+\lambda^2 R_{1414}+\mu^2 R_{2323}+\lambda^2\mu^2 R_{2424} - 2\lambda\mu R_{1234} \geq 0\] 
for every orthonormal four-frame $\{e_1,e_2,e_3,e_4\}$ and all $\lambda,\mu \in [0,1]$ satisfying $\lambda^2(1-\mu^2)^2+\mu^2(1-\lambda^2)^2+(1-\lambda^2)^2(1-\mu^2)^2 \leq \kappa$.
\end{proposition} 

\textbf{Proof.}
This follows immediately from the fact that $R - \theta \, \text{\rm scal} \, \text{\rm id} \owedge \text{\rm id} \in \text{\rm PIC}$ for all $R \in \mathcal{G}_t$. \\

We now turn to the proof of Theorem \ref{main.theorem}. To any curvature tensor $R \in \mathcal{G}_t$, we associate a curvature tensor $S \in \text{\rm PIC1}$ by $S = R + 2  \, g(\text{\rm scal}) \, \text{\rm id} \owedge \text{\rm id}$. Note that $|Q(R)-Q(S)|$ is bounded by a large multiple of $ (\text{\rm scal}+g(\text{\rm scal})) \, g(\text{\rm scal})$. In particular,
\[Q(R) - Q(S) + L(n) \, (\text{\rm scal}+g(\text{\rm scal})) \, g(\text{\rm scal}) \, \text{\rm id} \owedge \text{\rm id} \in \text{\rm PIC2},\] 
where $L(n)$ is a large constant that depends only on $n$. \\

\begin{definition}
Let $g$ be chosen as in Theorem \ref{pinching.towards.pic1}, let $\kappa$ be chosen as in Proposition \ref{positivity.for.lambda.mu.close.to.1}, and let $L(n)$ be chosen as above. Let $f(s) := \Lambda \, \sqrt{(s+g(s)) \, g(s)}$, where $\Lambda \geq 16 \sqrt{\frac{L(n)}{\kappa}}$. For each $t \in [0,T]$, we denote by $\mathcal{F}_t$ denote the set of all algebraic curvature tensors $R \in \mathcal{G}_t$ such that 
\begin{align*} 
&S_{1313}+\lambda^2 S_{1414}+\mu^2 S_{2323}+\lambda^2\mu^2 S_{2424} - 2\lambda\mu S_{1234} \\ 
&+ (1-\lambda^2)(1-\mu^2) \, f(\text{\rm scal}) \geq 0 
\end{align*}
for all $\lambda,\mu \in [0,1]$, where $S$ is defined by $S = R + 2 \, g(\text{\rm scal}) \, \text{\rm id} \owedge \text{\rm id}$.
\end{definition}

Clearly, $f$ is a smooth, increasing, and concave function satisfying $\lim_{s \to \infty} \frac{f(s)}{s} = 0$. Since $f$ and $g$ are concave functions, the set $\mathcal{F}_t$ is convex for each $t \in [0,T]$. 

\begin{theorem}
\label{pinching.towards.pic2} 
If $\Lambda \geq 16 \sqrt{\frac{L(n)}{\kappa}}$, then the family of sets $\{\mathcal{F}_t: t \in [0,T]\}$ is invariant under the Hamilton ODE $\frac{d}{dt} R = Q(R)$. Moreover, $\mathcal{F}_t \subset \{R: R + 2 \, (f(\text{\rm scal})+g(\text{\rm scal})) \, \text{\rm id} \owedge \text{\rm id} \in \text{\rm PIC2}\}$. Finally, if we choose $\Lambda$ sufficiently large, then $\mathcal{K} \subset \mathcal{F}_0$.
\end{theorem}

\textbf{Proof.} 
We first show that the family of sets $\{\mathcal{F}_t: t \in [0,T]\}$ is preserved by the Hamilton ODE $\frac{d}{dt} R = Q(R)$. To prove this, we assume that $R \in \mathcal{F}_t$, $S = R + 2  \, g(\text{\rm scal}) \, \text{\rm id} \owedge \text{\rm id}$, and  
\begin{align*} 
&S_{1313}+\lambda^2 S_{1414}+\mu^2 S_{2323}+\lambda^2\mu^2 S_{2424} - 2\lambda\mu S_{1234} \\ 
&+ (1-\lambda^2)(1-\mu^2) \, f(\text{\rm scal}) = 0 
\end{align*}
for one particular orthonormal four-frame $\{e_1,e_2,e_3,e_4\}$ and one particular pair $\lambda,\mu \in [0,1]$. Since $S$ lies in the interior of the PIC1 cone, it follows that $\lambda,\mu \in [0,1)$. Moreover, 
\begin{align*} 
&R_{1313}+\lambda^2 R_{1414}+\mu^2 R_{2323}+\lambda^2\mu^2 R_{2424} - 2\lambda\mu R_{1234} \\ 
&< S_{1313}+\lambda^2 S_{1414}+\mu^2 S_{2323}+\lambda^2\mu^2 S_{2424} - 2\lambda\mu S_{1234} \leq 0. 
\end{align*}
Hence, Proposition \ref{positivity.for.lambda.mu.close.to.1} implies that 
\[\lambda^2(1-\mu^2)^2+\mu^2(1-\lambda^2)^2+(1-\lambda^2)^2(1-\mu^2)^2 > \kappa.\]
Since $\lambda,\mu \in [0,1)$, the first order conditions corresponding to variations of $\lambda$ and $\mu$ give 
\begin{align*} 
\lambda(1-\mu^2) \, f(\text{\rm scal}) 
&= \lambda S_{1414}+\lambda\mu^2 S_{2424}-\mu S_{1234} \\ 
&= (\lambda S_{1414}+\mu S_{1423}) - \mu (S_{1324}-\lambda \mu S_{2424}) 
\end{align*} 
and
\begin{align*} 
(1-\lambda^2)\mu \, f(\text{\rm scal}) 
&= \mu S_{2323}+\lambda^2\mu S_{2424}-\lambda S_{1234} \\ 
&= (\lambda S_{1423}+\mu S_{2323}) - \lambda (S_{1324}-\lambda \mu S_{2424}). 
\end{align*} 
Moreover, 
\begin{align*} 
(1-\lambda^2)(1-\mu^2) \, f(\text{\rm scal}) 
&= -(S_{1313}-\lambda\mu S_{1324}) - \lambda(\lambda S_{1414} + \mu S_{1423}) \\ 
&- \mu(\lambda S_{1423} + \mu S_{2323}) + \lambda\mu (S_{1324} - \lambda\mu S_{2424}).
\end{align*}
Consequently, 
\begin{align*} 
&[\lambda^2(1-\mu^2)^2+\mu^2(1-\lambda^2)^2+(1-\lambda^2)^2(1-\mu^2)^2] \, f(\text{\rm scal})^2 \\ 
&\leq 2 \,  (\lambda S_{1414}+\mu S_{1423})^2 + 2\mu^2 (S_{1324}-\lambda \mu S_{2424})^2 \\ 
&+ 2 \, (\lambda S_{1423}+\mu S_{2323})^2 + 2\lambda^2 (S_{1324}-\lambda \mu S_{2424})^2 \\ 
&+  4 \, (S_{1313}-\lambda\mu S_{1324})^2 + 4\lambda^2(\lambda S_{1414} + \mu S_{1423})^2 \\ 
&+ 4\mu^2(\lambda S_{1423} + \mu S_{2323})^2 + 4\lambda^2\mu^2 (S_{1324} - \lambda\mu S_{2424})^2 \\ 
&\leq 16 \sum_{p,q=1}^n (S_{13pq}-\lambda\mu S_{24pq})^2 + 16 \sum_{p,q=1}^n (\lambda S_{14pq}+\mu S_{23pq})^2. 
\end{align*}
Putting these facts together, we obtain 
\[\sum_{p,q=1}^n (S_{13pq}-\lambda\mu S_{24pq})^2 + \sum_{p,q=1}^n (\lambda S_{14pq}+\mu S_{23pq})^2 \geq \frac{\kappa}{16} \, f(\text{\rm scal})^2.\] 
On the other hand, Proposition 7.26 in \cite{Brendle-book} implies  
\[(S^\#)_{1313}+\lambda^2 (S^\#)_{1414}+\mu^2 (S^\#)_{2323}+\lambda^2\mu^2 (S^\#)_{2424} + 2\lambda\mu (S^\#)_{1342} + 2\lambda\mu (S^\#)_{1423} \geq 0.\] 
Adding these inequalities gives 
\[Q(S)_{1313}+\lambda^2 Q(S)_{1414}+\mu^2 Q(S)_{2323}+\lambda^2\mu^2 Q(S)_{2424} - 2\lambda\mu Q(S)_{1234} \geq \frac{\kappa}{16} \, f(\text{\rm scal})^2.\] 
Since $Q(R) - Q(S) + L(n) \, (\text{\rm scal}+g(\text{\rm scal})) \, g(\text{\rm scal}) \, \text{\rm id} \owedge \text{\rm id} \in \text{\rm PIC2}$, it follows that 
\begin{align*} 
&Q(R)_{1313}+\lambda^2 Q(R)_{1414}+\mu^2 Q(R)_{2323}+\lambda^2\mu^2 Q(R)_{2424} - 2\lambda\mu Q(R)_{1234} \\ 
&\geq \frac{\kappa}{16} \, f(\text{\rm scal})^2 - 8L(n) \, (\text{\rm scal}+g(\text{\rm scal})) \, g(\text{\rm scal}). 
\end{align*}
Since $f$ and $g$ are monotone increasing, we conclude that 
\begin{align*} 
&\frac{d}{dt} \big ( S_{1313}+\lambda^2 S_{1414}+\mu^2 S_{2323}+\lambda^2\mu^2 S_{2424} - 2\lambda\mu S_{1234} \\ 
&\hspace{5mm} + (1-\lambda^2)(1-\mu^2) \, f(\text{\rm scal}) \big ) \\ 
&\geq \frac{\kappa}{16} \, f(\text{\rm scal})^2 - 8L(n) \, (\text{\rm scal}+g(\text{\rm scal})) \, g(\text{\rm scal}), 
\end{align*} 
and the right hand side is positive in view of our choice of $f$. This shows that the family of sets $\{\mathcal{F}_t: t \in [0,T]\}$ is invariant under the Hamilton ODE. 

If $R \in \mathcal{F}_t$ and $S = R + 2 \, g(\text{\rm scal}) \, \text{\rm id} \owedge \text{\rm id}$, then $S + 2 \, f(\text{\rm scal}) \, \text{\rm id} \owedge \text{\rm id} \in \text{\rm PIC2}$ (cf. Proposition 7.18 in \cite{Brendle-book}), and consequently $R + 2 \, (f(\text{\rm scal})+g(\text{\rm scal})) \, \text{\rm id} \owedge \text{\rm id} \in \text{\rm PIC2}$.

Finally, if we choose $\Lambda$ sufficiently large, we can arrange that the following holds: if $R \in \mathcal{K}$ and $S = R + 2 \, g(\text{\rm scal}) \, \text{\rm id} \owedge \text{\rm id}$, then 
\begin{align*} 
&S_{1313}+\lambda^2 S_{1414}+\mu^2 S_{2323}+\lambda^2\mu^2 S_{2424} - 2\lambda\mu S_{1234} \\ 
&+ (1-\lambda^2)(1-\mu^2) \, f(\text{\rm scal}) \geq 0 
\end{align*}
for $\lambda,\mu \in [0,1]$. This ensures that $\mathcal{K} \subset \mathcal{F}_0$.

\section{Ancient solutions which are weakly PIC2 and uniformly PIC}

\label{analysis.of.ancient.solutions}

In this section, we study ancient solutions to the Ricci flow of dimension $n \geq 5$ which are weakly PIC2 and uniformly PIC. An important ingredient is the differential Harnack inequality for the curvature tensor. This inequality was proved in a fundamental paper by Hamilton \cite{Hamilton3} for solutions to the Ricci flow with nonnegative curvature operator. In \cite{Brendle2}, we showed that the differential Harnack inequality holds on any solution to the Ricci flow which is weakly PIC2. \\

\begin{theorem}[cf. R.~Hamilton \cite{Hamilton3}] 
\label{harnack}
Assume that $(M,g(t))$, $t \in (0,T)$, is a solution to the Ricci flow which is complete with bounded curvature, and is weakly PIC2. Then 
\[\frac{\partial}{\partial t} \text{\rm scal} + 2 \, \langle \nabla \text{\rm scal},v \rangle + 2 \, \text{\rm Ric}(v,v) + \frac{1}{t} \, \text{\rm scal} \geq 0\] 
for every tangent vector $v$. In particular, the product $t \, \text{\rm scal}$ is monotone increasing at each point in space. 
\end{theorem} 

On an ancient solution, the Harnack inequality gives the following estimate:

\begin{corollary}
\label{harnack.for.ancient.solutions}
Assume that $(M,g(t))$ is an ancient solution to the Ricci flow which is complete with bounded curvature, and is weakly PIC2. Then 
\[\frac{\partial}{\partial t} \text{\rm scal} + 2 \, \langle \nabla \text{\rm scal},v \rangle + 2 \, \text{\rm Ric}(v,v) \geq 0\] 
for every tangent vector $v$. In particular, the scalar curvature is monotone increasing at each point in space. 
\end{corollary} 

The following inequality is obtained by integrating the differential Harnack inequality along paths in space-time: 

\begin{corollary}
\label{integrated.harnack.for.ancient.solutions}
Assume that $(M,g(t))$ is an ancient solution to the Ricci flow which is complete with bounded curvature, and is weakly PIC2. Then
\[\text{\rm scal}(x_1,t_1) \leq \exp \Big ( \frac{d_{g(t_1)}(x_1,x_2)^2}{2(t_2-t_1)} \Big ) \, \text{\rm scal}(x_2,t_2)\] 
whenever $t_1 < t_2$.
\end{corollary}

The following result from \cite{Brendle-Huisken-Sinestrari} plays a central role in our argument: 

\begin{theorem}[S.~Brendle, G.~Huisken, C.~Sinestrari \cite{Brendle-Huisken-Sinestrari}]
\label{ancient.solutions.bounded.curvature.uniformly.pic1}
Assume that $(M,g(t))$ is a complete, non-flat ancient solution to the Ricci flow with bounded curvature. Suppose that $(M,g(t))$ is uniformly PIC1, so that $R - \theta \, \text{\rm scal} \, \text{\rm id} \owedge \text{\rm id} \in \text{\rm PIC1}$ for some uniform constant $\theta>0$. Then $(M,g(t))$ has constant curvature for each $t$. 
\end{theorem}

\textbf{Proof.} 
If $M$ is compact, the assertion follows from Theorem 2 in \cite{Brendle-Huisken-Sinestrari}. The proof of Theorem 2 in \cite{Brendle-Huisken-Sinestrari} relies exclusively on the maximum principle (see Corollary 7 in that paper). Since the maximum principle works on complete manifolds with bounded curvature (see \cite{Chow-et-al}, Theorem 12.35), the arguments in \cite{Brendle-Huisken-Sinestrari} carry over to that setting. \\

We note that the conclusion of Theorem \ref{ancient.solutions.bounded.curvature.uniformly.pic1} still holds if $(M,g(t))$ has unbounded curvature (see \cite{Yokota}). 

We next state several splitting theorems, which are based on the strict maximum principle (see \cite{Hamilton2} and \cite{Brendle-book}). The first version works locally:

\begin{proposition}
\label{splitting.a}
Let $(M,g(t))$, $t \in (0,T]$, be a (possibly incomplete) solution to the Ricci flow which is weakly PIC2 and strictly PIC. Moreover, suppose that there exists a point $(p_0,t_0)$ in space-time and a unit vector $v \in T_{p_0} M$ with the property that $\text{\rm Ric}(v,v)=0$. Then, for each $t \leq t_0$, the flow $(M,g(t))$ locally splits as a product of an $(n-1)$-dimensional manifold with an interval. 
\end{proposition}

\textbf{Proof.} 
Suppose that the assertion is false. Then there exists a time $\tau \in (0,t_0)$ with the property that $(M,g(\tau))$ does not locally split as a product of an $(n-1)$-dimensional manifold with an interval. Since $(M,g(\tau))$ is strictly PIC, we conclude that $(M,g(\tau))$ is locally irreducible. 

The Ricci tensor of $(M,g(t))$ satisfies the evolution equation 
\[D_t \text{\rm Ric} = \Delta \text{\rm Ric} + 2 \, R*\text{\rm Ric},\] 
where $(R * \text{\rm Ric})_{ik} = \sum_{p,q=1}^n R_{ipkq} \text{\rm Ric}_{pq}$. Since $R$ is weakly PIC2, the term $R * \text{\rm Ric}$ is weakly positive definite. Using the strict maximum principle (see \cite{Brendle-book}, Section 9, or \cite{Hamilton2}), we conclude that the null space of $\text{\rm Ric}_{g(\tau)}$ defines a parallel subbundle of the tangent bundle of $(M,g(\tau))$. Since $(M,g(\tau))$ is locally irreducible, this subbundle must have rank $0$. Consequently, the Ricci curvature of $(M,g(\tau))$ is strictly positive.

Let $\Omega$ be a bounded open neighborhood of the point $p_0$ with smooth boundary. Let us choose a smooth function $f: \bar{\Omega} \to \mathbb{R}$ with the property that $f>0$ in $\Omega$, $f=0$ on $\partial \Omega$, and $\text{\rm Ric}_{g(\tau)} - f \, \text{\rm id}$ is weakly positive definite. Let $F: \bar{\Omega} \times [\tau,t_0] \to \mathbb{R}$ denote the solution of the linear heat equation with respect to the evolving metric $g(t)$ with initial data $F(\cdot,\tau)=f$ and Dirichlet boundary condition $F=0$ on $\partial \Omega \times [\tau,t_0]$. Using the maximum principle, we conclude that $\text{\rm Ric}_{g(t)} - F(\cdot,t) \, \text{\rm id}$ is weakly positive definite at each point in $\Omega \times [\tau,t_0]$. Since $F(p_0,t_0)>0$, the Ricci curvature at $(p_0,t_0)$ is strictly positive, contrary to our assumption. \\

\begin{proposition}
\label{splitting.b}
Let $(M,g(t))$, $t \in (0,T]$, be a complete solution to the Ricci flow (possibly with unbounded curvature). Moreover, we assume that $(M,g(t))$ is weakly PIC2 and strictly PIC. Suppose that there exists a point $(p_0,t_0)$ in space-time with the property that the curvature tensor at $(p_0,t_0)$ lies on the boundary of the PIC2 cone. Then, for each $t \leq t_0$, the universal cover of $(M,g(t))$ splits off a line.
\end{proposition}

\textbf{Proof.}
Suppose that the assertion is false. Then there exists a time $\tau \in (0,t_0)$ with the property that the universal cover of $(M,g(\tau))$ does not split off a line. Since $(M,g(\tau))$ is strictly PIC, we conclude that $(M,g(\tau))$ is locally irreducible. In view of Berger's holonomy classification theorem, there are four possible cases:

\textit{Case 1:} Suppose first that $\text{\rm Hol}^0(M,g(\tau)) = SO(n)$. If we fix an arbitrary pair of real numbers $\lambda,\mu \in [0,1]$, then Theorem 9.13 in \cite{Brendle-book} implies that the set of all orthonormal four-frames $\{e_1,e_2,e_3,e_4\}$ satisfying 
\begin{align*} 
&R_{g(\tau)}(e_1,e_3,e_1,e_3) + \lambda^2 R_{g(\tau)}(e_1,e_4,e_1,e_4) \\ 
&+ \mu^2 R_{g(\tau)}(e_2,e_3,e_2,e_3) + \lambda^2\mu^2 \, R_{g(\tau)}(e_2,e_4,e_2,e_4) \\ 
&- 2\lambda\mu \, R_{g(\tau)}(e_1,e_2,e_3,e_4) = 0 
\end{align*} 
is invariant under parallel transport. Since $\text{\rm Hol}^0(M,g(\tau)) = SO(n)$, we conclude that 
\begin{align*} 
&R_{g(\tau)}(e_1,e_3,e_1,e_3) + \lambda^2 R_{g(\tau)}(e_1,e_4,e_1,e_4) \\ 
&+ \mu^2 R_{g(\tau)}(e_2,e_3,e_2,e_3) + \lambda^2\mu^2 \, R_{g(\tau)}(e_2,e_4,e_2,e_4) \\ 
&- 2\lambda\mu \, R_{g(\tau)}(e_1,e_2,e_3,e_4) > 0 
\end{align*} 
for all orthonormal four-frames $\{e_1,e_2,e_3,e_4\}$ and all $\lambda,\mu \in [0,1]$. In other words, the curvature tensor of $(M,g(\tau))$ lies in the interior of the PIC2 cone. 

Let $\Omega$ be a bounded open neighborhood of the point $p_0$ with smooth boundary. Let us choose a smooth function $f: \bar{\Omega} \to \mathbb{R}$ with the property that $f>0$ in $\Omega$, $f=0$ on $\partial \Omega$, and $R_{g(\tau)} - \frac{1}{2} \, f \, \text{\rm id} \owedge \text{\rm id} \in \text{\rm PIC2}$. Let $F: \bar{\Omega} \times [\tau,t_0] \to \mathbb{R}$ denote the solution of the linear heat equation with respect to the evolving metric $g(t)$ with initial data $F(\cdot,\tau)=f$ and Dirichlet boundary condition $F=0$ on $\partial \Omega \times [\tau,t_0]$. Let $S := R_{g(t)} - \frac{1}{2} \, F(\cdot,t) \, \text{\rm id} \owedge \text{\rm id}$. Then $D_t S = \Delta S + Q(S) + 2F \, \text{\rm Ric}(S) \owedge \text{\rm id} + (n-1) \, F^2 \, \text{\rm id} \owedge \text{\rm id}$. Using the maximum principle, we conclude that $R_{g(t)} - \frac{1}{2} \, F(\cdot,t) \, \text{\rm id} \owedge \text{\rm id} \in \text{\rm PIC2}$ at each point in $\Omega \times [\tau,t_0]$. Since $F(p_0,t_0)>0$, the curvature tensor at $(p_0,t_0)$ lies in the interior of the PIC2 cone, contrary to our assumption. 

\textit{Case 2:} Suppose that $n=2m$ and $\text{\rm Hol}^0(M,g(\tau)) = U(m)$. In this case, the universal cover of $(M,g(\tau))$ is a K\"ahler manifold. This is impossible since $(M,g(\tau))$ is strictly PIC. 

\textit{Case 3:} Suppose next that $n=4m \geq 8$ and $\text{\rm Hol}^0(M,g(\tau)) = \text{\rm Sp}(m) \cdot \text{\rm Sp}(1)$. In this case, the universal cover of $(M,g(\tau))$ is a quaternionic K\"ahler manifold. In particular, $(M,g(\tau))$ is an Einstein manifold. Since $(M,g(\tau))$ is strictly PIC, we conclude that the universal cover of $(M,g(\tau))$ is a round sphere (cf. \cite{Brendle3}). This is impossible.

\textit{Case 4:} Suppose finally that $\text{\rm Hol}^0(M,g(\tau)) \neq SO(n)$ and $(M,g(\tau))$ is locally symmetric. In this case, $(M,g(\tau))$ is an Einstein manifold, for otherwise the eigenspaces of the Ricci tensor give a splitting of the tangent bundle into parallel subbundles, which is impossible since $(M,g(\tau))$ is locally irreducible. Since $(M,g(\tau))$ is strictly PIC, it follows that the universal cover of $(M,g(\tau))$ is a round sphere (cf. \cite{Brendle3}). This contradicts the fact that $\text{\rm Hol}^0(M,g(\tau)) \neq SO(n)$. This completes the proof of Theorem \ref{splitting.b}. \\

By combining Proposition \ref{splitting.b} and Theorem \ref{ancient.solutions.bounded.curvature.uniformly.pic1}, we can draw the following conclusion:

\begin{corollary}
\label{splitting.c}
Let $(M,g(t))$, $t \in (-\infty,T]$, be a complete, non-flat ancient solution to the Ricci flow with bounded curvature. Moreover, we assume that $(M,g(t))$ is weakly PIC2 and satisfies $R - \theta \, \text{\rm scal} \, \text{\rm id} \owedge \text{\rm id} \in \text{\rm PIC}$ for some uniform constant $\theta>0$. Suppose that there exists a point $(p_0,t_0)$ in space-time with the property that the curvature tensor at $(p_0,t_0)$ lies on the boundary of the PIC2 cone. Then, for each $t \leq t_0$, the universal cover of $(M,g(t))$ is isometric to a family of shrinking cylinders $S^{n-1} \times \mathbb{R}$.
\end{corollary}

\textbf{Proof.} 
By the strict maximum principle, the scalar curvature of $(M,g(t))$ is strictly positive. Consequently, $(M,g(t))$ is strictly PIC. By Proposition \ref{splitting.b}, the universal cover of $(M,g(t))$ is isometric to a product $(X,g_X(t)) \times \mathbb{R}$ for each $t \leq t_0$. Clearly, $(X,g_X(t))$, $t \leq t_0$, is a complete, non-flat ancient solution to the Ricci flow of dimension $n-1$, which is weakly PIC2 and uniformly PIC1. Since $(X,g_X(t))$ has bounded curvature, Theorem \ref{ancient.solutions.bounded.curvature.uniformly.pic1} implies that the universal cover of $(X,g_X(t))$ is isometric to a family of shrinking spheres. This completes the proof of Corollary \ref{splitting.c}. \\

We next recall two results due to Perelman, which play a fundamental role in the argument: 

\begin{proposition}[G.~Perelman]
\label{splitting.d}
Assume that $(M,g)$ is a complete noncompact manifold which is weakly PIC2. Fix a point $p \in M$ and let $p_j$ be a sequence of points such that $d(p,p_j) \to \infty$. Moreover, suppose that $\lambda_j$ is a sequence of positive real numbers satisfying $\lambda_j \, d(p,p_j)^2 \to \infty$. If the rescaled manifolds $(M,\lambda_j g,p_j)$ converge in the Cheeger-Gromov sense to a smooth limit, then the limit splits off a line.
\end{proposition} 

\begin{proposition}[G.~Perelman]
\label{no.small.necks}
Let $(M,g)$ be a complete noncompact Riemannian manifold which is weakly PIC2. Then $(M,g)$ does not contain a sequence of necks with radii converging to $0$.
\end{proposition}

The proofs of Proposition \ref{splitting.d} and Proposition \ref{no.small.necks} are based on Toponogov's theorem. For a detailed exposition of these results of Perelman we refer to \cite{Chen-Zhu}, Lemma 2.2 and Proposition 2.3.

We now define the class of ancient solutions that we will study. 

\begin{definition}
\label{def.of.kappa.solution}
An ancient $\kappa$-solution is a non-flat ancient solution to the Ricci flow of dimension $n$ which is complete; has bounded curvature; is weakly PIC2; and is $\kappa$-noncollapsed on all scales. 
\end{definition}

The following result is a consequence of Proposition \ref{no.small.necks} and Corollary \ref{splitting.c}.

\begin{proposition}
\label{boundedness.of.curvature}
Suppose that $(M,g(t))$, $t \in (-\infty,0]$, is a complete ancient solution to the Ricci flow which is $\kappa$-noncollapsed on all scales; is weakly PIC2; and satisfies $R - \theta \, \text{\rm scal} \, \text{\rm id} \owedge \text{\rm id} \in \text{\rm PIC}$ for some uniform constant $\theta>0$. Moreover, suppose that $(M,g(t))$ satisfies the Harnack inequality 
\[\frac{\partial}{\partial t} \text{\rm scal} + 2 \, \langle \nabla \text{\rm scal},v \rangle + 2 \, \text{\rm Ric}(v,v) \geq 0\] 
for every tangent vector $v$. Then $(M,g(t))$ has bounded curvature.
\end{proposition}

\textbf{Proof.} 
Since $(M,g(t))$ satisfies the Harnack inequality, it suffices to show that $(M,g(0))$ has bounded curvature. The proof is by contradiction. Suppose that $(M,g(0))$ has unbounded curvature. 

By the strict maximum principle, we can find a real number $\delta>0$ with the property that the scalar curvature $(M,g(t))$ is strictly positive for $t \in (-\delta,0]$. Consequently, $(M,g(t))$ is strictly PIC for $t \in (-\delta,0]$. We distinguish two cases: 

\textit{Case 1:} Suppose first that $(M,g(0))$ is strictly PIC2. In this case, $M$ is diffeomorphic to $\mathbb{R}^n$ by the soul theorem (cf. \cite{Cheeger-Gromoll2}). By a standard point-picking argument, there exists a sequence of points $x_j \in M$ such that $Q_j := \text{\rm scal}(x_j,0) \geq j^4$ and 
\[\sup_{x \in B_{g(0)}(x_j,2j \, Q_j^{-\frac{1}{2}})} \text{\rm scal}(x,0) \leq 4 Q_j.\] 
Since $(M,g(t))$ satisfies the Harnack inequality, we obtain 
\[\sup_{(x,t) \in B_{g(0)}(x_j,2j \, Q_j^{-\frac{1}{2}}) \times [-4j^2 \, Q_j^{-1},0]} \text{\rm scal}(x,t) \leq 4 Q_j.\] 
Using Shi's interior derivative estimate, we obtain bounds for all the derivatives of curvature on $B_{g(0)}(x_j,j \, Q_j^{-\frac{1}{2}}) \times [-j^2 \, Q_j^{-1},0]$. We now dilate the flow $(M,g(t))$ around the point $(x_j,0)$ by the factor $Q_j$. Using the noncollapsing assumption and the curvature derivative estimates, we conclude that, after passing to a subsequence, the rescaled flows converge in the Cheeger-Gromov sense to a smooth non-flat ancient solution $(M^\infty,g^\infty(t))$, $t \in (-\infty,0]$. The limit $(M^\infty,g^\infty(t))$ is complete; has bounded curvature; is weakly PIC2; and satisfies $R - \theta \, \text{\rm scal} \, \text{\rm id} \owedge \text{\rm id} \in \text{\rm PIC}$. By Proposition \ref{splitting.d}, the manifold $(M^\infty,g^\infty(0))$ splits off a line. By Corollary \ref{splitting.c}, universal cover of $(M^\infty,g^\infty(t))$ is isometric to a family of shrinking cylinders $S^{n-1} \times \mathbb{R}$. Therefore, $(M^\infty,g^\infty(0))$ is isometric to a quotient $(S^{n-1}/\Gamma) \times \mathbb{R}$. If $\Gamma$ is non-trivial, then a result of Hamilton implies that $M$ contains a non-trivial incompressible $(n-1)$-dimensional space form $S^{n-1}/\Gamma$ (cf. \cite{Brendle4}, Theorem A.2), but this is impossible since $M$ is diffeomorphic to $\mathbb{R}^n$. Thus, $\Gamma$ is trivial, and $(M^\infty,g^\infty(0))$ is isometric to a standard cylinder $S^{n-1} \times \mathbb{R}$. Consequently, $(M,g(0))$ contains a sequence of necks with radii converging to $0$. This contradicts Proposition \ref{no.small.necks}. 

\textit{Case 2:} Suppose finally that $(M,g(0))$ is not strictly PIC2. By Proposition \ref{splitting.b}, the universal cover of $(M,g(t))$ is isometric to a product $(X,g_X(t)) \times \mathbb{R}$ for each $t \in (-\delta,0]$. Clearly, $(X,g_X(t))$, $t \in (-\delta,0]$, is a complete solution to the Ricci flow of dimension $n-1$, which is $\kappa$-noncollapsed on all scales; weakly PIC2; and uniformly PIC1. By assumption, $(X,g_X(0))$ has unbounded curvature. Using a point-picking lemma, we can construct a sequence of points $x_j \in X$ such that $Q_j := \text{\rm scal}(x_j,0) \geq j^4$ and 
\[\sup_{x \in B_{g_X(0)}(x_j,2j \, Q_j^{-\frac{1}{2}})} \text{\rm scal}(x,0) \leq 4 Q_j.\] 
As in Case 1, the Harnack inequality and Shi's interior derivative estimate give bounds for all the derivatives of curvature on $B_{g_X(0)}(x_j,j \, Q_j^{-\frac{1}{2}}) \times [-j^2 \, Q_j^{-1},0]$. We now dilate the manifold $(X,g_X(0))$ around the point $x_j$ by the factor $Q_j$. Passing to the limit as $j \to \infty$, we obtain a smooth non-flat limit which is uniformly PIC1 and which must split off a line by Proposition \ref{splitting.d}. This is a contradiction. \\

We next recall a key result from Perelman's first paper:

\begin{theorem}[cf. G.~Perelman \cite{Perelman1}, Corollary 11.6]
\label{perelman.cor.11.6}
Given a positive real number $w>0$, we can find positive constants $B$ and $C$ (depending on $w$ and $n$) such that the following holds: Let $(M,g(t))$, $t \in [0,T]$, be a solution to the Ricci flow which is weakly PIC2. Suppose that the ball $B_{g(T)}(x_0,r_0)$ is compactly contained in $M$, and $r_0^{-n} \, \text{\rm vol}_{g(t)}(B_{g(t)}(x_0,r_0)) \geq w$ for each $t \in [0,T]$. Then $\text{\rm scal}(x,t) \leq C r_0^{-2} + Bt^{-1}$ for all $t \in (0,T]$ and all $x \in B_{g(t)}(x_0,\frac{1}{4} \, r_0)$.
\end{theorem}

Note that Perelman imposes the stronger assumption $(M,g(t))$ has nonnegative curvature operator. However, his proof works under the weaker assumption that $(M,g(t))$ is weakly PIC2. One main ingredient in Perelman's work is the trace Harnack inequality (see Theorem \ref{harnack}). The proof also relies on the fact that a solution to the Ricci flow which has evolved for some positive time cannot contain an open set which is isometric to a piece of a non-flat cone. This argument relies on the strict maximum principle, and works if the solution is weakly PIC2 (see Proposition \ref{splitting.a}).

One of the main tools in Perelman's theory is the longrange curvature estimate for ancient $\kappa$-solutions in dimension $3$. In the next step, we verify that this estimate holds in our situation. \\

\begin{theorem}[cf. G.~Perelman \cite{Perelman1}, Section 11.7]
\label{longrange.curvature.estimate}
Given $\kappa>0$, we can find a positive function $\omega: [0,\infty) \to (0,\infty)$ (depending on $n$ and $\kappa$) with the following property: Let $(M,g(t))$ be a an ancient $\kappa$-solution. Then 
\[\text{\rm scal}(x,t) \leq \text{\rm scal}(y,t) \, \omega(\text{\rm scal}(y,t) \, d_{g(t)}(x,y)^2)\] 
for all points $x,y \in M$ and all $t$. 
\end{theorem}

\textbf{Proof.} 
The proof is essentially the same as in Section 11.7 in Perelman's paper \cite{Perelman1} (see also \cite{Chen-Zhu}). We sketch the argument for the convenience of the reader. Let us fix a point $y \in M$. By rescaling, we can arrange that $\text{\rm scal}(y,0) = 1$. For abbreviation, let $A = \{x \in M: \text{\rm scal}(x,0) \, d_{g(0)}(y,x)^2 \geq 1\}$. We distinguish two cases: 

\textit{Case 1:} Suppose that $A = \emptyset$. In this case, we can find a point $z \in M$ such that $\sup_{x \in M} \text{\rm scal}(x,0) = \text{\rm scal}(z,0)$. Using the Harnack inequality, we obtain  
\[\sup_{x \in M} \text{\rm scal}(x,t) \leq \text{\rm scal}(z,0)\] 
for all $t \in (-\infty,0]$. By Shi's derivative estimates, $\frac{\partial}{\partial t} \text{\rm scal}(z,t) \leq C(n) \, \text{\rm scal}(z,0)^2$ for all $t \in [-\text{\rm scal}(z,0)^{-1},0]$. Moreover, $d_{g(0)}(y,z) \leq \text{\rm scal}(z,0)^{-\frac{1}{2}}$ since $A = \emptyset$. Hence, we can find a small positive constant $\beta$, depending only on $n$, such that $\text{\rm scal}(z,t) \geq \frac{1}{2} \, \text{\rm scal}(z,0)$ and $d_{g(t)}(y,z) \leq 2 \, \text{\rm scal}(z,0)^{-\frac{1}{2}}$ for all $t \in [-\beta \, \text{\rm scal}(z,0)^{-1},0]$. If we apply the Harnack inequality (cf. Corollary \ref{integrated.harnack.for.ancient.solutions}) with $t = -\beta \, \text{\rm scal}(z,0)^{-1}$, then we obtain 
\begin{align*} 
\frac{1}{2} \, \text{\rm scal}(z,0) 
&\leq \text{\rm scal}(z,t) \\ 
&\leq \exp \Big ( \frac{d_{g(t)}(y,z)^2}{(-2t)} \Big ) \, \text{\rm scal}(y,0) \\ 
&\leq \exp \Big ( \frac{2}{(-t) \, \text{\rm scal}(z,0)} \Big ) \, \text{\rm scal}(y,0) \\ 
&= \exp \Big ( \frac{2}{\beta} \Big ). 
\end{align*} 
Putting these facts together, we conclude that $\sup_{x \in M} \text{\rm scal}(x,0) \leq 2 \, \exp \big ( \frac{2}{\beta} \big )$. 

\textit{Case 2:} Suppose now that $A \neq \emptyset$. In this case, we choose a point $z \in A$ which has minimal distance from $y$ with respect to the metric $g(0)$. Note that $\text{\rm scal}(z,0) = d_{g(0)}(y,z)^{-2}$ since $z$ lies on the boundary of $A$. Let $p$ be the mid-point of the minimizing geodesic in $(M,g(0))$ joining $y$ and $z$. Clearly, $B_{g(0)}(p,\frac{1}{4} \, d_{g(0)}(y,z)) \cap A = \emptyset$. This implies 
\[\sup_{x \in B_{g(0)}(p,\frac{1}{4} \, d_{g(0)}(y,z))} \text{\rm scal}(x,0) \leq 16 \, d_{g(0)}(y,z)^{-2}.\]
By the Harnack inequality, 
\[\sup_{x \in B_{g(t)}(p,\frac{1}{4} \, d_{g(0)}(y,z))} \text{\rm scal}(x,t) \leq 16 \, d_{g(0)}(y,z)^{-2}\]  
for all $t \in (-\infty,0]$. The noncollapsing assumption gives 
\[\text{\rm vol}_{g(t)} \big ( B_{g(t)}(p,\frac{1}{4} \, d_{g(0)}(y,z)) \big ) \geq \kappa \, (\frac{1}{4} \, d_{g(0)}(y,z))^n\] 
for all $t \in (-\infty,0]$. This implies 
\[(4r)^{-n} \, \text{\rm vol}_{g(t)}(B_{g(t)}(p,4r)) \geq \kappa \, (\frac{1}{16} \, r^{-1} \, d_{g(0)}(y,z))^n\] 
for all $t \in (-\infty,0]$ and all $r \geq d_{g(0)}(y,z)$. Applying Theorem \ref{perelman.cor.11.6} with $x_0 := p$, $r_0 := 4r$, and $w := \kappa \, (\frac{1}{16} \, r^{-1} \, d_{g(0)}(y,z))^n$, we obtain 
\[\sup_{x \in B_{g(0)}(p,r)} \text{\rm scal}(x,0) \leq d_{g(0)}(y,z)^{-2} \, \omega(d_{g(0)}(y,z)^{-1} \, r)\] 
for all $r \geq d_{g(0)}(y,z)$, where $\omega: [1,\infty) \to [0,\infty)$ is a positive and increasing function that may depend on $n$ and $\kappa$. In particular, if we put $r = d_{g(0)}(y,z)$ and apply the Harnack inequality, we obtain  
\[\sup_{x \in B_{g(0)}(p,d_{g(0)}(y,z))} \text{\rm scal}(x,t) \leq d_{g(0)}(y,z)^{-2} \, \omega(1),\] 
for all $t \in (-\infty,0]$. By Shi's derivative estimates, $\frac{\partial}{\partial t} \text{\rm scal}(z,t) \leq C(n,\kappa) \, d_{g(0)}(y,z)^{-4}$ for all $t \in [-d_{g(0)}(y,z)^2,0]$. Moreover, $\text{\rm scal}(z,0) = d_{g(0)}(y,z)^{-2}$ by our choice of $z$. Therefore, we can find a small positive constant $\beta$, depending only on $n$ and $\kappa$, such that $\text{\rm scal}(z,t) \geq \frac{1}{2} \, d_{g(0)}(y,z)^{-2}$ and $d_{g(t)}(y,z) \leq 2 \, d_{g(0)}(y,z)$ for all $t \in [-\beta \, d_{g(0)}(y,z)^2,0]$. If we apply the Harnack inequality (cf. Corollary \ref{integrated.harnack.for.ancient.solutions}) with $t = -\beta \, d_{g(0)}(y,z)^2$, then we obtain 
\begin{align*} 
\frac{1}{2} \, d_{g(0)}(y,z)^{-2} 
&\leq \text{\rm scal}(z,t) \\ 
&\leq \exp \Big ( \frac{d_{g(t)}(y,z)^2}{(-2t)} \Big ) \, \text{\rm scal}(y,0) \\ 
&\leq \exp \Big ( \frac{2 \, d_{g(0)}(y,z)^2}{(-t)} \Big ) \, \text{\rm scal}(y,0) \\ 
&= \exp \Big ( \frac{2}{\beta} \Big ). 
\end{align*} 
This finally implies 
\begin{align*} 
\sup_{x \in B_{g(0)}(y,r)} \text{\rm scal}(x,0) 
&\leq \sup_{x \in B_{g(0)}(p,r+d_{g(0)}(y,z))} \text{\rm scal}(x,0) \\ 
&\leq d_{g(0)}(y,z)^{-2} \, \omega(d_{g(0)}(y,z)^{-1} \, r+1) \\ 
&\leq 2 \, e^{\frac{2}{\beta}} \, \omega(\sqrt{2} \, e^{\frac{1}{\beta}} \, r+1) 
\end{align*} 
for all $r \geq 0$. This completes the proof of Theorem \ref{longrange.curvature.estimate}. \\

Combining the longrange curvature estimate in Theorem \ref{longrange.curvature.estimate} with Shi's interior derivative estimates, we can bound the $m$-th covariant derivative of the Riemann curvature tensor by a constant times $\text{\rm scal}^{\frac{m+2}{2}}$ at each point in space-time. In particular, we can draw the following conclusion:

\begin{corollary}
\label{pointwise.derivative.estimate}
Given $\kappa>0$, we can find a large positive constant $\eta = \eta(n,\kappa)$ with the following property: Let $(M,g(t))$ be a an ancient $\kappa$-solution. Then $|D \text{\rm scal}| \leq \eta \, \text{\rm scal}^{\frac{3}{2}}$ and $|\frac{\partial}{\partial t} \text{\rm scal}| \leq \eta \, \text{\rm scal}^2$ at each point in space-time. 
\end{corollary}

We next establish an analogue of Perelman's compactness theorem for ancient $\kappa$-solutions:

\begin{corollary}[cf. G.~Perelman \cite{Perelman1}, Section 11.7]
\label{compactness.thm.for.ancient.solutions}
Fix $\kappa>0$ and $\theta > 0$. Assume that $(M^{(j)},g^{(j)}(t))$ is a sequence of ancient $\kappa$-solutions satisfying $R - \theta \, \text{\rm scal} \, \text{\rm id} \owedge \text{\rm id} \in \text{\rm PIC}$. Suppose that $x_j$ is a point on $M^{(j)}$ satisfying $\text{\rm scal}(x_j,0) = 1$. Then, after passing to a subsequence if necessary, the sequence $(M^{(j)},g^{(j)}(t),x_j)$ converges in the Cheeger-Gromov sense to an ancient $\kappa$-solution satisfying $R - \theta \, \text{\rm scal} \, \text{\rm id} \owedge \text{\rm id} \in \text{\rm PIC}$.
\end{corollary}

\textbf{Proof.} 
It follows from the noncollapsing assumption and the longrange curvature estimate in Theorem \ref{longrange.curvature.estimate} that, after passing to a subsequence if necessary, the sequence $(M^{(j)},g^{(j)}(t),x_j)$ converges in the Cheeger-Gromov sense to a smooth non-flat ancient solution $(M^\infty,g^\infty(t))$. Clearly, $(M^\infty,g^\infty(t))$ is $\kappa$-noncollapsed on all scales; is weakly PIC2; and satisfies $R - \theta \, \text{\rm scal} \, \text{\rm id} \owedge \text{\rm id} \in \text{\rm PIC}$. It remains to show that $(M^\infty,g^\infty(t))$ has bounded curvature. Since $(M^{(j)},g^{(j)}(t))$ has bounded curvature, the flow $(M^{(j)},g^{(j)}(t))$ satisfies the Harnack inequality 
\[\frac{\partial}{\partial t} \text{\rm scal} + 2 \, \langle \nabla \text{\rm scal},v \rangle + 2 \, \text{\rm Ric}(v,v) \geq 0\] 
(cf. Corollary \ref{harnack.for.ancient.solutions}). Passing to the limit as $j \to \infty$, we conclude that the Harnack inequality holds on the limit $(M^\infty,g^\infty(t))$. Consequently, $(M^\infty,g^\infty(t))$ has bounded curvature by Proposition \ref{boundedness.of.curvature}. \\

In the remainder of this section, we establish a structure theorem and a universal noncollapsing theorem for ancient $\kappa$-solutions. This was first established by Perelman \cite{Perelman2} in the three-dimensional case, and adapted to dimension $4$ in \cite{Chen-Zhu}. 

\begin{definition}
\label{definition.of.cap}
Let us fix a small number $\varepsilon_0 = \varepsilon_0(n)$, and let $0 < \varepsilon < \frac{\varepsilon_0}{4}$. We say that a compact domain $\Omega \subset (M,g)$ is an $\varepsilon$-cap if the following holds: 
\begin{itemize}
\item The domain $\Omega$ is diffeomorphic to $B^n$ and the boundary $\partial \Omega$ is a cross-sectional sphere of an $\varepsilon$-neck. \item If $\tilde{\Omega} \subset \Omega$ is a compact domain diffeomorphic to $B^n$ and the boundary $\partial \tilde{\Omega}$ is a cross-sectional sphere of an $(\varepsilon_0-\varepsilon)$-neck, then there exists a diffeomorphism $F: \tilde{\Omega} \to B^n$ and an $(\varepsilon_0+\varepsilon)$-isometry $f: \partial \tilde{\Omega} \to S^{n-1}$ with the property that $F|_{\partial \tilde{\Omega}}: \partial \tilde{\Omega} \to S^{n-1}$ is isotopic to $f$.
\end{itemize}
\end{definition}

The second condition in Definition \ref{definition.of.cap} will play a crucial role when we analyze how the diffeomorphism type of the manifold changes under surgery (see Proposition \ref{change.of.topology} below).

\begin{proposition} 
\label{topology.of.cap}
Let $(M,g)$ be a complete, noncompact manifold of dimension $n \geq 5$ which is strictly PIC2. Suppose that $\Sigma$ is a cross-sectional sphere of an $\varepsilon_0$-neck $N$ in $(M,g)$. Then $\Sigma$ bounds a compact domain $\Omega$. Moreover, there exists a diffeomorphism $F: \Omega \to B^n$ and an $\varepsilon_0$-isometry $f: \partial \Omega \to S^{n-1}$ with the property that $F|_{\partial \Omega}: \partial \Omega \to S^{n-1}$ is isotopic to $f$.
\end{proposition}

\textbf{Proof.} 
By the soul theorem (cf. \cite{Cheeger-Gromoll2}), $M$ is diffeomorphic to $\mathbb{R}^n$. By the solution of the Schoenflies conjecture in dimension $n \geq 5$, $\Sigma$ bounds a compact domain $\Omega$ which is diffeomorphic to $B^n$. Let us choose a diffeomorphism $F: \Omega \to B^n$, and let $f: \partial \Omega \to S^{n-1}$ be an $\varepsilon_0$-isometry. Without loss of generality, we may assume that the map $F|_{\partial \Omega} \circ f^{-1}: S^{n-1} \to S^{n-1}$ is orientation-preserving. We claim that the map $F|_{\partial \Omega} \circ f^{-1}: S^{n-1} \to S^{n-1}$ is isotopic to the identity. To prove this, let $z$ denote the height function on the neck $N$, so that $\{z=0\} = \Sigma$ and $\{z \in [-10,0)\} \subset \Omega$. Moreover, let $\varphi(z) = e^{-\frac{1}{z}}$ for $z \in (0,\frac{1}{10}]$. Let us define a modified metric $\tilde{g}$ as follows. The metric $\tilde{g}$ agrees with $g$ inside $\Omega$. In the region $\{z \in (0,\frac{1}{20}]\}$, we change the metric conformally by $\tilde{g} = e^{-2\varphi} \, g$. In the region $\{z \in (\frac{1}{20},\frac{1}{10}]\}$, we define $\tilde{g} = e^{-2\varphi} \, (\chi(z) \, g + (1-\chi(z)) \, \bar{g})$, where $\bar{g}$ denotes the standard metric on the cylinder and $\chi: (\frac{1}{20},\frac{1}{10}] \to [0,1]$ is a smooth cutoff function satisfying $\chi(z)=1$ for $z \in [\frac{1}{20},\frac{1}{18}]$ and $\chi(z)=0$ for $z \in [\frac{1}{12},\frac{1}{10}]$. Since the original metric $g$ is strictly PIC2, the modified metric $\tilde{g}$ will be strictly PIC2 in the region $\{z \in (0,\frac{1}{10}]\}$. Moreover, the surgically modified metric $\tilde{g}$ is rotationally symmetric for $z \in [\frac{1}{12},\frac{1}{10}]$. Consequently, we can extend $\tilde{g}$ by gluing in a rotationally symmetric cap. This can be done in such a way that the glued metric is strictly PIC2. To summarize, we obtain a closed manifold which is obtained by gluing two balls, and which admits a metric which is strictly PIC2. If the map $F|_{\partial \Omega} \circ f^{-1}$ is not isotopic to the identity, then this glued manifold is an exotic sphere, which contradicts the fact that exotic spheres do not admit metrics which are strictly PIC2 (cf. \cite{Brendle-Schoen}). \\

We now state the structure theorem in the noncompact case.

\begin{theorem}[cf. G.~Perelman \cite{Perelman1}, Corollary 11.8; Chen-Zhu \cite{Chen-Zhu}, Proposition 3.4]
\label{structure.theorem.noncompact.case}
Given $\varepsilon>0$ and $\theta>0$, we can find large positive constants $C_1 = C_1(n,\theta,\varepsilon)$ and $C_2 = C_2(n,\theta,\varepsilon)$ with the following property: Suppose that $(M,g(t))$ is a noncompact ancient $\kappa$-solution satisfying $R - \theta \, \text{\rm scal} \, \text{\rm id} \owedge \text{\rm id} \in \text{\rm PIC}$ which is not locally isometric to a family of shrinking cylinders. Then, for each point $(x_0,t_0)$ in space-time, we can find a neighborhood $B$ of $x_0$ satisfying $B_{g(t_0)}(x_0,C_1^{-1} \, \text{\rm scal}(x_0,t_0)^{-\frac{1}{2}}) \subset B \subset B_{g(t_0)}(x_0,C_1 \, \text{\rm scal}(x_0,t_0)^{-\frac{1}{2}})$ and $C_2^{-1} \, \text{\rm scal}(x_0,t_0) \leq \text{\rm scal}(x,t_0) \leq C_2 \, \text{\rm scal}(x_0,t_0)$ for all $x \in B$. Moreover, $B$ satisfies one of the following conditions: 
\begin{itemize}
\item $B$ is an $\varepsilon$-neck with center at $x_0$.
\item $B$ is an $\varepsilon$-cap in the sense of Definition \ref{definition.of.cap}. 
\end{itemize}
In particular, $(M,g(t_0))$ is $\kappa_0$-noncollapsed for some universal constant $\kappa_0 = \kappa_0(n,\theta)$.
\end{theorem}

A key point is that the constants $C_1$ and $C_2$ in Theorem \ref{structure.theorem.noncompact.case} do not depend on $\kappa$. \\

\textbf{Proof.} 
It suffices to prove the assertion for $t_0=0$. Suppose that $x_0$ does not lie at the center of an $\varepsilon$-neck in $(M,g(0))$. Since $(M,g(t))$ is not locally isometric to a family of shrinking cylinders, Corollary \ref{splitting.c} implies that $(M,g(t))$ is strictly PIC2. By the soul theorem (cf. \cite{Cheeger-Gromoll2}), $M$ is diffeomorphic to $\mathbb{R}^n$. We denote by $M_\varepsilon$ the set of all points which do not lie at the center of an $\frac{\varepsilon}{4}$-neck in $(M,g(0))$. Note that $x_0$ lies in the interior of $M_\varepsilon$.

\textit{Step 1:} We claim that the closure of $M_\varepsilon$ is compact. If this is false, we can find a sequence of points $x_j \in M_\varepsilon$ such that $d_{g(0)}(x_0,x_j) \to \infty$. The longrange curvature estimate in Theorem \ref{longrange.curvature.estimate} implies that $\lim_{j \to \infty} \text{\rm scal}(x_j,0) \, d_{g(0)}(x_0,x_j)^2 = \infty$. We now dilate the flow $(M,g(t))$ around the point $(x_j,0)$ by the factor $\text{\rm scal}(x_j,0)$. By Corollary \ref{compactness.thm.for.ancient.solutions}, the rescaled flows converge in the Cheeger-Gromov sense to an ancient $\kappa$-solution $(M^\infty,g^\infty(t))$ satisfying $R - \theta \, \text{\rm scal} \, \text{\rm id} \owedge \text{\rm id} \in \text{\rm PIC}$. By Proposition \ref{splitting.d}, the manifold $(M^\infty,g^\infty(0))$ splits off a line. Since the limit $(M^\infty,g^\infty(t))$ has bounded curvature, Corollary \ref{splitting.c} implies that the universal cover of $(M^\infty,g^\infty(t))$ is isometric to a family of shrinking cylinders $S^{n-1} \times \mathbb{R}$. Therefore, $(M^\infty,g^\infty(0))$ is isometric to a quotient $(S^{n-1}/\Gamma) \times \mathbb{R}$. If $\Gamma$ is non-trival, then Theorem A.2 in \cite{Brendle4} implies that $M$ contains a non-trivial incompressible $(n-1)$-dimensional space form $S^{n-1}/\Gamma$, but this is impossible since $M$ is diffeomorphic to $\mathbb{R}^n$. Thus, $\Gamma$ is trivial, and $(M^\infty,g^\infty(0))$ is isometric to a standard cylinder $S^{n-1} \times \mathbb{R}$. Consequently, $x_j$ lies at the center of an $\frac{\varepsilon}{8}$-neck if $j$ is sufficiently large. This contradicts the assumption that $x_j \in M_\varepsilon$. 

\textit{Step 2:} We now consider an arbitrary point $y \in \partial M_\varepsilon$. By definition of $M_\varepsilon$, $y$ lies at the center of an $\frac{\varepsilon}{2}$-neck in $(M,g(0))$. The Harnack inequality gives $\text{\rm scal}(x,t) \leq \text{\rm scal}(x,0) \leq 2 \, \text{\rm scal}(y,0)$ for all $x \in B_{g(0)}(y,\text{\rm scal}(y,0)^{-\frac{1}{2}})$ and all $t \leq 0$. Using Shi's derivative estimates, we obtain $\frac{\partial}{\partial t} \text{\rm scal}(y,t) \leq C(n) \, \text{\rm scal}(y,0)^2$ for all $t \in [-\text{\rm scal}(y,0)^{-1},0]$. Hence, we can find a small constant $\beta=\beta(n)>0$ such that $\text{\rm scal}(y,t) \geq \frac{1}{2} \, \text{\rm scal}(y,0)$ and 
\[\text{\rm vol}_{g(t)} \big ( B_{g(t)}(y,\text{\rm scal}(y,0)^{-\frac{1}{2}}) \big ) \geq \beta \, \text{\rm scal}(y,0)^{-\frac{n}{2}}\] 
for all $t \in [-\beta \, \text{\rm scal}(y,0)^{-1},0]$. Applying Theorem \ref{perelman.cor.11.6}, we obtain 
\[\text{\rm scal}(x,0) \leq \text{\rm scal}(y,0) \, \omega(\text{\rm scal}(y,0) \, d_{g(0)}(x,y)^2)\] 
for all $x \in M$, where $\omega: [0,\infty) \to [0,\infty)$ denotes a positive function that does not depend on $\kappa$. Using the Harnack inequality, we conclude that 
\[\text{\rm scal}(x,t) \leq \text{\rm scal}(y,0) \, \omega(\text{\rm scal}(y,0) \, d_{g(0)}(x,y)^2)\] 
for all $x \in M$ and all $t \leq 0$.

\textit{Step 3:} We again consider an arbitrary point $y \in \partial M_\varepsilon$. Recall that $y$ lies at the center of an $\frac{\varepsilon}{2}$-neck in $(M,g(0))$. By work of Hamilton \cite{Hamilton5}, a neck admits a unique foliation by CMC spheres. Let $\Sigma_y$ denote the leaf of the CMC foliation of $(M,g(0))$ which passes through $y$. By Proposition \ref{topology.of.cap}, $\Sigma_y$ bounds a compact domain $\Omega_y$, and $\Omega_y$ is an $\varepsilon$-cap in the sense of Definition \ref{definition.of.cap}. 

We next show that $\text{\rm scal}(y,0) \, \text{\rm diam}_{g(0)}(\Omega_y)^2 \leq C$, where $C$ depends on $n$, $\theta$, and $\varepsilon$, but not on $\kappa$. The proof is by contradiction. Suppose that $(M^{(j)},g^{(j)}(t))$ is a sequence of noncompact ancient $\kappa_j$-solutions which satisfy $R - \theta \, \text{\rm scal} \, \text{\rm id} \owedge \text{\rm id} \in \text{\rm PIC}$ and which are not locally isometric to a family of shrinking cylinders. Suppose further that $y_j \in M^{(j)}$ is a sequence of points such that $y_j \in \partial M_\varepsilon^{(j)}$ and $\text{\rm scal}(y_j,0) \, \text{\rm diam}_{g^{(j)}(0)}(\Omega_{y_j})^2 \to \infty$, where $\Omega_{y_j}$ denotes the region in $(M^{(j)},g^{(j)}(0))$ which is bounded by the CMC sphere passing through $y_j$. We dilate the flow $(M^{(j)},g^{(j)}(t))$ around the point $(y_j,0)$ by the factor $\text{\rm scal}(y_j,0)$. Using the longrange curvature estimate established in Step 2, we conclude that, after passing to a subsequence if necessary, the rescaled manifolds converge to a complete, smooth, non-flat ancient solution $(M^\infty,g^\infty(t))$ which is weakly PIC2 and satisfies $R - \theta \, \text{\rm scal} \, \text{\rm id} \owedge \text{\rm id} \in \text{\rm PIC}$. Since $y_j \in \partial M_\varepsilon^{(j)}$ for each $j$, the limit $y_\infty := \lim_{j \to \infty} y_j$ lies at the center of an $\frac{\varepsilon}{2}$-neck in $(M^\infty,g^\infty(0))$. Since $\text{\rm scal}(y_j,0) \, \text{\rm diam}_{g^{(j)}(0)}(\Omega_{y_j})^2 \to \infty$, the manifold $(M^\infty,g^\infty(0))$ contains a minimizing geodesic line. By the Cheeger-Gromoll splitting theorem (cf. \cite{Cheeger-Gromoll1}), the manifold $(M^\infty,g^\infty(0))$ is isometric to a product $X \times \mathbb{R}$. Since $y_\infty$ lies at the center of an $\frac{\varepsilon}{2}$-neck, the cross-section $X$ is compact, and is nearly isometric to $S^{n-1}$. In particular, $(M^\infty,g^\infty(0))$ has bounded curvature. By Corollary \ref{splitting.c}, $(M^\infty,g^\infty(t))$ is isometric to family of shrinking cylinders $S^{n-1} \times \mathbb{R}$. Therefore, if $j$ is sufficiently large, then $y_j$ lies at the center of an $\frac{\varepsilon}{8}$-neck. This contradicts the fact that $y_j$ lies on the boundary of $M_\varepsilon^{(j)}$.

\textit{Step 4:} By Step 3, we have $\text{\rm scal}(y,0) \, d_{g(0)}(x,y)^2 \leq C$ for all $y \in \partial M_\varepsilon$ and all $x \in \Omega_y$, where $C$ is a positive constant that depends only on $n$, $\theta$, and $\varepsilon$, but not on $\kappa$. Combining this estimate with the longrange curvature estimate in Step 2, we conclude that $\text{\rm scal}(x,t) \leq C \, \text{\rm scal}(y,0)$ and $\text{\rm scal}(y,0) \, d_{g(t)}(x,y)^2 \leq C$ for all $y \in \partial M_\varepsilon$, all $x \in \Omega_y$, and all $t \in [-\beta \, \text{\rm scal}(y,0)^{-1},0]$. Here, $C$ is a positive constant that depends only on $n$, $\theta$, and $\varepsilon$, but not on $\kappa$. Using the Harnack inequality (cf. Corollary \ref{integrated.harnack.for.ancient.solutions}), we obtain 
\[\text{\rm scal}(x,0) \geq \frac{1}{C} \, \text{\rm scal}(y,-\beta \, \text{\rm scal}(y,0)^{-1}) \geq \frac{1}{2C} \, \text{\rm scal}(y,0)\] 
for all $y \in \partial M_\varepsilon$ and all $x \in \Omega_y$, where $C$ depends only on $n$, $\theta$, and $\varepsilon$, but not on $\kappa$. 

\textit{Step 5:} Finally, we claim that we can find a point $y \in \partial M_\varepsilon$ such that $x_0 \in \Omega_y$. To prove this, let $\gamma: [0,\infty) \to (M,g(0))$ be a minimizing unit-speed geodesic with $\gamma(0) = x_0$. We define $\bar{s} := \sup \{s \in [0,\infty): \gamma(s) \in M_\varepsilon\}$. Note that $\bar{s} > 0$ since $x_0$ lies in the interior of $M_\varepsilon$, and $\bar{s} < \infty$ since the closure of $M_\varepsilon$ is compact. Let $\bar{y} := \gamma(\bar{s})$. Clearly, $\bar{y} \in \partial M_\varepsilon$. In particular, $\bar{y}$ lies at the center of an $\frac{\varepsilon}{2}$-neck. Since $x_0$ does not lie at the center of an $\varepsilon$-neck, we must have $d_{g(0)}(x_0,\bar{y}) > 100n \, \text{\rm scal}(\bar{y},0)^{-\frac{1}{2}}$. We claim that $x_0 \in \Omega_{\bar{y}}$. Suppose this is false. Then the curve $\gamma$ must enter and exit the cap $\Omega_{\bar{y}}$. More precisely, we have $\gamma(0) \notin \Omega_{\bar{y}} \cup B_{g(0)}(\bar{y},100n \, \text{\rm scal}(\bar{y},0)^{-\frac{1}{2}})$, $\gamma(\bar{s}) \in \Omega_{\bar{y}}$, and $\gamma(s) \notin \Omega_{\bar{y}} \cup B_{g(0)}(\bar{y},100n \, \text{\rm scal}(\bar{y},0)^{-\frac{1}{2}})$ when $s$ is sufficiently large. This is impossible since $\gamma$ minimizes length. Consequently, $x_0 \in \Omega_{\bar{y}}$. The inequality $d_{g(0)}(x_0,\bar{y}) > 100n \, \text{\rm scal}(\bar{y},0)^{-\frac{1}{2}}$ implies $B_{g(0)}(x_0,C^{-1} \, \text{\rm scal}(\bar{y},0)^{-\frac{1}{2}}) \subset \Omega_{\bar{y}}$. Moreover, the diameter estimate in Step 3 gives $\Omega_{\bar{y}} \subset B_{g(0)}(x_0,C \, \text{\rm scal}(\bar{y},0)^{-\frac{1}{2}})$. Here, $C$ is a positive constant that depends only on $n$, $\theta$, and $\varepsilon$, but not on $\kappa$. Finally, by Step 4, we can estimate $\text{\rm scal}(x_0,0)$ from above and below in terms of $\text{\rm scal}(\bar{y},0)$. Putting these facts together, we conclude that $B_{g(0)}(x_0,C^{-1} \, \text{\rm scal}(x_0,0)^{-\frac{1}{2}}) \subset \Omega_{\bar{y}}$ and $\Omega_{\bar{y}} \subset B_{g(0)}(x_0,C \, \text{\rm scal}(\bar{y},0)^{-\frac{1}{2}})$, where $C$ is a positive constant that depends only on $n$, $\theta$, and $\varepsilon$, but not on $\kappa$. To summarize, the set $B := \Omega_{\bar{y}}$ has all the required properties. \\

\begin{theorem}[cf. G.~Perelman \cite{Perelman1}; Chen-Zhu \cite{Chen-Zhu}]
\label{universal.noncollapsing} 
Fix $\theta>0$. We can find a constant $\kappa_0 = \kappa_0(n,\theta)$ such that the following holds: Suppose that $(M,g(t))$ is an ancient $\kappa$-solution for some $\kappa>0$, which in addition satisfies $R - \theta \, \text{\rm scal} \, \text{\rm id} \owedge \text{\rm id} \in \text{\rm PIC}$. Then either $(M,g(t))$ is $\kappa_0$-noncollapsed for all $t$; or $(M,g(t))$ is a metric quotient of the round sphere $S^n$; or $(M,g(t))$ is a noncompact metric quotient of the standard cylinder $S^{n-1} \times \mathbb{R}$.
\end{theorem}

\textbf{Proof.} 
If $M$ is noncompact, the assertion follows from Theorem \ref{structure.theorem.noncompact.case}. Hence, we may assume that $M$ is compact. Moreover, we may assume that $(M,g(t))$ is not a metric quotient of the round sphere $S^n$. The noncollapsing assumption implies that $(M,g(t))$ cannot be a compact quotient of a standard cylinder. By Corollary \ref{splitting.c}, $(M,g(t))$ is strictly PIC2. 

Consider a sequence of times $t_j \to -\infty$ and a sequence of points $q_j \in M$ such that $\ell(q_j,t_j) \leq n$, where $\ell$ denotes the reduced distance. By Perelman's work, $\ell(x,t_j) + (-t_j) \, R(x,t_j) \leq C(n)$ for all $x \in B_{g(t_j)}(q_j,(-t_j)^{\frac{1}{2}})$. If we dilate the flow $(M,g(t))$ around the point $(q_j,t_j)$ by the factor $(-t_j)^{-1}$, then the rescaled flows converge to a complete ancient solution $(\bar{M},\bar{g}(t))$ (cf. \cite{Perelman1}, Section 11.2). Perelman proved that the limit $(\bar{M},\bar{g}(t))$ is a non-flat shrinking gradient Ricci soliton. Clearly, $(\bar{M},\bar{g}(t))$ is $\kappa$-noncollapsed, weakly PIC2, and satisfies $R - \theta \, \text{\rm scal} \, \text{\rm id} \owedge \text{\rm id} \in \text{\rm PIC}$. Moreover, since the Harnack inequality holds on $(M,g(t))$, it holds on the limit $(\bar{M},\bar{g}(t))$. By Proposition \ref{boundedness.of.curvature}, $(\bar{M},\bar{g}(t))$ has bounded curvature. We distinguish two cases:

\textit{Case 1:} We first consider the case that $\bar{M}$ is compact. A shrinking soliton $(\bar{M},\bar{g}(t))$ cannot be isometric to a compact quotient of a standard cylinder. By Corollary \ref{splitting.c}, $(\bar{M},\bar{g}(t))$ is strictly PIC2. Since $(\bar{M},\bar{g}(t))$ is a shrinking soliton, results in \cite{Brendle-Schoen} imply that $(\bar{M},\bar{g}(t))$ must be isometric to a metric quotient of the round sphere $S^n$. Consequently, the original ancient solution $(M,g(t))$ is a metric quotient of $S^n$, contrary to our assumption. 

\textit{Case 2:} Suppose next that $\bar{M}$ is noncompact. In view of Theorem \ref{structure.theorem.noncompact.case}, there are two possibilities: either $(\bar{M},\bar{g}(t))$ is $\kappa_0$-noncollapsed for some universal constant $\kappa_0=\kappa_0(n,\theta)$, or else $(\bar{M},\bar{g}(t))$ is isometric to a noncompact quotient of the standard cylinder. The second case can be divided into two subcases: 
\begin{itemize} 
\item If the dimension $n$ is odd and $\bar{M}$ is a noncompact quotient $(S^{n-1} \times \mathbb{R})/\Gamma$, then there are only finitely many possibilities for the group $\Gamma$, and the resulting quotients are all noncollapsed with a universal constant. 
\item If the dimension $n$ is even and $\bar{M}$ is a noncompact quotient $(S^{n-1} \times \mathbb{R})/\Gamma$, then the center slice $(S^{n-1} \times \{0\})/\Gamma$ is incompressible in $M$ by Theorem A.1 in \cite{Brendle4}. Moreover, since $(M,g(t))$ is strictly PIC2, the fundamental group of $M$ has order at most $2$ by Synge's theorem. Hence, there are only finitely many possibilities for the group $\Gamma$, and the resulting quotients $(S^{n-1} \times \mathbb{R})/\Gamma$ are noncollapsed with a universal constant. 
\end{itemize} 
To summarize, we have shown that the asymptotic shrinking soliton $(\bar{M},\bar{g}(t))$ is noncollapsed with a universal constant that depends only on $n$ and $\theta$. This implies $\liminf_{j \to \infty} (-t_j)^{-\frac{n}{2}} \, \text{\rm vol}_{g(t_j)} \big ( B_{g(t_j)}(q_j,(-t_j)^{\frac{1}{2}}) \big ) \geq \frac{1}{C(n,\theta)}$. Hence, if $V(t) = \int_M (-t)^{-\frac{n}{2}} \, e^{-\ell(\cdot,t)} \, d\text{\rm vol}_{g(t)}$ denotes the reduced volume, then $\liminf_{j \to \infty} V(t_j) \geq \frac{1}{C(n,\theta)}$. Perelman's monotonicity formula gives $V(t) \geq \frac{1}{C(n,\theta)}$. Thus, $(M,g(t))$ is noncollapsed with a universal constant that depends only on $n$ and $\theta$. \\

\begin{corollary}
\label{pointwise.derivative.estimate.with.universal.constant}
Fix $\theta>0$. We can find a constant $\eta = \eta(n,\theta)$ such that the following holds: Suppose that $(M,g(t))$ is an ancient $\kappa$-solution for some $\kappa>0$, which in addition satisfies $R - \theta \, \text{\rm scal} \, \text{\rm id} \owedge \text{\rm id} \in \text{\rm PIC}$. Then $|D \text{\rm scal}| \leq \eta \, \text{\rm scal}^{\frac{3}{2}}$ and $|\frac{\partial}{\partial t} \text{\rm scal}| \leq \eta \, \text{\rm scal}^2$ at each point in space-time. 
\end{corollary} 

\textbf{Proof.} 
If $(M,g(t))$ is a metric quotient of $S^n$ or $S^{n-1} \times \mathbb{R}$, the assertion is trivial. Otherwise, Theorem \ref{universal.noncollapsing} implies that $(M,g(t))$ is an ancient $\kappa_0$-solution, where $\kappa_0$ depends only on $n$ and $\theta$. Hence, the assertion follows from Corollary \ref{pointwise.derivative.estimate}. \\

\begin{theorem}[cf. G.~Perelman \cite{Perelman2}, Section 1.5]
\label{structure.theorem.general.case}
Given $\varepsilon>0$ and $\theta>0$, there exist positive constants $C_1=C_1(n,\theta,\varepsilon)$ and $C_2=C_2(n,\theta,\varepsilon)$ such that the following holds: Assume that $(M,g(t))$ is an ancient $\kappa$-solution satisfying $R - \theta \, \text{\rm scal} \, \text{\rm id} \owedge \text{\rm id} \in \text{\rm PIC}$. Then, for each point $(x_0,t_0)$ in space-time, there exists a neighborhood $B$ of $x_0$ such that $B_{g(t_0)}(x_0,C_1^{-1} \, \text{\rm scal}(x_0,t_0)^{-\frac{1}{2}}) \subset B \subset B_{g(t_0)}(x_0,C_1 \, \text{\rm scal}(x_0,t_0)^{-\frac{1}{2}})$ and $C_2^{-1} \, \text{\rm scal}(x_0,t_0) \leq \text{\rm scal}(x,t_0) \leq C_2 \, \text{\rm scal}(x_0,t_0)$ for all $x \in B$. Finally, $B$ satisfies one of the following conditions: 
\begin{itemize}
\item $B$ is an $\varepsilon$-neck with center at $x_0$.
\item $B$ is an $\varepsilon$-cap. 
\item $B$ is a closed manifold diffeomorphic to $S^n/\Gamma$.
\item $B$ is an $\varepsilon$-quotient neck of the form $(S^{n-1} \times [-L,L])/\Gamma$.
\end{itemize}
\end{theorem}

\textbf{Proof.} 
If $M$ is noncompact, the assertion follows from Theorem \ref{structure.theorem.noncompact.case}. Hence, it suffices to consider the case when $M$ is compact. As usual, it is enough to consider the case $t_0=0$. Suppose that the assertion is false. Then we can find a sequence of compact ancient $\kappa_j$-solutions $(M^{(j)},g^{(j)}(t))$ satisfying $R - \theta \, \text{\rm scal} \, \text{\rm id} \owedge \text{\rm id} \in \text{\rm PIC}$ and a sequence of points $x_j \in M^{(j)}$ with the following property: There does not exist a neighborhood $B$ of $x_j$ with the property that $B_{g^{(j)}(0)}(x_j,j^{-1} \, \text{\rm scal}(x_j,0)^{-\frac{1}{2}}) \subset B \subset B_{g^{(j)}(0)}(x_j,j \, \text{\rm scal}(x_j,0)^{-\frac{1}{2}})$, $j^{-1} \, \text{\rm scal}(x_j,0) \leq \text{\rm scal}(x,0) \leq j \, \text{\rm scal}(x_j,0)$ for all $x \in B$, and such that $B$ is either an $\varepsilon$-neck with center at $x_j$; or an $\varepsilon$-cap; or a closed manifold diffeomorphic to $S^n/\Gamma$; or an $\varepsilon$-quotient neck. By scaling, we may assume $\text{\rm scal}(x_j,0) = 1$ for each $j$. 

The noncollapsing assumption implies that $(M^{(j)},g^{(j)}(t))$ cannot be isometric to a compact quotient of the standard cylinder. By Corollary \ref{splitting.c}, $(M^{(j)},g^{(j)}(t))$ is strictly PIC2. Clearly, $(M^{(j)},g^{(j)}(t))$ cannot be isometric to a quotient of a round sphere. By Theorem \ref{universal.noncollapsing}, $(M^{(j)},g^{(j)}(t))$ is $\kappa_0$-noncollapsed for some uniform constant $\kappa_0$ which is independent of $j$.

We now apply the compactness theorem for ancient $\kappa_0$-solutions (cf. Corollary \ref{compactness.thm.for.ancient.solutions}) to the sequence $(M^{(j)},g^{(j)}(t),x_j)$. Consequently, after passing to a subsequence if necessary, the sequence $(M^{(j)},g^{(j)}(t),x_j)$ will converge in the Cheeger-Gromov sense to an ancient $\kappa_0$-solution satisfying $R - \theta \, \text{\rm scal} \, \text{\rm id} \owedge \text{\rm id} \in \text{\rm PIC}$. Let us denote this limiting ancient $\kappa_0$-solution by $(M^\infty,g^\infty(t))$, and let $x_\infty$ denote the limit of the sequence $x_j$. There are two possibilities:

\textit{Case 1:} We first consider the case that $M^\infty$ is compact. In this case, the diameter of $(M^{(j)},g^{(j)}(0))$ has a uniform upper bound independent of $j$. Therefore, if $j$ is sufficiently large, then $B^{(j)} := M^{(j)}$ is a neighborhood of the point $x_j$ satisfying $B_{g^{(j)}(0)}(x_j,j^{-1}) \subset B^{(j)} \subset B_{g^{(j)}(0)}(x_j,j)$ and $j^{-1} \leq \text{\rm scal}(x,0) \leq j$ for all $x \in B^{(j)}$. Since $(M^{(j)},g^{(j)}(t))$ is strictly PIC2, results in \cite{Brendle-Schoen} imply that $B^{(j)} = M^{(j)}$ is diffeomorphic to $S^n/\Gamma$. This contradicts our choice of $x_j$.

\textit{Case 2:} We next consider the case that $M^\infty$ is noncompact. If $(M^\infty,g^\infty(t))$ is isometric to a noncompact quotient of the standard cylinder, then, for $j$ large enough, the point $x_j$ lies at the center of an $\varepsilon$-neck or it lies on an $\varepsilon$-quotient neck. This contradicts our choice of $x_j$. Consequently, $(M^\infty,g^\infty(t))$ is not isometric to a quotient of the standard cylinder. At this point, we apply Theorem \ref{structure.theorem.noncompact.case} to $(M^\infty,g^\infty(t))$ (and with $\varepsilon$ replaced by $\frac{\varepsilon}{2}$). Therefore, we can find a neighborhood $B^\infty \subset M^\infty$ of the point $x_\infty$ satisfying $B_{g^\infty(0)}(x_\infty,C_1^{-1}) \subset B^\infty \subset B_{g^\infty(0)}(x_\infty,C_1)$ and $C_2^{-1} \leq \text{\rm scal}(x,0) \leq C_2$ for all $x \in B^\infty$. Furthermore, $B^\infty$ is either an $\frac{\varepsilon}{2}$-neck with center at $x_\infty$ or an $\frac{\varepsilon}{2}$-cap. Hence, if we choose $j$ sufficiently large, then we can find a neighborhood $B^{(j)} \subset M^{(j)}$ of the point $x_j$ satisfying $B_{g^{(j)}(0)}(x_j,(2C_1)^{-1}) \subset B^{(j)} \subset B_{g^{(j)}(0)}(x_j,2C_1)$ and $(2C_2)^{-1} \leq \text{\rm scal}(x,0) \leq 2C_2$ for all $x \in B^{(j)}$. Furthermore, $B^{(j)}$ is either an $\varepsilon$-neck with center at $x_j$ or an $\varepsilon$-cap. This contradicts our choice of $x_j$. \\

For the purpose of the surgery construction, we will need the following refinement of Theorem \ref{structure.theorem.general.case}: 

\begin{corollary}[cf. G.~Perelman \cite{Perelman2}, Section 1.5]
\label{structure.theorem.refined.version}
Given $\varepsilon>0$ and $\theta>0$, there exist positive constants $C_1=C_1(n,\theta,\varepsilon)$ and $C_2=C_2(n,\theta,\varepsilon)$ such that the following holds: Assume that $(M,g(t))$ is an ancient $\kappa$-solution satisfying $R - \theta \, \text{\rm scal} \, \text{\rm id} \owedge \text{\rm id} \in \text{\rm PIC}$. Then, for each point $(x_0,t_0)$ in space-time, there exists a neighborhood $B$ of $x_0$ such that $B_{g(t_0)}(x_0,C_1^{-1} \, \text{\rm scal}(x_0,t_0)^{-\frac{1}{2}}) \subset B \subset B_{g(t_0)}(x_0,C_1 \, \text{\rm scal}(x_0,t_0)^{-\frac{1}{2}})$ and $C_2^{-1} \, \text{\rm scal}(x_0,t_0) \leq \text{\rm scal}(x,t_0) \leq C_2 \, \text{\rm scal}(x_0,t_0)$ for all $x \in B$. Finally, $B$ satisfies one of the following conditions: 
\begin{itemize}
\item $B$ is a strong $\varepsilon$-neck (in the sense of \cite{Perelman2}) with center at $x_0$.
\item $B$ is an $\varepsilon$-cap. 
\item $B$ is a closed manifold diffeomorphic to $S^n/\Gamma$.
\item $B$ is an $\varepsilon$-quotient neck of the form $(S^{n-1} \times [-L,L])/\Gamma$.
\end{itemize}
\end{corollary}

\textbf{Proof.} 
Given $\varepsilon>0$, we can find a positive real number $\tilde{\varepsilon}<\varepsilon$, depending only on $n$, $\theta$, and $\varepsilon$ such that, if $(x_0,t_0)$ lies at the center of an $\tilde{\varepsilon}$-neck, then $(x_0,t_0)$ lies at the center of a strong $\varepsilon$-neck. Hence, the assertion follows from Theorem \ref{structure.theorem.general.case}.

\section{A Canonical Neighborhood Theorem for Ricci flows starting from initial metrics with positive isotropic curvature}

\label{canonical.neighborhood.theorem}

In this section, we consider a solution of the Ricci flow starting from a compact manifold of dimension $n \geq 12$ with positive isotropic curvature. Our goal is to establish an analogue of Perelman's Canonical Neighborhood Theorem. We begin with a definition:

\begin{definition}
\label{f.theta.pinching}
Assume that $f: [0,\infty) \to [0,\infty)$ is a concave and increasing function satisfying $\lim_{s \to \infty} \frac{f(s)}{s} = 0$, and $\theta$ is a positive real number. We say that a Riemannian manifold has $(f,\theta)$-pinched curvature if $R + f(\text{\rm scal}) \, \text{\rm id} \owedge \text{\rm id} \in \text{\rm PIC2}$ and $R - \theta \, \text{\rm scal} \, \text{\rm id} \owedge \text{\rm id} \in \text{\rm PIC}$.
\end{definition}

If $(M,g_0)$ is a compact manifold of dimension $n \geq 12$ with positive isotropic curvature, then Corollary \ref{pinching.estimates.for.ricci.flow} implies that the subsequent solution to the Ricci flow has $(f,\theta)$-pinched curvature for some suitable choice of $f$ and $\theta$. 

\begin{theorem}[cf. G.~Perelman \cite{Perelman1}, Theorem 12.1]
\label{high.curvature.regions.1}
Let $(M,g_0)$ be a compact manifold with positive isotropic curvature of dimension $n \geq 12$, which does not contain any non-trivial incompressible $(n-1)$-dimensional space forms. Let $g(t)$, $t \in [0,T)$, denote the solution to the Ricci flow with initial metric $g_0$. Given a small number $\tilde{\varepsilon}>0$ and a large number $A_0$, we can find a positive number $\hat{r}$ with the following property. If $(x_0,t_0)$ is a point in space-time with $Q := \text{\rm scal}(x_0,t_0) \geq \hat{r}^{-2}$, then the parabolic neighborhood $B_{g(t_0)}(x_0,A_0 Q^{-\frac{1}{2}}) \times [t_0-A_0 Q^{-1},t_0]$ is, after scaling by the factor $Q$, $\tilde{\varepsilon}$-close to the corresponding subset of an ancient $\kappa$-solution satisfying $R - \theta \, \text{\rm scal} \, \text{\rm id} \owedge \text{\rm id} \in \text{\rm PIC}$.
\end{theorem}

The proof of Theorem \ref{high.curvature.regions.1} is an adaptation of Perelman's work \cite{Perelman1} (see also \cite{Chen-Zhu}, where the four-dimensional case is treated). Let us fix $\varepsilon > 0$ and $\theta > 0$, and let $C_1=C_1(n,\theta,\varepsilon)$ and $C_2=C_2(n,\theta,\varepsilon)$ denote the constants in Corollary \ref{structure.theorem.refined.version}. It suffices to prove the assertion when $A_0 \geq 8C_1$ and $\tilde{\varepsilon}$ is much smaller than $\varepsilon$. To do that, we argue by contradiction. If the assertion is false, then we can find a sequence of points $(x_j,t_j)$ in space-time with the following properties: 
\begin{itemize} 
\item[(i)] $Q_j := \text{\rm scal}(x_j,t_j) \geq j^2$. 
\item[(ii)] After dilating by the factor $Q_j$, the parabolic neighborhood $B_{g(t_j)}(x_j,A_0 Q_j^{-\frac{1}{2}}) \times [t_j-A_0 Q_j^{-1},t_j]$ is not $\tilde{\varepsilon}$-close to the corresponding subset of any ancient $\kappa$-solution satisfying $R - \theta \, \text{\rm scal} \, \text{\rm id} \owedge \text{\rm id} \in \text{\rm PIC}$. 
\end{itemize} 
By a point-picking argument, we can arrange that $(x_j,t_j)$ satisfies the following condition: 
\begin{itemize}
\item[(iii)] If $(\tilde{x},\tilde{t})$ is a point in space-time such that $\tilde{t} \leq t_j$ and $\tilde{Q} := \text{\rm scal}(\tilde{x},\tilde{t}) \geq 4Q_j$, then the parabolic neighborhood $B_{g(\tilde{t})}(\tilde{x},A_0 \tilde{Q}^{-\frac{1}{2}}) \times [\tilde{t}-A_0 \tilde{Q}^{-1},\tilde{t}]$ is, after dilating by the factor $\tilde{Q}$, $\tilde{\varepsilon}$-close to the corresponding subset of an ancient $\kappa$-solution satisfying $R - \theta \, \text{\rm scal} \, \text{\rm id} \owedge \text{\rm id} \in \text{\rm PIC}$.
\end{itemize}
Our strategy is to rescale the flow $(M,g(t))$ around the point $(x_j,t_j)$ by the factor $Q_j$. We will show that the rescaled flows converge to an ancient $\kappa$-solution satisfying $R - \theta \, \text{\rm scal} \, \text{\rm id} \owedge \text{\rm id} \in \text{\rm PIC}$. To that end, we proceed in several steps: 

\textit{Step 1:} We first establish a pointwise curvature derivative estimate. By Corollary \ref{pointwise.derivative.estimate.with.universal.constant}, we can find a large constant $\eta$, depending only on $n$ and $\theta$, such that $|D \text{\rm scal}| \leq \eta \, \text{\rm scal}^{\frac{3}{2}}$ and $|\frac{\partial}{\partial t} \text{\rm scal}| \leq \eta \, \text{\rm scal}^2$ on every ancient $\kappa$-solution. Using property (iii) above, we conclude that $|D \text{\rm scal}| \leq 2\eta \, \text{\rm scal}^{\frac{3}{2}}$ and $|\frac{\partial}{\partial t} \text{\rm scal}| \leq 2\eta \, \text{\rm scal}^2$ for each point $(x,t)$ in space-time satisfying $t \leq t_j$ and $\text{\rm scal}(x,t) \geq 4Q_j$. 

\textit{Step 2:} We next prove bounds for the higher order covariant derivatives of the curvature tensor. Suppose that $(x,t)$ is a point in space-time satisfying $\text{\rm scal}(x,t)+Q_j \leq r_0^{-2}$. The pointwise curvature derivative estimate in Step 1 implies that $\text{\rm scal} \leq 8r_0^{-2}$ in the parabolic neighborhood $P(x_0,t_0,\frac{r_0}{100\eta},-\frac{r_0^2}{100\eta})$. Using Shi's interior derivative estimates, we conclude that $|D^m R| \leq C(n,m,\eta) \, r_0^{-m-2}$ at the point $(x_0,t_0)$. 

\textit{Step 3:} We next prove a longrange curvature estimate. Given any $\rho>0$, we define 
\[\mathbb{M}(\rho) = \limsup_{j \to \infty} \sup_{x \in B_{g(t_j)}(x_j,\rho Q_j^{-\frac{1}{2}})} Q_j^{-1} \, \text{\rm scal}(x,t_j).\] 
Here, we allow the possibility that $\mathbb{M}(\rho)=\infty$. The pointwise curvature derivative estimate in Step 1 implies that $\mathbb{M}(\rho) \leq 8$ for $0 < \rho < \frac{1}{100\eta}$. 

We claim that $\mathbb{M}(\rho) < \infty$ for all $\rho>0$. Suppose this is false. Let 
\[\rho^* = \sup \{\rho \geq 0: \mathbb{M}(\rho) < \infty\} < \infty.\] 
Clearly, $\rho^* \geq \frac{1}{100\eta}$. By definition of $\rho^*$, we have an upper bound for the curvature in the geodesic ball $B_{g(t_j)}(x_j,\rho Q_j^{-\frac{1}{2}})$ for each $\rho<\rho^*$. By Step 2, we obtain bounds for the covariant derivatives of the curvature tensor in the geodesic ball $B_{g(t_j)}(x_j,\rho Q_j^{-\frac{1}{2}})$ for each $\rho<\rho^*$. Moreover, Perelman's noncollapsing estimate gives a lower bound for the volume. We rescale around $(x_j,t_j)$ by the factor $Q_j$ and pass to the limit as $j \to \infty$. In the limit, we obtain an incomplete manifold $(B^\infty,g^\infty)$ which is weakly PIC2 (cf. \cite{Morgan-Tian}, Theorem 5.6).

By definition of $\rho^*$, we can find a sequence of points $y_j$ such that $\rho_j := Q_j^{\frac{1}{2}} \, d_{g(t_j)}(x_j,y_j) \to \rho^*$ and $Q_j^{-1} \, \text{\rm scal}(y_j,t_j) \to \infty$. For each $j$, we can find a unit speed geodesic $\gamma_j: [0,\rho_j Q_j^{-\frac{1}{2}}] \to (M,g(t_j))$ such that $\gamma_j(0)=x_j$ and $\gamma_j(\rho_j Q_j^{-\frac{1}{2}}) = y_j$. Let $\gamma_\infty: [0,\rho^*) \to (B^\infty,g^\infty)$ denote the limit of $\gamma_j$. Using the pointwise curvature derivative estimate in Step 1, we obtain 
\[\text{\rm scal}_{g^\infty}(\gamma_\infty(s)) = \lim_{j \to \infty} Q_j^{-1} \, \text{\rm scal}(\gamma_j(sQ_j^{-\frac{1}{2}}),t_j) \geq (\eta(\rho^*-s))^{-2} \geq 100\] 
if $s \in [\rho^*-\frac{1}{100\eta},\rho^*)$.

Let us consider a real number $\bar{s} \in [\rho^*-\frac{1}{100\eta},\rho^*)$ such that $8C_1 \eta (\rho^*-\bar{s}) \leq \bar{s}$. We claim that $\gamma_j(\bar{s}Q_j^{-\frac{1}{2}})$ lies at the center of a $2\varepsilon$-neck if $j$ is sufficiently large (depending on $\bar{s}$). Indeed, if $j$ is sufficiently large, it follows from property (iii) and Corollary \ref{structure.theorem.refined.version} that the point $(\gamma_j(\bar{s}Q_j^{-\frac{1}{2}}),t_j)$ has a Canonical Neighborhood which is either a $2\varepsilon$-neck; or a $2\varepsilon$-cap; or a closed manifold diffeomorphic to $S^n/\Gamma$; or a $2\varepsilon$-quotient neck. Furthermore, the Canonical Neighborhood is contained in a geodesic ball around $\gamma_j(\bar{s}Q_j^{-\frac{1}{2}})$ of radius $2C_1 \, \text{\rm scal}(\gamma_j(\bar{s}Q_j^{-\frac{1}{2}}),t_j)^{-\frac{1}{2}}$, and the scalar curvature is at most $2C_2 \, \text{\rm scal}(\gamma_j(\bar{s}Q_j^{-\frac{1}{2}}),t_j)$ at each point in the Canonical Neighborhood. Since $\mathbb{M}(\bar{s}) < \infty$, we obtain $\lim_{j \to \infty} \text{\rm scal}(\gamma_j(\bar{s}Q_j^{-\frac{1}{2}}),t_j)^{-1} \, \text{\rm scal}(y_j,t_j) = \infty$; consequently, $\text{\rm scal}(y_j,t_j) \geq 4C_2 \, \text{\rm scal}(\gamma_j(\bar{s}Q_j^{-\frac{1}{2}}),t_j)$ if $j$ is sufficiently large. Hence, if $j$ is sufficiently large, then the Canonical Neighborhood does not contain the point $y_j$. We next observe that $8C_1 \, \text{\rm scal}_{g^\infty}(\gamma_\infty(\bar{s}))^{-\frac{1}{2}} \leq 8C_1 \eta (\rho^*-\bar{s}) \leq \bar{s}$. This implies $4C_1 \, \text{\rm scal}(\gamma_j(\bar{s}Q_j^{-\frac{1}{2}}),t_j)^{-\frac{1}{2}} \leq \bar{s} Q_j^{-\frac{1}{2}}$ if $j$ is sufficiently large. Hence, if $j$ is sufficiently large, then the Canonical Neighborhood does not contain the point $x_j$. In particular, if $j$ is sufficiently large, then the Canonical Neighborhood of $(\gamma_j(\bar{s}Q_j^{-\frac{1}{2}}),t_j)$ cannot be a closed manifold diffeomorphic to $S^n/\Gamma$. Moreover, if the Canonical Neighborhood of $(\gamma_j(\bar{s}Q_j^{-\frac{1}{2}}),t_j)$ is a $2\varepsilon$-cap, then the geodesic $\gamma_j$ must enter and exit this $2\varepsilon$-cap, but this is impossible since $\gamma_j$ minimizes length. Finally, if the Canonical Neighborhood of $(\gamma_j(\bar{s}Q_j^{-\frac{1}{2}}),t_j)$ is a quotient neck, then Theorem A.1 in \cite{Brendle4} implies that $M$ contains a non-trivial incompressible $(n-1)$-dimensional space form, contrary to our assumption. To summarize, if $j$ is sufficiently large (depending on $\bar{s}$), then the point $(\gamma_j(\bar{s}Q_j^{-\frac{1}{2}}),t_j)$ has a Canonical Neighborhood which is a $2\varepsilon$-neck. In particular, if $j$ is sufficiently large (depending on $\bar{s}$), then we have $|D \text{\rm scal}| \leq C(n)\varepsilon \, \text{\rm scal}^{\frac{3}{2}}$ at the point $(\gamma_j(\bar{s}Q_j^{-\frac{1}{2}}),t_j)$.

Passing to the limit as $j \to \infty$, we conclude that $|D \text{\rm scal}_{g^\infty}| \leq C(n)\varepsilon \, \text{\rm scal}_{g^\infty}^{\frac{3}{2}}$ at the point $\gamma_\infty(\bar{s})$. Integrating this estimate along $\gamma_\infty$ gives $\text{\rm scal}_{g^\infty}(\gamma_\infty(\bar{s})) \geq (C(n)\varepsilon (\rho^*-\bar{s}))^{-2}$. Moreover, since $(\gamma_j(\bar{s}Q_j^{-\frac{1}{2}}),t_j)$ lies at the center of a $2\varepsilon$-neck for $j$ sufficiently large, the point $\gamma_\infty(\bar{s})$ must lie on a $C(n)\varepsilon$-neck in $(B^\infty,g^\infty)$.

As in \cite{Perelman1}, Section 12.1, there is a sequence of rescalings which converges to a piece of a non-flat metric cone in the limit. Using the pointwise curvature derivative estimate established in Step 1, we can extend the metric backwards in time. This gives a locally defined solution to the Ricci flow which is weakly PIC2 and which, at the final time, is a piece of non-flat metric cone. This contradicts Proposition \ref{splitting.a}.

\textit{Step 4:} We now rescale the manifold $(M,g(t_j))$ around the point $x_j$ by the factor $Q_j$. By Step 3, we have uniform bounds for the curvature at bounded distance. Using the curvature derivative estimate in Step 1 together with Shi's interior derivative estimates, we conclude that the covariant derivatives of the curvature tensor are bounded at bounded distance. Combining this with Perelman's noncollapsing estimate, we conclude that (after passing to a subsequence) the rescaled manifolds converge in the Cheeger-Gromov sense to a complete smooth limit manifold $(M^\infty,g^\infty)$. Since $(M,g(t_j))$ has $(f,\theta)$-pinched curvature, the curvature tensor of $(M^\infty,g^\infty)$ is weakly PIC2 and satisfies $R - \theta \, \text{\rm scal} \, \text{\rm id} \owedge \text{\rm id} \in \text{\rm PIC}$. Using property (iii) and Corollary \ref{structure.theorem.refined.version}, we conclude that every point in $(M^\infty,g^\infty)$ with scalar curvature greater than $4$ has a Canonical Neighborhood which is either a $2\varepsilon$-neck; or a $2\varepsilon$-cap; or a $2\varepsilon$-quotient neck. Note that the last possibility cannot occur; indeed, if $(M^\infty,g^\infty)$ contains a quotient neck, then $(M,g(t_j))$ contains a quotient neck for $j$ sufficiently large, and Theorem A.1 in \cite{Brendle4} then implies that $M$ contains a non-trivial incompressible $(n-1)$-dimensional space form, contrary to our assumption.

We claim that the limit manifold $(M^\infty,g^\infty)$ has bounded curvature. Indeed, if there is a sequence of points in $(M^\infty,g^\infty)$ with curvature going to infinity, then $(M^\infty,g^\infty)$ contains a sequence of necks with radii converging to $0$, contradicting Proposition \ref{no.small.necks}. Thus, $(M^\infty,g^\infty)$ has bounded curvature.

\textit{Step 5:} We now extend the limit $(M^\infty,g^\infty)$ backwards in time. By Step 4, the scalar curvature of $(M^\infty,g^\infty)$ is bounded from above by a constant $\Lambda>4$. Using the pointwise curvature derivative estimate in Step 1, we conclude that 
\[\limsup_{j \to \infty} \sup_{(x,t) \in B_{g(t_j)}(x_j,AQ_j^{-\frac{1}{2}}) \times [t_j-\frac{1}{100\eta\Lambda} \, Q_j^{-1},t_j]} Q_j^{-1} \, \text{\rm scal}(x,t) \leq 2\Lambda\] 
for each $A>1$. Hence, if we put $\tau_1 := -\frac{1}{200\eta \Lambda}$, then we can extend $(M^\infty,g^\infty)$ backwards in time to a complete solution $(M^\infty,g^\infty(t))$ of the Ricci flow which is defined for $t \in [\tau_1,0]$ and satisfies $\Lambda_1 := \sup_{t \in [\tau_1,0]} \sup_{M^\infty} \text{\rm scal}_{g^\infty(t)} \leq 2\Lambda$. In the next step, we put $\tau_2 := \tau_1 - \frac{1}{200\eta \Lambda_1}$. Using the pointwise curvature derivative estimate in Step 1, we can extend the solution $(M^\infty,g^\infty)$, $t \in [\tau_1,0]$, backwards in time to a solution $(M^\infty,g^\infty(t))$, $t \in [\tau_2,0]$. Moreover, $\Lambda_2 := \sup_{t \in [\tau_2,0]} \sup_{M^\infty} \text{\rm scal}_{g^\infty(t)} \leq 2\Lambda_1$. Continuing this process, we can extend the solution backwards in time to the interval $[\tau_m,0]$, where $\tau_{m+1} := \tau_m - \frac{1}{200\eta \Lambda_m}$ and $\Lambda_{m+1} := \sup_{t \in [\tau_{m+1},0]} \sup_{M^\infty} \text{\rm scal}_{g^\infty(t)} \leq 2\Lambda_m$.

Let $\tau^* = \lim_{m \to \infty} \tau_m \leq -\frac{1}{200\eta\Lambda}$. Using a standard diagonal sequence argument, we obtain a complete, smooth limit flow $(M^\infty,g^\infty(t))$ which is defined on the interval $(\tau^*,0]$ and which has bounded curvature for each $t \in (\tau^*,0]$. Since $(M,g(t))$ has $(f,\theta)$-pinched curvature, the curvature tensor of the limit flow $(M^\infty,g^\infty(t))$ is weakly PIC2 and satisfies $R - \theta \, \text{\rm scal} \, \text{\rm id} \owedge \text{\rm id} \in \text{\rm PIC}$. 

\textit{Step 6:} We claim that $\tau^* = -\infty$. To prove this, we argue by contradiction. Suppose that $\tau^* > -\infty$. Clearly, $\lim_{m \to \infty} (\tau_m-\tau_{m+1}) = 0$, hence $\lim_{m \to \infty} \Lambda_m = \infty$. 

By the Harnack inequality (cf. Theorem \ref{harnack} above), the function $t \mapsto (t-\tau^*) \, \text{\rm scal}_{g^\infty(t)}(x)$ is monotone increasing at each point $x \in M^\infty$. Since $\text{\rm scal}_{g^\infty(0)}(x) \leq \Lambda$ for all $x \in M^\infty$, we obtain 
\[\text{\rm scal}_{g^\infty(t)}(x) \leq \frac{-\tau^*}{t-\tau^*} \, \Lambda\] 
for all $x \in M^\infty$ and all $t \in (\tau^*,0]$. Using Lemma 8.3(b) in \cite{Perelman1}, we conclude that 
\[0 \leq -\frac{d}{dt} d_{g^\infty(t)}(x,y) \leq C(n) \, \sqrt{\frac{-\tau^*}{t-\tau^*} \, \Lambda}\] 
for all $x,y \in M^\infty$ and all $t \in (\tau^*,0]$. Integrating over $t$ gives 
\[d_{g^\infty(0)}(x,y) \leq d_{g^\infty(t)}(x,y) \leq d_{g^\infty(0)}(x,y) + C(n) \, (-\tau^*) \, \sqrt{\Lambda} \] 
for all $x,y \in M^\infty$ and all $t \in (\tau^*,0]$. 

By the maximum principle, 
\[\inf_{M^\infty} \text{\rm scal}_{g^\infty(t)} \leq \inf_{M^\infty} \text{\rm scal}_{g^\infty(0)} \leq 1\] 
for all $t \in (\tau^*,0]$. Hence, we can find a point $y_\infty \in M^\infty$ such that $\text{\rm scal}_{g^\infty(t)}(y_\infty) \leq 4$ for $t = \tau^* + \frac{1}{1000\eta\Lambda} \in (\tau^*,0]$. Using the pointwise curvature derivative estimate in Step 1, we obtain $\text{\rm scal}_{g^\infty(t)}(y_\infty) \leq 8$ for all $t \in (\tau^*,\tau^*+\frac{1}{1000\eta\Lambda}]$. In particular, $\text{\rm scal}_{g^\infty(\tau_m)}(y_\infty) \leq 8$ if $m$ is sufficiently large. Arguing as in Step 3, we can show that 
\[\limsup_{m \to \infty} \sup_{B_{g^\infty(\tau_m)}(y_\infty,A)} \text{\rm scal}_{g^\infty(\tau_m)} < \infty\] 
for every $A>1$. Consequently, a subsequence of the manifolds $(M^\infty,g^\infty(\tau_m),y_\infty)$ converges in the Cheeger-Gromov sense to a complete, smooth limit. If this limit manifold has unbounded curvature, then (by property (iii) above) it contains a sequence of necks with radii converging to $0$, contradicting Proposition \ref{no.small.necks}. Therefore, a subsequence of the manifolds $(M^\infty,g^\infty(\tau_m),y_\infty)$ converges in the Cheeger-Gromov sense to a complete, smooth limit with bounded curvature. Consequently, we can find a constant $\Lambda^*>\Lambda$ (independent of $A$) such that 
\[\liminf_{m \to \infty} \sup_{B_{g^\infty(\tau_m)}(y_\infty,A)} \text{\rm scal}_{g^\infty(\tau_m)} \leq \Lambda^*\]
for every $A>1$. Using the distance estimate, we obtain $B_{g^\infty(0)}(y_\infty,A) \subset B_{g^\infty(\tau_m)}(y_\infty,A + C(n) \, (-\tau^*) \, \sqrt{\Lambda})$. Putting these facts together, we conclude that 
\[\liminf_{m \to \infty} \sup_{B_{g^\infty(0)}(y_\infty,A)} \text{\rm scal}_{g^\infty(\tau_m)} \leq \Lambda^*\]
for every $A>1$. Hence, for each $A>1$, we can find a large integer $m$ (depending on $A$) such that $\tau_m \in (\tau^*,\tau^*+\frac{1}{1000\eta\Lambda^*}]$ and 
\[\sup_{B_{g^\infty(0)}(y_\infty,A)} \text{\rm scal}_{g^\infty(\tau_m)} \leq 2\Lambda^*.\]
Using the pointwise derivative estimate in Step 1, we obtain 
\[\sup_{t \in (\tau^*,\tau^*+\frac{1}{1000\eta\Lambda^*}]} \sup_{B_{g^\infty(0)}(y_\infty,A)} \text{\rm scal}_{g^\infty(t)} \leq 4\Lambda^*\] 
for every $A>1$. Since $\Lambda^*$ is independent of $A$, we conclude that 
\[\sup_{t \in (\tau^*,\tau^*+\frac{1}{1000\eta\Lambda^*}]} \sup_{M^\infty} \text{\rm scal}_{g^\infty(t)} \leq 4\Lambda^*.\] 
Therefore, the flow $(M^\infty,g^\infty(t))$, $t \in (\tau^*,0]$, has bounded curvature. This contradicts the fact that $\lim_{m \to \infty} \Lambda_m = \infty$. Thus, $\tau^* = -\infty$.

To summarize, if we dilate the flow $(M,g(t))$ around the point $(x_j,t_j)$ by the factor $Q_j$, then (after passing to a subsequence), the rescaled flows converge in the Cheeger-Gromov sense to an ancient $\kappa$-solution $(M^\infty,g^\infty(t))$, $t \in (-\infty,0]$, satisfying $R - \theta \, \text{\rm scal} \, \text{\rm id} \owedge \text{\rm id} \in \text{\rm PIC}$. Here, $\kappa$ depends only on the initial data. This contradicts statement (ii). This completes the proof of Theorem \ref{high.curvature.regions.1}. \\

Finally, by combining Theorem \ref{high.curvature.regions.1} with Theorem \ref{structure.theorem.refined.version}, we can draw the following conclusion:

\begin{corollary}[cf. G.~Perelman \cite{Perelman1}, Theorem 12.1]
\label{high.curvature.regions.2}
Let $(M,g_0)$ be a compact manifold with positive isotropic curvature of dimension $n \geq 12$, which does not contain any non-trivial incompressible $(n-1)$-dimensional space forms. Let $g(t)$, $t \in [0,T)$, denote the solution to the Ricci flow with initial metric $g_0$. Given any $\varepsilon>0$, there exists a positive number $\hat{r}$ with the following property. If $(x_0,t_0)$ is a point in space-time with $Q := \text{\rm scal}(x_0,t_0) \geq \hat{r}^{-2}$, then we can find a neighborhood $B$ of $x_0$ such that $B_{g(t_0)}(x_0,(2C_1)^{-1} \, \text{\rm scal}(x_0,t_0)^{-\frac{1}{2}}) \subset B \subset B_{g(t_0)}(x_0,2C_1 \, \text{\rm scal}(x_0,t_0)^{-\frac{1}{2}})$ and $(2C_2)^{-1} \, \text{\rm scal}(x_0,t_0) \leq \text{\rm scal}(x,t_0) \leq 2C_2 \, \text{\rm scal}(x_0,t_0)$ for all $x \in B$. Furthermore, $B$ satisfies one of the following conditions:
\begin{itemize}
\item $B$ is a strong $2\varepsilon$-neck (in the sense of \cite{Perelman2}) with center at $x_0$.
\item $B$ is a $2\varepsilon$-cap in the sense of Definition \ref{definition.of.cap}. 
\item $B$ is a closed manifold diffeomorphic to $S^n/\Gamma$.
\end{itemize}
Here, $C_1=C_1(n,\theta,\varepsilon)$ and $C_2=C_2(n,\theta,\varepsilon)$ are the constants appearing in Corollary \ref{structure.theorem.refined.version}. Finally, we have $|D \text{\rm scal}| \leq 2\eta \, \text{\rm scal}^{\frac{3}{2}}$ and $|\frac{\partial}{\partial t} \text{\rm scal}| \leq 2\eta \, \text{\rm scal}^2$ at the point $(x_0,t_0)$, where $\eta$ is a constant that depends only on $n$ and $\theta$.
\end{corollary}

\section{The behavior of the flow at the first singular time}

\label{first.singular.time}

Throughout this section, we fix a compact initial manifold $(M,g_0)$ of dimension $n \geq 12$ which has positive isotropic curvature and does not contain any non-trivial incompressible $(n-1)$-dimensional space forms. Let $(M,g(t))$ be the solution of the Ricci flow with initial metric $g_0$, and let $[0,T)$ denote the maximal time interval on which the solution is defined. Note that $T \leq \frac{n}{2 \, \inf_{x \in M} \text{\rm scal}(x,0)}$. By Theorem \ref{main.theorem}, we can find a continuous family of closed, convex, $O(n)$-invariant sets $\{\mathcal{F}_t: t \in [0,T]\}$ such that the family $\{\mathcal{F}_t: t \in [0,T]\}$ is invariant under the Hamilton ODE $\frac{d}{dt} R = Q(R)$; the curvature tensor of $(M,g_0)$ lies in the set $\mathcal{F}_0$; and 
\begin{align*} 
\mathcal{F}_t &\subset \{R: R - \theta \, \text{\rm scal} \, \text{\rm id} \owedge \text{\rm id} \in \text{\rm PIC}\} \\ 
&\cap \{R: R + f(\text{\rm scal}) \, \text{\rm id} \owedge \text{\rm id} \in \text{\rm PIC2}\}
\end{align*} 
for all $t \in [0,T]$. Here, $f$ is a concave and increasing function satisfying $\lim_{s \to \infty} \frac{f(s)}{s} = 0$, and $\theta$ and $N$ are positive numbers. Note that $f$, $\theta$, and $N$ depend only on the initial data. By Hamilton's PDE-ODE principle (cf. \cite{Chow-Lu}, Theorem 3, or \cite{Chow-et-al}, Theorem 10.16), the curvature tensor of $(M,g(t))$ lies in the set $\mathcal{F}_t$ for each $t \in [0,T)$. 

By Corollary \ref{high.curvature.regions.2}, every point in space-time where the scalar curvature is sufficiently large admits a Canonical Neighborhood which is either a $2\varepsilon$-neck; or a $2\varepsilon$-cap; or a closed manifold diffeomorphic to $S^n/\Gamma$. Let $\rho$ be a small positive number with the property that every point with $\text{\rm scal} \geq \frac{1}{4} \, \rho^{-2}$ satisfies the conclusion of the Canonical Neighborhood Theorem. In particular, we have $|D \text{\rm scal}| \leq 2\eta \, \text{\rm scal}^{\frac{3}{2}}$ and $|\frac{\partial}{\partial t} \text{\rm scal}| \leq 2\eta \, \text{\rm scal}^2$ whenever $\text{\rm scal} \geq \frac{1}{4} \, \rho^{-2}$. We define 
\[\Omega := \{x \in M: \limsup_{t \to T} \text{\rm scal}(x,t) < \infty\}.\] 
The pointwise curvature derivative estimate implies that $\Omega$ is an open subset of $M$. Using the pointwise curvature derivative estimate together with Shi's interior estimates, we conclude the metrics $g(t)$ converge to a smooth metric $g(T)$ on $\Omega$. Following \cite{Perelman1}, we consider the set 
\[\Omega_\rho := \{x \in M: \limsup_{t \to T} \text{\rm scal}(x,t) \leq \rho^{-2}\} = \{x \in \Omega: \text{\rm scal}(x,T) \leq \rho^{-2}\}.\] 
We distinguish two cases:

\textit{Case 1:} Suppose that $\Omega_\rho = \emptyset$. Using the inequality $|\frac{\partial}{\partial t} \text{\rm scal}| \leq 2\eta \, \text{\rm scal}^2$, we obtain $\inf_{x \in M} \text{\rm scal}(x,t) \geq \frac{1}{2} \, \rho^{-2}$ if $t$ is sufficiently close to $T$. Hence, if $t$ is sufficiently close to $T$, then every point in $(M,g(t))$ admits a Canonical Neighborhood which is either a $2\varepsilon$-neck; or a $2\varepsilon$-cap; or a closed manifold diffeomorphic to $S^n/\Gamma$. 

\textit{Case 2:} Suppose now that $\Omega_\rho \neq \emptyset$. The Canonical Neighborhood Theorem guarantees that every point in $\Omega \setminus \Omega_\rho$ lies either on a $2\varepsilon$-tube with boundary components in $\Omega_\rho$; or on a $2\varepsilon$-cap with boundary in $\Omega_\rho$; or on a $2\varepsilon$-horn with boundary in $\Omega_\rho$; or on a double $2\varepsilon$-horn; or on a capped $2\varepsilon$-horn; or on a closed manifold diffeomorphic to $S^n/\Gamma$. (Here, we use the definitions from Perelman's paper \cite{Perelman2}.) Following Perelman \cite{Perelman2}, we perform surgery on each $2\varepsilon$-horn with boundary in $\Omega_\rho$. We discard all double $2\varepsilon$-horns, all capped $2\varepsilon$-horns, and all closed manifolds diffeomorphic to $S^n/\Gamma$. We leave unchanged all the $2\varepsilon$-tubes with boundary in $\Omega_\rho$, and all $2\varepsilon$-caps with boundary in $\Omega_\rho$. \\

\begin{proposition}
\label{change.of.topology}
The pre-surgery manifold $M$ is diffeomorphic to a connected sum of the post-surgery manifold with a finite collection of standard spaces, each of which is a quotient of $S^n$ or $S^{n-1} \times \mathbb{R}$ by standard isometries.
\end{proposition}

\textbf{Proof.} 
Suppose first that $\Omega_\rho = \emptyset$. In this case, $M$ is diffeomorphic to either a quotient of $S^n$ by standard isometries; or a tube with caps attached on both sides; or an $S^{n-1}$-bundle over $S^1$ with a fiberwise round metric. In the second case, Definition \ref{definition.of.cap} ensures that $M$ is diffeomorphic to $S^n$. To handle the third case, we note that that there are two $S^{n-1}$-bundles over $S^1$ with a fiberwise round metric. One of them is orientable, the other one is not. Both are diffeomorphic to quotients of $S^{n-1} \times \mathbb{R}$ by standard isometries. To summarize, $M$ is diffeomorphic to a quotient of $S^n$ or a quotient of $S^{n-1} \times \mathbb{R}$ by standard isometries.

Suppose next that $\Omega_\rho \neq \emptyset$. In this case, we can recover the pre-surgery manifold $M$ from the post-surgery manifold as follows. We first reinstate the components that were discarded after surgery. More precisely, we form a disjoint union of the post-surgery manifold and a finite collection of standard spaces, each of which is a quotient of $S^n$ or $S^{n-1} \times \mathbb{R}$ by standard isometries. In the next step, we reverse the surgery by gluing in finitely many handles of the form $S^{n-1} \times I$. Note that, as we glue in these handles, the attaching maps are nearly isometric. Thus, the pre-surgery manifold is diffeomorphic to a connected sum of the post-surgery manifold with finitely many quotients of $S^n$ and $S^{n-1} \times \mathbb{R}$. This completes the proof of Proposition \ref{change.of.topology}. \\

In the remainder of this section, we show that the surgery procedure preserves our curvature pinching estimates, provided that the surgery parameters are sufficiently fine.

\begin{proposition}
\label{curvature.under.surgery}
Suppose that the curvature tensor of a $\delta$-neck lies in the set $\mathcal{F}_t$ prior to surgery. If $\delta$ is sufficiently small and the curvature of the neck is sufficiently large, then the curvature tensor of the surgically modified manifold lies in the set $\mathcal{F}_t$. Moreover, the scalar curvature is pointwise increasing under surgery.
\end{proposition}

\textbf{Proof.} 
Suppose that the scalar curvature of the neck is close to $h^{-2}$, where $h$ is small. Let us rescale by the factor $h^{-1}$ so that the scalar curvature of the neck is close to $1$ after rescaling. Let us, therefore, assume that $g$ is a Riemannian metric on $S^{n-1} \times [-10,10]$ which is close to the round metric with scalar curvature $1$, and which has curvature in the set $h^2 \mathcal{F}_t$.  We first recall the definition of the surgically modified metric $\tilde{g}$ (cf. \cite{Hamilton5}, Section 4.1). To that end, let $z$ denote the height function on $S^{n-1} \times [-10,10]$, and let $\varphi(z) = e^{-\frac{1}{z}}$ for $z \in (0,\frac{1}{10}]$. In the region $S^{n-1} \times [-10,0]$, the metric is unchanged under surgery, i.e. $\tilde{g} = g$. In the region $S^{n-1} \times (0,\frac{1}{20}]$, we change the metric conformally by $\tilde{g} = e^{-2\varphi} \, g$. In the region $S^{n-1} \times (\frac{1}{20},\frac{1}{10}]$, we define $\tilde{g} = e^{-2\varphi} \, (\chi(z) \, g + (1-\chi(z)) \, \bar{g})$, where $\bar{g}$ denotes the standard metric on the cylinder and $\chi: (\frac{1}{20},\frac{1}{10}] \to [0,1]$ is a smooth cutoff function satisfying $\chi(z)=1$ for $z \in [\frac{1}{20},\frac{1}{18}]$ and $\chi(z)=0$ for $z \in [\frac{1}{12},\frac{1}{10}]$. In particular, the surgically modified metric $\tilde{g}$ is rotationally symmetric for $z \in [\frac{1}{12},\frac{1}{10}]$. Hence, we may extend $\tilde{g}$ by gluing in a rotationally symmetric cap.

We now analyze the curvature of the surgically modified metric $\tilde{g}$. It suffices to consider the case when $z>0$ is small. In this region, $\tilde{g} = e^{-2\varphi} g$. Let $\{e_1,\hdots,e_n\}$ denote a local orthonormal frame with respect to the metric $g$. If we put $\tilde{e}_i = e^\varphi e_i$, then $\{\tilde{e}_1,\hdots,\tilde{e}_n\}$ is an orthonormal frame with respect to the metric $\tilde{g}$. We will express geometric quantities associated with the metric $g$ relative to the frame $\{e_1,\hdots,e_n\}$, while geometric quantities associated with $\tilde{g}$ will be expressed in terms of $\{\tilde{e}_1,\hdots,\tilde{e}_n\}$. With this understood, the curvature tensor after surgery is related to the curvature tensor before surgery by the formula 
\[\tilde{R} = e^{2\varphi} \, R + e^{2\varphi} \, \Big ( D^2 \varphi + d\varphi \otimes d\varphi - \frac{1}{2} \, |d\varphi|^2 \, \text{\rm id} \Big ) \owedge \text{\rm id}.\] 
This implies 
\[|\tilde{R} - R - z^{-4} \, e^{-\frac{1}{z}} \, (dz \otimes dz) \owedge \text{\rm id}| \ll z^{-4} \, e^{-\frac{1}{z}}\] 
for $z>0$ sufficiently small. Consequently, $\text{\rm scal}(\tilde{R}) > \text{\rm scal}(R)$ if $z>0$ is sufficiently small. Since the metric $g$ is close to the cylindrical metric, we obtain
\[\Big | R - \frac{1}{2} \, (\text{\rm id} - 2 \, z \otimes z) \owedge \text{\rm id} \Big | \ll 1,\] 
hence 
\begin{align*} 
&\Big | \tilde{R} - (1-z^{-4} \, e^{-\frac{1}{z}}) \, R - \frac{1}{2} \, z^{-4} \, e^{-\frac{1}{z}} \, \text{\rm id} \owedge \text{\rm id} \Big | \\ 
&\leq |\tilde{R} - R - z^{-4} \, e^{-\frac{1}{z}} \, (dz \otimes dz) \owedge \text{\rm id}| + z^{-4} \, e^{-\frac{1}{z}} \, \Big | R - \frac{1}{2} \, (\text{\rm id} - 2 \, z \otimes z) \owedge \text{\rm id} \Big | \\ 
&\ll z^{-4} \, e^{-\frac{1}{z}} 
\end{align*} 
for $z>0$ sufficiently small. Therefore, we may write 
\[\tilde{R} = (1-z^{-4} \, e^{-\frac{1}{z}}) \, R + z^{-4} \, e^{-\frac{1}{z}} \, S,\] 
where $|S - \frac{1}{2} \, \text{\rm id} \owedge \text{\rm id}| \ll 1$ for $z>0$ sufficiently small. Consequently, $S \in h^2 \mathcal{F}_t$ if $z>0$ is sufficiently small. Moreover, $R \in h^2 \mathcal{F}_t$ in view of our assumption. Since $\mathcal{F}_t$ is a convex set, we conclude that $\tilde{R} \in h^2 \mathcal{F}_t$ if $z>0$ is sufficiently small. This easily implies that $\tilde{R} \in h^2 \mathcal{F}_t$ for all $z \in (0,10)$.

\section{The standard solution}

\label{standard.solution}

In this section, we recall some basic facts concerning the so-called standard solution. The standard solution is used to model the evolution of a cap that is glued in during a surgery procedure. More precisely, suppose that $(S^{n-1} \times \mathbb{R},g(t))$, $t<0$, is a family of shrinking cylinders, normalized so that $\text{\rm scal}_{g(t)} = \frac{1}{1-\frac{2t}{n-1}}$ for $t<0$. Suppose that we perform surgery at time $t=0$: that is, we remove a half-cylinder and glue in a cap which is rotationally symmetric and has positive curvature. This gives a rotationally symmetric metric $g(0)$ on $\mathbb{R}^n$. The standard solution is obtained by evolving the manifold $(\mathbb{R}^n,g(0))$ under the Ricci flow.

The following results were proved by Perelman \cite{Perelman2} in dimension $3$ and were extended to higher dimensions in \cite{Chen-Zhu}.

\begin{theorem}[G.~Perelman \cite{Perelman2}, Section 2; B.L.~Chen, X.P.~Zhu \cite{Chen-Zhu}, Theorem A.1]
\label{properties.of.standard.solution}
There exists a complete solution $(\mathbb{R}^n,g(t))$, $t \in [0,\frac{n-1}{2})$, to the Ricci flow with the following properties: 
\begin{itemize}
\item[(i)] The initial manifold $(\mathbb{R}^n,g(0))$ is isometric to a standard cylinder with scalar curvature $1$ outside of a compact set, and this compact set is isometric to the cap used in the surgery procedure. 
\item[(ii)] For each $t \in [0,\frac{n-1}{2})$, the manifold $(\mathbb{R}^n,g(t))$ is rotationally symmetric. 
\item[(iii)] For each $t \in [0,\frac{n-1}{2})$, the manifold $(\mathbb{R}^n,g(t))$ is asymptotic to a cylinder with scalar curvature $\frac{1}{1-\frac{2t}{n-1}}$ at infinity. 
\item[(iv)] The scalar curvature is bounded from below by $\frac{1}{K_{\text{\rm stnd}} \, (1-\frac{2t}{n-1})}$, where $K_{\text{\rm stnd}}$ is a positive constant that depends only on $n$. 
\item[(v)] For each $t \in [0,\frac{n-1}{2})$, the manifold $(\mathbb{R}^n,g(t))$ is weakly PIC2 and satisfies $R - \theta \, \text{\rm scal} \, \text{\rm id} \owedge \text{\rm id} \in \text{\rm PIC}$ for some constant $\theta>0$ which depends only on $n$. 
\item[(vi)] The flow $(\mathbb{R}^n,g(t))$ is $\kappa$-noncollapsed for some constant $\kappa>0$ which depends only on $n$.
\item[(vii)] There exists a function $\omega: [0,\infty) \to (0,\infty)$ such that $\text{\rm scal}(x,t) \leq \text{\rm scal}(y,t) \, \omega(\text{\rm scal}(y,t) \, d_{g(t)}(x,y)^2)$ for all points $x,y$ and all $t \in [0,\frac{n-1}{2})$.
\end{itemize}
\end{theorem}

\textbf{Proof.} 
The statements (i), (ii), (iii), (iv), (vi), (vii) are established in \cite{Chen-Zhu}, Appendix A. Moreover, it is shown in \cite{Chen-Zhu} that $(\mathbb{R}^n,g(t))$ has nonnegative curvature operator. Hence, it remains to show that $R - \theta \, \text{\rm scal} \, \text{\rm id} \owedge \text{\rm id} \in \text{\rm PIC}$. To see this, we observe that the initial manifold $(\mathbb{R}^n,g(0))$ is uniformly PIC. Moreover, on the initial manifold $(\mathbb{R}^n,g(0))$, the sum of the two smallest eigenvalues of the Ricci tensor is bounded from below by a small multiple of the scalar curvature. Hence, we can find a small constant $b \in (0,b_{\text{\rm max}})$ such that the curvature tensor of $(\mathbb{R}^n,g(0))$ lies in the cone $\mathcal{C}(b)$. By Hamilton's PDE-ODE principle (see \cite{Chow-et-al}, Theorem 12.34), the curvature tensor of $(\mathbb{R}^n,g(t))$ lies in $\mathcal{C}(b)$ for each $t \geq 0$. Consequently, the curvature tensor of $(\mathbb{R}^n,g(t))$ satisfies $R - \theta \, \text{\rm scal} \, \text{\rm id} \owedge \text{\rm id} \in \text{\rm PIC}$ for each $t \geq 0$. \\

It turns out that the standard solution satisfies a Canonical Neighborhood Property: 

\begin{theorem}[cf. G.~Perelman \cite{Perelman2}; B.L.~Chen, X.P.~Zhu \cite{Chen-Zhu}, Corollary A.2]
\label{canonical.neighborhood.theorem.for.standard.solution}
Given a small number $\tilde{\varepsilon}>0$ and a large number $A_0$, we can find a number $\alpha \in [0,\frac{n-1}{2})$ with the following property. If $(x_0,t_0)$ is a point on the standard solution such that $t_0 \in [\alpha,\frac{n-1}{2})$, then the parabolic neighborhood $P(x_0,t_0,A_0 \, \text{\rm scal}(x_0,t_0)^{-\frac{1}{2}},-A_0 \, \text{\rm scal}(x_0,t_0)^{-1})$ is, after scaling by the factor $\text{\rm scal}(x_0,t_0)$, $\tilde{\varepsilon}$-close to the corresponding subset of a noncompact ancient $\kappa_0$-solution satisfying $R - \theta \, \text{\rm scal} \, \text{\rm id} \owedge \text{\rm id} \in \text{\rm PIC}$.
\end{theorem}

\textbf{Proof.} 
Suppose that the assertion is false. Then we can find a sequence of points $(x_j,t_j)$ on the standard solution such that $t_j \to \frac{n-1}{2}$ and the parabolic neighborhood $P(x_j,t_j,A_0 \, \text{\rm scal}(x_j,t_j)^{-\frac{1}{2}},-A_0 \, \text{\rm scal}(x_j,t_j)^{-1})$ is not $\tilde{\varepsilon}$-close to the corresponding subset of a noncompact ancient $\kappa_0$-solution satisfying $R - \theta \, \text{\rm scal} \, \text{\rm id} \owedge \text{\rm id} \in \text{\rm PIC}$. We dilate the solution around the point $(x_j,t_j)$ by the factor $\text{\rm scal}(x_j,t_j)$. Using statement (vii) in Theorem \ref{properties.of.standard.solution} together with the Harnack inequality (cf. Theorem \ref{harnack}), we conclude that the rescaled flows converge to a complete, noncompact ancient solution $(M^\infty,g^\infty(t))$. The limiting ancient solution $(M^\infty,g^\infty(t))$ is weakly PIC2 and satisfies $R - \theta \, \text{\rm scal} \, \text{\rm id} \owedge \text{\rm id} \in \text{\rm PIC}$. Moreover, the limiting ancient solution is $\kappa_0$-noncollapsed. 

By Theorem \ref{harnack}, the standard solution satisfies the Harnack inequality 
\[\frac{\partial}{\partial t} \text{\rm scal} + 2 \, \langle \nabla \text{\rm scal},v \rangle + 2 \, \text{\rm Ric}(v,v) + \frac{1}{t} \, \text{\rm scal} \geq 0\] 
for $t \in (0,\frac{n-1}{2})$. Consequently, the limiting ancient solution $(M^\infty,g^\infty(t))$ satisfies 
\[\frac{\partial}{\partial t} \text{\rm scal} + 2 \, \langle \nabla \text{\rm scal},v \rangle + 2 \, \text{\rm Ric}(v,v) \geq 0.\] 
Using Proposition \ref{boundedness.of.curvature}, we conclude that $(M^\infty,g^\infty(t))$ has bounded curvature. Thus, $(M^\infty,g^\infty(t))$ is a noncompact ancient $\kappa_0$-solution satisfying $R - \theta \, \text{\rm scal} \, \text{\rm id} \owedge \text{\rm id} \in \text{\rm PIC}$. This is a contradiction. \\

\begin{corollary}[cf. G.~Perelman \cite{Perelman2}; B.L.~Chen, X.P.~Zhu \cite{Chen-Zhu}, Corollary A.2]
\label{canonical.neighborhood.theorem.for.standard.solution.2}
Given $\varepsilon>0$, there exist positive constants $C_1=C_1(n,\varepsilon)$ and $C_2=C_2(n,\varepsilon)$ such that the following holds: For each point $(x_0,t_0)$ on the standard solution, there exists a neighborhood $B$ of $x_0$ such that $B_{g(t_0)}(x_0,C_1^{-1} \, \text{\rm scal}(x_0,t_0)^{-\frac{1}{2}}) \subset B \subset B_{g(t_0)}(x_0,C_1 \, \text{\rm scal}(x_0,t_0)^{-\frac{1}{2}})$ and $C_2^{-1} \, \text{\rm scal}(x_0,t_0) \leq \text{\rm scal}(x,t_0) \leq C_2 \, \text{\rm scal}(x_0,t_0)$ for all $x \in B$. Furthermore, $B$ satisfies one of the following conditions: 
\begin{itemize}
\item $B$ is a strong $\varepsilon$-neck (in the sense of \cite{Perelman2}) with center at $x_0$. In particular, if $t_0-R(x_0,t_0)^{-1} \leq 0$, then $B$ is disjoint from the surgical cap that was glued in at time $0$.
\item $B$ is an $\varepsilon$-cap in the sense of Definition \ref{definition.of.cap}. 
\end{itemize}
Finally, we have $|D \text{\rm scal}| \leq \eta \, \text{\rm scal}^{\frac{3}{2}}$ and $|\frac{\partial}{\partial t} \text{\rm scal}| \leq \eta \, \text{\rm scal}^2$.
\end{corollary}

\textbf{Proof.} If $t_0$ is sufficiently close to $\frac{n-1}{2}$ (depending on $\varepsilon$), this follows from Theorem \ref{canonical.neighborhood.theorem.for.standard.solution} and Theorem \ref{structure.theorem.noncompact.case}. If $t_0$ is bounded away from $\frac{n-1}{2}$, this follows from the fact that the standard solution is asymptotic to a cylinder at infinity. \\

Finally, we state a lemma which will be needed later. 

\begin{lemma}
\label{no.escape}
Given $\alpha \in [0,\frac{n-1}{2})$ and $l>0$, we can find a large number $A$ (depending on $\alpha$ and $l$) with the following property. Suppose that $t_1 \in [0,\alpha]$ and $\gamma$ is a space-time curve on the standard solution (parametrized by the interval $[0,t_1]$) such that $\gamma(0)$ lies on the cap at time $0$, and $\int_0^{t_1} |\gamma'(t)|_{g(t)}^2 \, dt \leq l$. Then the curve $\gamma$ is contained in the parabolic neighborhood $P(\gamma(0),0,\frac{A}{2},t_1)$.
\end{lemma}

\textbf{Proof.} Using the inequality $\int_0^{t_1} |\gamma'(t)|_{g(t)}^2 \, dt \leq l$ and H\"older's inequality, we obtain $\int_0^{t_1} |\gamma'(t)|_{g(t)} \, dt \leq \alpha^{\frac{1}{2}} \, l^{\frac{1}{2}}$. From this, the assertion follows easily.

\section{A priori estimates for Ricci flow with surgery}

\label{a.priori.estimates.under.surgery}

In this section, we give the definition of Ricci flow with surgery. Moreover, we discuss how Perelman's noncollapsing estimate and the Canonical Neighborhood Theorem can be extended to Ricci flow with surgery. Throughout this section, we fix a compact initial manifold $(M,g_0)$ of dimension $n \geq 12$ which has positive isotropic curvature and does not contain any non-trivial incompressible $(n-1)$-dimensional space forms. As above, let $\{\mathcal{F}_t: t \in [0,T]\}$ be a family of closed, convex, $O(n)$-invariant sets such that the family $\{\mathcal{F}_t: t \in [0,T]\}$ is invariant under the Hamilton ODE $\frac{d}{dt} R = Q(R)$; the curvature tensor of $(M,g_0)$ lies in the set $\mathcal{F}_0$; and 
\begin{align*} 
\mathcal{F}_t &\subset \{R: R - \theta \, \text{\rm scal} \, \text{\rm id} \owedge \text{\rm id} \in \text{\rm PIC}\} \\ 
&\cap \{R: R + f(\text{\rm scal}) \, \text{\rm id} \owedge \text{\rm id} \in \text{\rm PIC2}\}
\end{align*} 
for all $t \in [0,T]$. Here, $f$ is a concave and increasing function satisfying $\lim_{s \to \infty} \frac{f(s)}{s} = 0$, and $\theta$ and $N$ are positive numbers. 

Having fixed $\theta$, we can find a universal constant $\kappa_0$ such that the conclusion of Theorem \ref{universal.noncollapsing} holds. \textbf{Morevoer, we fix a constant $\eta$ such that the conclusions of Corollary \ref{pointwise.derivative.estimate.with.universal.constant} and Corollary \ref{canonical.neighborhood.theorem.for.standard.solution.2} hold.} In other words, we have $|D \text{\rm scal}| \leq \eta \, \text{\rm scal}^{\frac{3}{2}}$ and $|\frac{\partial}{\partial t} \text{\rm scal}| \leq \eta \, \text{\rm scal}^2$ on any ancient $\kappa$-solution, and the same inequalities hold on the standard solution. 

Let us fix a small positive number $\varepsilon>0$. \textbf{Moreover, we fix constants $C_1=C_1(n,\theta,\varepsilon)$ and $C_2=C_2(n,\theta,\varepsilon)$ such that the conclusions of Corollary \ref{structure.theorem.refined.version} and Corollary \ref{canonical.neighborhood.theorem.for.standard.solution.2} hold.} 

\begin{definition}
\label{definition.ricci.flow.with.surgery}
A Ricci flow with surgery on the interval $[0,T)$ consists of the following data: 
\begin{itemize}
\item A decomposition of $[0,T)$ into a disjoint union of finitely many subintervals $[t_k^-,t_k^+)$, $k \in \{0,1,\hdots,l\}$. In other words, $t_0^-=0$, $t_l^+=T$, and $t_k^- = t_{k-1}^+$ for $k \in \{1,\hdots,l\}$.
\item A collection of smooth Ricci flows $(M^{(k)},g^{(k)}(t))$, defined for $t \in [t_k^-,t_k^+)$ and going singular as $t \to t_k^+$ for $k \in \{0,1,\hdots,l-1\}$. 
\item Positive numbers $\varepsilon,r,\delta,h$, where $\delta \leq \varepsilon$ and $h \leq \delta r$. These are referred to as the surgery parameters.
\end{itemize}
For each $k \in \{0,1,\hdots,l-1\}$, we put $\Omega^{(k)} = \{x \in M^{(k)}: \limsup_{t \to t_k^+} \text{\rm scal}(x,t) < \infty\}$. We assume that the following conditions are satisfied: 
\begin{itemize} 
\item The manifold $(M^{(0)},g^{(0)}(0))$ is isometric to the given initial manifold $(M,g_0)$.
\item The manifold $(M^{(k)},g^{(k)}(t_k^-))$ is obtained from $(\Omega^{(k-1)},g^{(k-1)}(t_{k-1}^+))$ by performing surgery on finitely many necks. For each neck on which we perform surgery, we can find a point $(x_0,t_0)$ at the center of that neck such that $\text{\rm scal}(x_0,t_0)=h^{-2}$; moreover, the parabolic neighborhood $P(x_0,t_0,\delta^{-1} h,-\delta^{-1} h^2)$ is free of surgeries and is a $\delta$-neck.
\item After each surgery, we discard all double $4\varepsilon$-horns, all capped $4\varepsilon$-horns. Moreover, we remove all connected components which are diffeomorphic to $S^n/\Gamma$. 
\item Each flow $(M^{(k)},g^{(k)}(t))$ satisfies the Canonical Neighborhood Property with accuracy $4\varepsilon$ on all scales less than $r$. That is to say, if $(x_0,t_0)$ is an arbitrary point in space-time satisfying $\text{\rm scal}(x_0,t_0) \geq r^{-2}$, then there exists a neighborhood $B$ of $x_0$ with the property that $B_{g(t_0)}(x_0,(8C_1)^{-1} \, \text{\rm scal}(x_0,t_0)^{-\frac{1}{2}}) \subset B \subset B_{g(t_0)}(x_0,8C_1 \, \text{\rm scal}(x_0,t_0)^{-\frac{1}{2}})$ and $(8C_2)^{-1} \, \text{\rm scal}(x_0,t_0) \leq \text{\rm scal}(x,t_0) \leq 8C_2 \, \text{\rm scal}(x_0,t_0)$ for all $x \in B$. Moreover, $B$ is either a strong $4\varepsilon$-neck (in the sense of \cite{Perelman2}) with center at $x_0$ or a $4\varepsilon$-cap (in the sense of Definition \ref{definition.of.cap}). 
\item If $(x_0,t_0)$ is an arbitrary point in space-time satisfying $\text{\rm scal}(x_0,t_0) \geq r^{-2}$, then $|D \text{\rm scal}| \leq 4\eta \, \text{\rm scal}^{\frac{3}{2}}$ and $|\frac{\partial}{\partial t} \text{\rm scal}| \leq 4\eta \, \text{\rm scal}^2$ at $(x_0,t_0)$. 
\end{itemize}
\end{definition}

Note that the manifold $M^{(k)}$ may have multiple connected components. In the following, we will write the surgically modified solution simply as $g(t)$. However, it is important to remember that the underlying manifold changes across surgery times.

In the first step, we prove an upper bound for the length of the time interval on which the solution is defined.

\begin{proposition} 
\label{finite.time.extinction}
Suppose that we have a Ricci flow with surgery starting from $(M,g_0)$ which is defined on $[0,T)$. Then $T \leq \frac{n}{2 \, \inf_{x \in M} \text{\rm scal}(x,0)}$.
\end{proposition}

\textbf{Proof.} 
By the maximum principle, the function 
\[t \mapsto \frac{n}{2 \, \inf_{x \in M} \text{\rm scal}(x,t)} + t\] 
is monotone decreasing under smooth Ricci flow. By Proposition \ref{curvature.under.surgery}, this function is monotone decreasing across surgery times. From this, the assertion follows. \\

\begin{proposition}
\label{curvature.pinching.for.surgically.modified.flows}
Let $f$ and $\theta$ be as above. Moreover, let $g(t)$ be a Ricci flow with surgery starting from $(M,g_0)$. Then $(M,g(t))$ has $(f,\theta)$-pinched curvature.
\end{proposition}

\textbf{Proof.} 
By Theorem \ref{main.theorem} and Hamilton's PDE-ODE principle (cf. \cite{Chow-Lu}, Theorem 3, or \cite{Chow-et-al}, Theorem 10.16), the property that the curvature tensor of $g(t)$ lies in $\mathcal{F}_t$ is preserved by the Ricci flow. By Proposition \ref{curvature.under.surgery}, the property that the curvature tensor lies in $\mathcal{F}_t$ is preserved under surgery. Therefore, the property that the curvature tensor of $g(t)$ lies in $\mathcal{F}_t$ is preserved under Ricci flow with surgery. \\

\begin{proposition}
\label{shi.estimate.in.the.presence.of.surgeries}
Let $g(t)$ be a Ricci flow with surgery with surgery parameters $\varepsilon,r,\delta,h$. Let $(x_0,t_0)$ be a point in space time and let $r_0$ be a positive real number such that $t_0 \geq r_0^2$ and $\text{\rm scal}(x,t) \leq r_0^{-2}$ for all points $(x,t) \in P(x_0,t_0,r_0,-r_0^2)$. Then $|D^m R| \leq C(n,m) \, r_0^{-m-2}$ at the point $(x_0,t_0)$.
\end{proposition}

\textbf{Proof.} 
If the parabolic neighborhood $P(x_0,t_0,\frac{r_0}{2},-\frac{r_0^2}{4})$ is free of surgeries, this follows from the classical Shi estimate. Suppose next that the parabolic neighborhood $P(x_0,t_0,\frac{r_0}{2},-\frac{r_0^2}{4})$ does contain surgeries. At each point modified by surgery, the scalar curvature is at least $\frac{1}{4} \, h^{-2}$. Consequently, $\frac{1}{4} \, h^{-2} \leq r_0^{-2}$. The classical Shi estimate implies $|D^m R| \leq C(n,m) \, h^{-m-2} \leq C(n,m) \, r_0^{-m-2}$ on each strong neck on which we perform surgery. Moreover, $|D^m R| \leq C(n,m) \, h^{-m-2} \leq C(n,m) \, r_0^{-m-2}$ at each point modified by surgery. The assertion now follows from Theorem 3.29 in \cite{Morgan-Tian}. \\

\begin{proposition}[cf. G.~Perelman \cite{Perelman2}, Lemma 4.5]
\label{evolution.of.surgical.cap}
Fix $\varepsilon>0$ small, $\alpha \in [0,\frac{n-1}{2})$, and $A>1$. Then there exists $\bar{\delta}>0$ (depending on $\alpha$ and $A$) with the following property. Suppose that we have a Ricci flow with surgery with parameters $\varepsilon,r,\delta,h$, where $\delta \leq \bar{\delta}$. Suppose that $T_0 \in [0,T)$ is a surgery time, and let $x_0$ be a point which lies on a gluing cap at time $T_0$. Let $T_1 = \min \{T,T_0+\alpha h^2\}$. Then one of the following statements holds:
\begin{itemize}
\item[(i)] The flow is defined on $P(x_0,T_0,Ah,T_1-T_0)$. Moreover, after dilating the flow by $h^{-2}$ and shifting time $T_0$ to $0$, the parabolic neighborhood $P(x_0,T_0,Ah,T_1-T_0)$ is $A^{-1}$-close to the corresponding subset of the standard solution. 
\item[(ii)] There exists a surgery time $t^+ \in (T_0,T_1)$ such that the flow is defined on $P(x_0,T_0,Ah,t^+-T_0)$. Moreover, the parabolic neighborhood $P(x_0,T_0,Ah,t^+-T_0)$ is, after dilation by the factor $h^{-2}$, $A^{-1}$-close to the corresponding subset of the standard solution. Finally, for each point $x \in B_{g(T_0)}(x_0,Ah)$, the flow exists exactly until time $t^+$.
\end{itemize}
\end{proposition}

The proof is the same as the proof of Lemma 4.5 in Perelman's paper \cite{Perelman2}. We omit the details. \\

As in Perelman's work \cite{Perelman2}, it is crucial to establish a noncollapsing estimate in the presence of surgeries. 

\begin{definition}
\label{definition.of.noncollapsing.for.ricci.flow.with.surgery}
Suppose we are given a Ricci flow with surgery. We say that the flow is $\kappa$-noncollapsed on scales less than $\rho$ if the following holds. If $(x_0,t_0)$ is a point in space-time and $r_0$ is a positive number such that $r_0 \leq \rho$ and $\text{\rm scal}(x,t) \leq r_0^{-2}$ for all points $(x,t) \in P(x_0,t_0,r_0,-r_0^2)$ for which the flow is defined, then $\text{\rm vol}_{g(t_0)}(B_{g(t_0)}(x_0,r_0)) \geq \kappa \, r_0^n$. 
\end{definition}

As in Perelman's work \cite{Perelman2}, the noncollapsing estimate for Ricci flow with surgery will follow from the monotonicity formula for the reduced volume. 

\begin{definition}
Suppose we are given a Ricci flow with surgery. A curve in space-time is said to be admissible if it stays in the region unaffected by surgery. A curve in space-time is called barely admissible if it is on the boundary of the set of admissible curves.
\end{definition}

\begin{lemma}[cf. G.~Perelman \cite{Perelman2}, Lemma 5.3]
\label{lower.bound.for.L.length}
Fix $\varepsilon,r,L$. Then we can find a real number $\bar{\delta}>0$ (depending on $r$ and $L$) with the following property. Suppose that we have a Ricci flow with surgery with parameters $\varepsilon,r,\delta,h$, where $\delta \leq \bar{\delta}$. Let $(x_0,t_0)$ be a point in space-time such that $\text{\rm scal}(x_0,t_0) \leq r^{-2}$, and let $T_0<t_0$ be a surgery time. Finally, let $\gamma$ be a barely admissible curve (parametrized by the interval $[T_0,t_0]$) such that $\gamma(T_0)$ lies on the boundary of a surgical cap at time $T_0$, and $\gamma(t_0)=x_0$. Then  
\[\int_{T_0}^{t_0} \sqrt{t_0-t} \, (\text{\rm scal}(\gamma(t),t) + |\gamma'(t)|_{g(t)}^2) \, dt \geq L.\] 
\end{lemma}

\textbf{Proof.} 
By Definition \ref{definition.ricci.flow.with.surgery}, we have $|D \text{\rm scal}| \leq 4\eta \, \text{\rm scal}^{\frac{3}{2}}$ and $|\frac{\partial}{\partial t} \text{\rm scal}| \leq 4\eta \, \text{\rm scal}^2$ whenever $\text{\rm scal} \geq r^{-2}$. Since $\text{\rm scal}(x_0,t_0) \leq r^{-2}$, it follows that $\text{\rm scal} \leq 4r^{-2}$ in $P(x_0,t_0,\frac{r}{100\eta},-\frac{r^2}{100\eta})$. Let $\gamma$ be a barely admissible curve in space-time satisfying the assumptions of Lemma \ref{lower.bound.for.L.length}, and suppose that 
\[\int_{T_0}^{t_0} \sqrt{t_0-t} \, (\text{\rm scal}(\gamma(t),t) + |\gamma'(t)|_{g(t)}^2) \, dt < L.\] 
Using H\"older's inequality and the positivity of the scalar curvature, we obtain $\int_{t_0-\tau}^{t_0} |\gamma'(t)|_{g(t)} \, dt < (2L)^{\frac{1}{2}} \, \tau^{\frac{1}{4}}$ for $\tau>0$. Hence, we can find a real number $\tau \in (0,\frac{r^2}{100\eta})$, depending only on $r$ and $L$, such that $\gamma|_{[t_0-\tau,t_0]}$ is contained in the parabolic neighborhood $P(x_0,t_0,\frac{r}{100\eta},-\frac{r^2}{100\eta})$. This implies 
\[\text{\rm scal}(\gamma(t),t) \leq 4r^{-2}\] 
for all $t \in [t_0-\tau,t_0]$.

Having chosen $\tau$, we define real numbers $\alpha \in [0,\frac{n-1}{2})$ and $l>0$ by the relations 
\[\frac{(n-1)\sqrt{\tau}}{4K_{\text{\rm stnd}}} \, \Big | \log \Big ( 1-\frac{2\alpha}{n-1} \Big ) \Big | = L\] 
and 
\[\frac{l}{2} \, \sqrt{\tau} = L.\] 
Having fixed $\alpha$ and $l$, we choose a large constant $A$ so that the conclusion of Lemma \ref{no.escape} holds. Having chosen $\alpha$ and $A$, we choose $\bar{\delta}$ so that the conclusion of Proposition \ref{evolution.of.surgical.cap} holds. Moreover, by choosing $\bar{\delta}$ small enough, we can arrange that $K_{\text{\rm stnd}} \, \bar{\delta}^2 \leq \frac{1}{16}$.

In the following, we assume that $\delta \leq \bar{\delta}$. Let $T_1 \in [T_0,T_0+\alpha h^2]$ denote the largest number with the property that $\gamma|_{[T_0,T_1]}$ is contained in the parabolic neighborhood $P(\gamma(T_0),T_0,Ah,\alpha h^2)$. By Proposition \ref{evolution.of.surgical.cap}, the parabolic neighborhood $P(\gamma(T_0),T_0,Ah,T_1-T_0)$ is close to the corresponding subset of the standard solution. Since $h \leq \delta r$, it follows that 
\[\text{\rm scal}(\gamma(t),t) \geq \frac{1}{2K_{\text{\rm stnd}} \, (h^2-\frac{2(t-T_0)}{n-1})} \geq \frac{1}{2K_{\text{\rm stnd}} \, \delta^2 r^2} \geq 8r^{-2}\] 
for all $t \in [T_0,T_1]$. Since $\text{\rm scal}(\gamma(t),t) \leq 4r^{-2}$ for all $t \in [t_0-\tau,t_0]$, the intervals $[T_0,T_1]$ and $[t_0-\tau,t_0]$ are disjoint. In other words, $T_1 \leq t_0-\tau$. We distinguish two cases: 

\textit{Case 1:} Suppose that $T_1 < T_0 + \alpha h^2$. In this case, the curve $\gamma|_{[T_0,T_1]}$ exits the parabolic neighborhood $P(\gamma(T_0),T_0,Ah,\alpha h^2)$ at time $T_1$. Since the parabolic neighborhood $P(\gamma(T_0),T_0,Ah,T_1-T_0)$ is close to the corresponding subset of the standard solution, Lemma \ref{no.escape} implies that $\int_{T_0}^{T_1} |\gamma'(t)|_{g(t)}^2 \, dt \geq \frac{l}{2}$. (Here, we have used the fact that $\int |\gamma'(t)|_{g(t)}^2 \, dt$ is invariant under scaling.) Consequently, 
\begin{align*} 
L &> \int_{T_0}^{T_1} \sqrt{t_0-t} \, (\text{\rm scal}(\gamma(t),t) + |\gamma'(t)|_{g(t)}^2) \, dt \\ 
&\geq \sqrt{\tau} \int_{T_0}^{T_1} |\gamma'(t)|_{g(t)}^2 \, dt \\ 
&\geq \frac{l}{2} \, \sqrt{\tau}, 
\end{align*} 
which contradicts our choice of $l$.

\textit{Case 2:} Suppose that $T_1 = T_0 + \alpha h^2$. In this case, 
\begin{align*} 
L &> \int_{T_0}^{T_1} \sqrt{t_0-t} \, (\text{\rm scal}(\gamma(t),t) + |\gamma'(t)|_{g(t)}^2) \, dt \\ 
&\geq \sqrt{\tau} \int_{T_0}^{T_1} \text{\rm scal}(\gamma(t),t) \, dt \\ 
&\geq \sqrt{\tau} \int_{T_0}^{T_1} \frac{1}{2K_{\text{\rm stnd}} \, (h^2-\frac{2(t-T_0)}{n-1})} \, dt \\ 
&= \frac{(n-1)\sqrt{\tau}}{4K_{\text{\rm stnd}}} \, \Big | \log \Big ( 1-\frac{2\alpha}{n-1} \Big ) \Big |,
\end{align*}  
which contradicts our choice of $\alpha$. \\

\begin{proposition}[cf. G.~Perelman \cite{Perelman2}, Lemma 5.2]
\label{noncollapsing.for.ricci.flows.with.surgery}
Fix a small number $\varepsilon>0$. Then we can find a positive number $\kappa$ and a positive function $\tilde{\delta}(\cdot)$ with the following property. Suppose that we have a Ricci flow with surgery with parameters $\varepsilon,r,\delta,h$, where $\delta \leq \tilde{\delta}(r)$. Then the flow is $\kappa$-noncollapsed on all scales less than $\varepsilon$.
\end{proposition}

Note that the constant $\kappa$ in the noncollapsing estimate may depend on the initial data, but it is independent of the surgery parameters $\varepsilon,r,\delta,h$. \\

\textbf{Proof.}
Consider a point $(x_0,t_0)$ in space-time and a positive number $r_0 \leq \varepsilon$ such that $\text{\rm scal}(x,t) \leq r_0^{-2}$ for all points $(x,t) \in P(x_0,t_0,r_0,-r_0^2)$ for which the flow is defined. We need to show that $\text{\rm vol}_{g(t_0)}(B_{g(t_0)}(x_0,r_0)) \geq \kappa \, r_0^n$ for some uniform constant $\kappa>0$. We distinguish three cases: 

\textit{Case 1:} Suppose first that $\text{\rm scal}(x_0,t_0) \geq r^{-2}$. In this case, the Canonical Neighborhood Assumption implies that $\text{\rm vol}_{g(t_0)}(B_{g(t_0)}(x_0,r_0)) \geq \kappa \, r_0^n$ for some uniform constant $\kappa>0$. 

\textit{Case 2:} Suppose next that the parabolic neighborhood $P(x_0,t_0,\frac{r_0}{2},-\frac{r_0^2}{4})$ contains points modified by surgery. Let $(x,t)$ be a point in $P(x_0,t_0,\frac{r_0}{2},-\frac{r_0^2}{4})$ which lies on a surgical cap. Clearly, $\frac{1}{4} \, h^{-2} \leq \text{\rm scal}(x,t) \leq r_0^{-2}$, hence $r_0 \leq 2h$. This implies $\text{\rm vol}_{g(t)}(B_{g(t)}(x,\frac{r_0}{100})) \geq \kappa \, r_0^n$ for some uniform constant $\kappa>0$. Since $B_{g(t)}(x,\frac{r_0}{100}) \subset B_{g(t_0)}(x,\frac{r_0}{4})$, we deduce that $\text{\rm vol}_{g(t)}(B_{g(t_0)}(x_0,r_0)) \geq \text{\rm vol}_{g(t)}(B_{g(t_0)}(x,\frac{r_0}{4})) \geq \kappa \, r_0^n$.

\textit{Case 3:} Suppose finally that $\text{\rm scal}(x_0,t_0) \leq r^{-2}$ and the parabolic neighborhood $P(x_0,t_0,\frac{r_0}{2},-\frac{r_0^2}{4})$ is free of surgeries. Note that $t_0$ is bounded from above by Proposition \ref{finite.time.extinction}. By Lemma \ref{lower.bound.for.L.length}, we can find a positive function $\tilde{\delta}(\cdot)$ such that the following holds: Suppose that the surgery parameters satisfy $\delta \leq \tilde{\delta}(r)$, and suppose further that $T_0<t_0$ is a surgery time and $\gamma$ is a barely admissible curve (parametrized by the interval $[T_0,t_0]$) such that $\gamma(T_0)$ lies on the boundary of a surgical cap at time $T_0$ and $\gamma(t_0)=x_0$. Then 
\[\int_{T_0}^{t_0} \sqrt{t_0-t} \, (\text{\rm scal}(\gamma(t),t) + |\gamma'(t)|_{g(t)}^2) \, dt \geq 8n \sqrt{t_0}.\] 
Thus, if $\delta \leq \tilde{\delta}(r)$, then every barely admissible curve has reduced length greater than $2n$.

In the following, we assume that $\delta \leq \tilde{\delta}(r)$. For $t<t_0$, we denote by $\ell(x,t)$ the reduced distance from $(x_0,t_0)$, i.e. the infimum of the reduced length over all admissible curves joining $(x,t)$ and $(x_0,t_0)$. We claim that $\inf_x \ell(x,t) \leq \frac{n}{2}$ for all $t < t_0$. This is clearly true if $t$ is sufficiently close to $t_0$. Now, if $\ell(x,t) < 2n$ for some point $(x,t)$ in space-time, then the reduced length is attained by a strictly admissible curve. Hence, we may apply results of Perelman (cf. \cite{Perelman1}, Section 7) to conclude that 
\[\frac{\partial}{\partial t} \ell \geq \Delta \ell + \frac{1}{t_0-t} \, \Big ( \ell-\frac{n}{2} \Big )\] 
whenever $\ell < 2n$. Using the maximum principle, we deduce that $\inf_x \ell(x,t) \leq \frac{n}{2}$ for all $t < t_0$. 

In particular, there exists a point $y \in M$ such that $\ell(y,\varepsilon) \leq \frac{n}{2}$. Hence, we can find a radius $\rho>0$ such that $\sup_{x \in B_{g(0)}(y,\rho)} \ell(x,0) \leq n$. Note that $\rho$ depends only on $\varepsilon$ and the initial data $(M,g_0)$, but not on the surgery parameters. Hence, for each point $x \in B_{g(0)}(y,\rho)$, the reduced distance is attained by a strictly admissible curve, and this curve must be an $\mathcal{L}$-geodesic. 

Given a tangent vector $v$ at $(x_0,t_0)$, we denote by $\gamma_v(t) = \mathcal{L}_{t,t_0} \exp_{x_0}(v)$ the $\mathcal{L}$-geodesic satisfying $\lim_{t \to t_0} \sqrt{t_0-t} \, \gamma_v'(t) = v$. Note that, due to the presence of surgeries, $\gamma_v(t)$ may not be defined on the entire interval $[0,t_0)$. Let $\mathcal{V}$ denote the set of all tangent vectors $v$ at $(x_0,t_0)$ with the property that $\gamma_v$ is defined on $[0,t_0)$; $\gamma_v$ has minimal $\mathcal{L}$-length; and $\gamma_v(0) \in B_{g(0)}(y,\rho)$. In view of the discussion above, the map $\mathcal{L}_{0,t_0} \exp_{x_0}: \mathcal{V} \to B_{g(0)}(y,\rho)$ is onto. For each $t \in [0,t_0)$, we define  
\[V(t) = \int_{\mathcal{V}} (t_0-t)^{-\frac{n}{2}} \, e^{-\ell(\gamma_v(t),t)} \, J_v(t),\] 
where $J_v(t) = \det (D \mathcal{L}_{t,t_0} \exp_{x_0})_v$ denotes the Jacobian determinant of the $\mathcal{L}$-exponential map, and the integration is with respect to the standard Lebesgue measure on the tangent space $(T_{x_0} M,g(t_0))$. For each tangent vector $v \in \mathcal{V}$, Perelman's Jacobian comparison theorem (cf. \cite{Perelman1}, Section 7) implies that the function $t \mapsto (t_0-t)^{-\frac{n}{2}} \, e^{-\ell(\gamma_v(t),t)} \, J_v(t)$ is monotone increasing. Moreover, $\lim_{t \to t_0} (t_0-t)^{-\frac{n}{2}} \, e^{-\ell(\gamma_v(t),t)} \, J_v(t) = 2^n \, e^{-|v|^2}$ for each $v \in \mathcal{V}$. The monotonicity property for the Jacobian determinant implies that the function $t \mapsto V(t)$ is monotone increasing. 

We first estimate the reduced volume from below in terms of the initial data. Since $\ell(x,0) \leq n$ for all points $x \in B_{g(0)}(y,\rho)$, we obtain a uniform lower bound for $V(0)$:
\begin{align*} 
V(0) 
&= \int_{\mathcal{V}} t_0^{-\frac{n}{2}} \, e^{-\ell(\gamma_v(0),0)} \, J_v(0) \\ 
&\geq \int_{B_{g(0)}(y,\rho)} t_0^{-\frac{n}{2}} \, e^{-\ell(x,0)} \, d\text{\rm vol}_{g(0)}(x) \\ 
&\geq t_0^{-\frac{n}{2}} \, e^{-n} \, \text{\rm vol}_{g(0)}(B_{g(0)}(y,\rho)). 
\end{align*} 
We next estimate the reduced volume from above. By assumption, the parabolic neighborhood $P(x_0,t_0,\frac{r_0}{2},-\frac{r_0^2}{4})$ is free of surgeries, and we have $\text{\rm scal} \leq r_0^{-2}$ in $P(x_0,t_0,\frac{r_0}{2},-\frac{r_0^2}{4})$. Using Shi's interior derivative estimates, we conclude that the covariant derivatives of the Riemann curvature tensor are bounded by $C(n) \, r_0^{-3}$ and the second covariant derivatives of the Riemann curvature tensor are bounded by $C(n) \, r_0^{-4}$ in $P(x_0,t_0,\frac{r_0}{4},-\frac{r_0^2}{16})$. Using the $\mathcal{L}$-geodesic equation, we conclude that there exists a small positive constant $\mu(n)$ (depending only on $n$) with the following property: if $\bar{t} \in [t_0-\mu(n) r_0^2,t_0)$ and $|v| \leq \frac{r_0}{32\sqrt{t_0-\bar{t}}}$, then $\sqrt{t_0-t} \, |\gamma_v'(t)|_{g(t)} \leq \frac{r_0}{16\sqrt{t_0-\bar{t}}}$ and $\gamma_v(t) \in B_{g(t_0)}(x_0,\frac{r_0 \sqrt{t_0-t}}{4\sqrt{t_0-\bar{t}}}) \subset B_{g(t_0)}(x_0,\frac{r_0}{4})$ for all $t \in [\bar{t},t_0)$. This implies 
\begin{align*} 
V(0) \leq V(\bar{t}) 
&\leq \int_{\{v \in \mathcal{V}: |v| \leq \frac{r_0}{32\sqrt{t_0-\bar{t}}}\}} (t_0-\bar{t})^{-\frac{n}{2}} \, e^{-\ell(\gamma_v(\bar{t}),\bar{t})} \, J_v(\bar{t}) \\ 
&+ \int_{\{v \in \mathcal{V}: |v| \geq \frac{r_0}{32\sqrt{t_0-\bar{t}}}\}} (t_0-\bar{t})^{-\frac{n}{2}} \, e^{-\ell(\gamma_v(\bar{t}),\bar{t})} \, J_v(\bar{t}) \\ 
&\leq \int_{\{v \in \mathcal{V}: |v| \leq \frac{r_0}{32\sqrt{t_0-\bar{t}}}\}} (t_0-\bar{t})^{-\frac{n}{2}} \, J_v(\bar{t}) \\ 
&+ \int_{\{|v| \geq \frac{r_0}{32\sqrt{t_0-\bar{t}}}\}} 2^n \, e^{-|v|^2} \\ 
&\leq (t_0-\bar{t})^{-\frac{n}{2}} \, \text{\rm vol}_{g(\bar{t})}(B_{g(t_0)}(x_0,\frac{r_0}{4})) \\ 
&+ \int_{\{|v| \geq \frac{r_0}{32\sqrt{t_0-\bar{t}}}\}} 2^n \, e^{-|v|^2} 
\end{align*} 
for all $\bar{t} \in [t_0-\mu(n) r_0^2,t_0)$. Putting these facts together gives 
\begin{align*} 
&(t_0-\bar{t})^{-\frac{n}{2}} \, \text{\rm vol}_{g(\bar{t})}(B_{g(t_0)}(x_0,\frac{r_0}{4})) \\ 
&\geq t_0^{-\frac{n}{2}} \, e^{-n} \, \text{\rm vol}_{g(0)}(B_{g(0)}(y,\rho)) - \int_{\{|v| \geq \frac{r_0}{32\sqrt{t_0-\bar{t}}}\}} 2^n \, e^{-|v|^2} 
\end{align*}
for all $\bar{t} \in [t_0-\mu(n) r_0^2,t_0)$. Finally, we choose $\bar{t} \in [t_0-\mu(n) r_0^2,t_0)$ so that $t_0-\bar{t}$ is a small, but fixed, multiple of $r_0^2$, and the quantity 
\[t_0^{-\frac{n}{2}} \, e^{-n} \, \text{\rm vol}_{g(0)}(B_{g(0)}(y,\rho)) - \int_{\{|v| \geq \frac{r_0}{32\sqrt{t_0-\bar{t}}}\}} 2^n \, e^{-|v|^2}\] 
is bounded from below by a positive constant. This gives a lower bound for $r_0^{-n} \, \text{\rm vol}_{g(\bar{t})}(B_{g(t_0)}(x_0,\frac{r_0}{4}))$, as desired. \\

We now state the main result of this section. This result guarantees that, for a suitable choice of $\varepsilon$, $\hat{r}$, $\hat{\delta}$, every Ricci flow with surgery with parameters $\varepsilon,\hat{r},\hat{\delta},h$ will satisfy the Canonical Neighborhood Property with accuracy $2\varepsilon$ on all scales less than $2\hat{r}$. 

\begin{theorem}[cf. G.~Perelman \cite{Perelman2}, Section 5]
\label{canonical.neighborhood.property.improved.accuracy}
Fix a small number $\varepsilon>0$. Then we can find positive numbers $\hat{r}$ and $\hat{\delta}$ with the following property. Suppose that we have a Ricci flow with surgery with parameters $\varepsilon,\hat{r},\hat{\delta},h$ which is defined on some interval $[0,T)$. Moreover, suppose that $(x_0,t_0)$ is an arbitrary point in space-time satisfying $\text{\rm scal}(x_0,t_0) \geq (2\hat{r})^{-2}$. Then there exists a neighborhood $B$ of $x_0$ such that $B_{g(t_0)}(x_0,(2C_1)^{-1} \, \text{\rm scal}(x_0,t_0)^{-\frac{1}{2}}) \subset B \subset B_{g(t_0)}(x_0,2C_1 \, \text{\rm scal}(x_0,t_0)^{-\frac{1}{2}})$ and $(2C_2)^{-1} \, \text{\rm scal}(x_0,t_0) \leq \text{\rm scal}(x,t_0) \leq 2C_2 \, \text{\rm scal}(x_0,t_0)$ for all $x \in B$. Moreover, $B$ is either a strong $2\varepsilon$-neck (in the sense of \cite{Perelman2}) with center at $x_0$ or a $2\varepsilon$-cap. Finally, we have $|D \text{\rm scal}| \leq 2\eta \, \text{\rm scal}^{\frac{3}{2}}$ and $|\frac{\partial}{\partial t} \text{\rm scal}| \leq 2\eta \, \text{\rm scal}^2$ at the point $(x_0,t_0)$.
\end{theorem}

\textbf{Proof.} 
We argue by contradiction. Suppose that the assertion is false. Then we can find a sequence of Ricci flows with surgery $\mathcal{M}^{(j)}$ and a sequence of points $(x_j,t_j)$ in space-time with the following properties: 
\begin{itemize} 
\item[(i)] The flow $\mathcal{M}^{(j)}$ is defined on the time interval $[0,T_j)$ and has surgery parameters $\varepsilon,\hat{r}_j,h_j,\hat{\delta}_j$, where $\hat{r}_j \leq \frac{1}{j}$ and $\hat{\delta}_j \leq \min \{\tilde{\delta}(\hat{r}_j),\frac{1}{j}\}$. Here, $\tilde{\delta}(\cdot)$ is the function introduced in Proposition \ref{noncollapsing.for.ricci.flows.with.surgery}.
\item[(ii)] $Q_j := \text{\rm scal}(x_j,t_j) \geq (2\hat{r}_j)^{-2}$. 
\item[(iii)] The point $(x_j,t_j)$ does not satisfy the conclusion of Theorem \ref{canonical.neighborhood.property.improved.accuracy}. 
\end{itemize}
The condition (iii) means that at least one of the following statements is true: 
\begin{itemize}
\item There does not exist a neighborhood $B$ of $x_j$ with the property that $B_{g(t_j)}(x_j,(2C_1)^{-1} \, \text{\rm scal}(x_j,t_j)^{-\frac{1}{2}}) \subset B \subset B_{g(t_j)}(x_j,2C_1 \, \text{\rm scal}(x_j,t_j)^{-\frac{1}{2}})$ and $(2C_2)^{-1} \, \text{\rm scal}(x_j,t_j) \leq \text{\rm scal}(x,t_j) \leq 2C_2 \, \text{\rm scal}(x_j,t_j)$ for all $x \in B$, and such that $B$ is either a strong $2\varepsilon$-neck (in the sense of \cite{Perelman2}) with center at $x_j$ or a $2\varepsilon$-cap. 
\item $|D \text{\rm scal}| > 2\eta \, \text{\rm scal}^{\frac{3}{2}}$ at $(x_j,t_j)$. 
\item $|\frac{\partial}{\partial t} \text{\rm scal}| > 2\eta \, \text{\rm scal}^2$ at $(x_j,t_j)$. 
\end{itemize}
We will proceed in several steps: 

\textit{Step 1:} By Definition \ref{definition.ricci.flow.with.surgery}, we have $|D \text{\rm scal}| \leq 4\eta \, \text{\rm scal}^{\frac{3}{2}}$ and $|\frac{\partial}{\partial t} \text{\rm scal}| \leq 4\eta \, \text{\rm scal}^2$ for each point $(x,t)$ in space-time satisfying $\text{\rm scal}(x,t) \geq 4Q_j$. Moreover, by Proposition \ref{noncollapsing.for.ricci.flows.with.surgery}, the flow $\mathcal{M}^{(j)}$ is $\kappa$-noncollapsed on scales less than $\varepsilon$ for some uniform constant $\kappa$ which may depend on the initial data, but is independent of $j$. 

\textit{Step 2:} Suppose that $(x_0,t_0)$ is a point in space-time satisfying $\text{\rm scal}(x_0,t_0)+Q_j \leq r_0^{-2}$. The pointwise curvature derivative estimate implies that $\text{\rm scal} \leq 8r_0^{-2}$ in the parabolic neighborhood $P(x_0,t_0,\frac{r_0}{100\eta},-\frac{r_0^2}{100\eta})$. Using Proposition \ref{shi.estimate.in.the.presence.of.surgeries}, we conclude that $|D^m R| \leq C(n,m,\eta) \, r_0^{-m-2}$ at the point $(x_0,t_0)$. Moreover, Proposition \ref{noncollapsing.for.ricci.flows.with.surgery} implies $\text{\rm vol}_{g(t_0)}(B_{g(t_0)}(x_0,r_0)) \geq \kappa r_0^n$ for some uniform constant $\kappa$ which is independent of $j$. 

\textit{Step 3:} We next prove a longrange curvature estimate. Given any $\rho>0$, we put 
\[\mathbb{M}(\rho) = \limsup_{j \to \infty} \sup_{x \in B_{g(t_j)}(x_j,\rho Q_j^{-\frac{1}{2}})} Q_j^{-1} \, \text{\rm scal}(x,t_j).\] 
The pointwise curvature derivative estimate implies that $\mathbb{M}(\rho) \leq 16$ for $0 < \rho < \frac{1}{100\eta}$. 

We claim that $\mathbb{M}(\rho) < \infty$ for all $\rho>0$. Suppose this is false. Let 
\[\rho^* = \sup \{\rho \geq 0: \mathbb{M}(\rho) < \infty\} < \infty.\] 
By definition of $\rho^*$, we have an upper bound for the curvature in the geodesic ball $B_{g(t_j)}(x_j,\rho Q_j^{-\frac{1}{2}})$ for each $\rho<\rho^*$. Using the results in Step 2, we obtain bounds for all the covariant derivatives of the curvature tensor in the geodesic ball $B_{g(t_j)}(x_j,\rho Q_j^{-\frac{1}{2}})$ for each $\rho<\rho^*$. Moreover, the noncollapsing estimate in Step 2 gives a lower bound for the volume. We rescale around $(x_j,t_j)$ by the factor $Q_j$ and pass to the limit as $j \to \infty$. In the limit, we obtain an incomplete manifold $(B^\infty,g^\infty)$ which is weakly PIC2 (cf. \cite{Morgan-Tian}, Theorem 5.6).

By definition of $\rho^*$, there exists a sequence of points $y_j$ such that $\rho_j := Q_j^{\frac{1}{2}} \, d_{g(t_j)}(x_j,y_j) \to \rho^*$ and $Q_j^{-1} \, \text{\rm scal}(y_j,t_j) \to \infty$. Let $\gamma_j: [0,\rho_j Q_j^{-\frac{1}{2}}] \to (M,g(t_j))$ be a unit-speed geodesic such that $\gamma_j(0)=x_j$ and $\gamma_j(\rho_j Q_j^{-\frac{1}{2}}) = y_j$, and let $\gamma_\infty: [0,\rho^*) \to (B^\infty,g^\infty)$ denote the limit of $\gamma_j$. Using the estimate $|D \text{\rm scal}| \leq 4\eta \, \text{\rm scal}^{\frac{3}{2}}$, we obtain 
\[\text{\rm scal}_{g^\infty}(\gamma_\infty(s)) = \lim_{j \to \infty} Q_j^{-1} \, \text{\rm scal}(\gamma_j(sQ_j^{-\frac{1}{2}}),t_j) \geq (2\eta(\rho^*-s))^{-2} \geq 100\] 
if $s \in [\rho^*-\frac{1}{100\eta},\rho^*)$.

Let us consider a real number $\bar{s} \in [\rho^*-\frac{1}{100\eta},\rho^*)$ such that $64C_1 \eta (\rho^*-\bar{s}) \leq \bar{s}$. We claim that $\gamma_j(\bar{s}Q_j^{-\frac{1}{2}})$ lies at the center of a strong $4\varepsilon$-neck if $j$ is sufficiently large (depending on $\bar{s}$). Indeed, if $j$ is sufficiently large, then the Canonical Neighborhood Assumption implies that the point $(\gamma_j(\bar{s}Q_j^{-\frac{1}{2}}),t_j)$ has a Canonical Neighborhood which is either a strong $4\varepsilon$-neck or a $4\varepsilon$-cap. Furthermore, the Canonical Neighborhood is contained in a geodesic ball around $\gamma_j(\bar{s}Q_j^{-\frac{1}{2}})$ of radius $8C_1 \, \text{\rm scal}(\gamma_j(\bar{s}Q_j^{-\frac{1}{2}}),t_j)^{-\frac{1}{2}}$, and the scalar curvature is at most $8C_2 \, \text{\rm scal}(\gamma_j(\bar{s}Q_j^{-\frac{1}{2}}),t_j)$ at each point in the Canonical Neighborhood. Since $\mathbb{M}(\bar{s}) < \infty$, we obtain $\lim_{j \to \infty} \text{\rm scal}(\gamma_j(\bar{s}Q_j^{-\frac{1}{2}}),t_j)^{-1} \, \text{\rm scal}(y_j,t_j) = \infty$; consequently, $\text{\rm scal}(y_j,t_j) \geq 16C_2 \, \text{\rm scal}(\gamma_j(\bar{s}Q_j^{-\frac{1}{2}}),t_j)$ if $j$ is sufficiently large. Hence, if $j$ is sufficiently large, then the Canonical Neighborhood does not contain the point $y_j$. We next observe that $32C_1 \, \text{\rm scal}_{g^\infty}(\gamma_\infty(\bar{s}))^{-\frac{1}{2}} \leq 64C_1 \eta (\rho^*-\bar{s}) \leq \bar{s}$. This implies $16C_1 \, \text{\rm scal}(\gamma_j(\bar{s}Q_j^{-\frac{1}{2}}),t_j)^{-\frac{1}{2}} \leq \bar{s} Q_j^{-\frac{1}{2}}$ if $j$ is sufficiently large. Hence, if $j$ is sufficiently large, then the Canonical Neighborhood does not contain the point $x_j$. If the Canonical Neighborhood of $(\gamma_j(\bar{s}Q_j^{-\frac{1}{2}}),t_j)$ is a $4\varepsilon$-cap, then the geodesic $\gamma_j$ must enter and exit this $4\varepsilon$-cap, but this is impossible since $\gamma_j$ minimizes length. To summarize, if $j$ is sufficiently large (depending on $\bar{s}$), then the point $(\gamma_j(\bar{s}Q_j^{-\frac{1}{2}}),t_j)$ has a Canonical Neighborhood which is a strong $4\varepsilon$-neck. In particular, if $j$ is sufficiently large (depending on $\bar{s}$), then we have $|D \text{\rm scal}| \leq C(n)\varepsilon \, \text{\rm scal}^{\frac{3}{2}}$ at the point $(\gamma_j(\bar{s}Q_j^{-\frac{1}{2}}),t_j)$.

Passing to the limit as $j \to \infty$, we conclude that $|D \text{\rm scal}_{g^\infty}| \leq C(n)\varepsilon \, \text{\rm scal}_{g^\infty}^{\frac{3}{2}}$ at the point $\gamma_\infty(\bar{s})$.  Integrating this estimate along $\gamma_\infty$ gives $\text{\rm scal}_{g^\infty}(\gamma_\infty(\bar{s})) \geq (C(n)\varepsilon (\rho^*-\bar{s}))^{-2}$. Moreover, since $(\gamma_j(\bar{s}Q_j^{-\frac{1}{2}}),t_j)$ lies at the center of a strong $4\varepsilon$-neck for $j$ sufficiently large, the point $\gamma_\infty(\bar{s})$ must lie on a strong $C(n)\varepsilon$-neck in $(B^\infty,g^\infty)$.

As in \cite{Perelman1}, Section 12.1, there is a sequence of rescalings which converges to a piece of a non-flat metric cone in the limit. Let us fix a point on this metric cone. In view of the preceding discussion, this point must lie on a strong $C(n)\varepsilon$-neck. This gives a locally defined solution to the Ricci flow which is weakly PIC2 and which, at the final time, is a piece of non-flat metric cone. This contradicts Proposition \ref{splitting.a}.

\textit{Step 4:} We now dilate the manifold $(M,g(t_j))$ around the point $x_j$ by the factor $Q_j$. By Step 3, we have uniform bounds for the curvature at bounded distance. Using the results in Step 2, we obtain bounds for all the covariant derivatives of the curvature tensor at bounded distance. Using these estimates together with the noncollapsing estimate in Step 2, we conclude that the rescaled manifolds converge in the Cheeger-Gromov sense to a complete limit manifold $(M^\infty,g^\infty)$. Since $(M,g(t_j))$ has $(f,\theta)$-pinched curvature, the curvature tensor of $(M^\infty,g^\infty)$ is weakly PIC2 and satisfies $R - \theta \, \text{\rm scal} \, \text{\rm id} \owedge \text{\rm id} \in \text{\rm PIC}$. Using the Canonical Neighborhood Assumption, we conclude that every point in $(M^\infty,g^\infty)$ with scalar curvature greater than $4$ has a neighborhood which is either a strong $8\varepsilon$-neck or a $8\varepsilon$-cap.

We claim that $(M^\infty,g^\infty)$ has bounded curvature. Indeed, if there is a sequence of points in $(M^\infty,g^\infty)$ with curvature going to infinity, then $(M^\infty,g^\infty)$ contains a sequence of necks with radii converging to $0$, contradicting Proposition \ref{no.small.necks}. This shows that $(M^\infty,g^\infty)$ has bounded curvature.

\textit{Step 5:} We now extend the limit $(M^\infty,g^\infty)$ backwards in time. By Step 4, the scalar curvature of $(M^\infty,g^\infty)$ is bounded from above by a constant $\Lambda>4$. We claim that, given any $A>1$, the parabolic neighborhood $P(x_j,t_j,A Q_j^{-\frac{1}{2}},-\frac{1}{100\eta \Lambda} \, Q_j^{-1})$ is free of surgeries if $j$ is sufficiently large. To prove this, fix $A>1$ and suppose that $P(x_j,t_j,A Q_j^{-\frac{1}{2}},-\frac{1}{100\eta \Lambda} \, Q_j^{-1})$ contains points modified by surgery. Let $s_j \in [0,\frac{1}{100\eta \Lambda}]$ be the largest number such that $P(x_j,t_j,A Q_j^{-\frac{1}{2}},-s_j Q_j^{-1})$ is free of surgeries. The pointwise curvature derivative estimate gives 
\[\sup_{P(x_j,t_j,A Q_j^{-\frac{1}{2}},-s_j Q_j^{-1})} \text{\rm scal} \leq 2\Lambda \, Q_j\] 
if $j$ is sufficiently large. Since the scalar curvature is greater than $\frac{1}{2} \, h_j^{-2}$ at each point modified by surgery, we deduce that $\frac{1}{2} \, h_j^{-2} \leq 2\Lambda \, Q_j$ if $j$ is sufficiently large. In particular, $s_j Q_j^{-1} \leq \frac{1}{10\eta} \, h_j^2$ if $j$ is sufficiently large. Since $\hat{\delta}_j \to 0$, Proposition \ref{evolution.of.surgical.cap} implies that the parabolic neighborhood $P(x_j,t_j,A Q_j^{-\frac{1}{2}},-s_j Q_j^{-1})$ is, after dilating by the factor $h_j$, arbitrarily close to a piece of the standard solution when $j$ is sufficiently large. Using Corollary \ref{canonical.neighborhood.theorem.for.standard.solution.2}, we conclude that $(x_j,t_j)$ lies on a $2\varepsilon$-neck or a $2\varepsilon$-cap when $j$ is sufficiently large. If $(x_j,t_j)$ lies on an $2\varepsilon$-neck, then this neck is actually a strong $2\varepsilon$-neck, since we are assuming that each $\hat{\delta}_j$-neck on which we perform surgery has a large backward parabolic neighborhood that is free of surgeries (cf. Definition \ref{definition.ricci.flow.with.surgery}). Moreover, Corollary \ref{canonical.neighborhood.theorem.for.standard.solution.2} implies that $|D \text{\rm scal}| \leq 2\eta \, \text{\rm scal}^{\frac{3}{2}}$ and $|\frac{\partial}{\partial t} \text{\rm scal}| \leq 2\eta \, \text{\rm scal}^2$ at the point $(x_j,t_j)$. Therefore, the point $(x_j,t_j)$ satisfies the conclusion of Theorem \ref{canonical.neighborhood.property.improved.accuracy}, and this contradicts statement (iii). Thus, given any $A>1$, the parabolic neighborhood $P(x_j,t_j,A Q_j^{-\frac{1}{2}},-\frac{1}{100\eta \Lambda} \, Q_j^{-1})$ is free of surgeries if $j$ is sufficiently large. 

Let $\tau_1 := -\frac{1}{200\eta\Lambda}$. In view of the preceding discussion, we may extend $(M^\infty,g^\infty)$ backwards in time to a complete solution $(M^\infty,g^\infty(t))$ which is defined for $t \in [\tau_1,0]$ and satisfies $\Lambda_1 := \sup_{t \in [\tau_1,0]} \sup_{M^\infty} \text{\rm scal}_{g^\infty(t)} \leq 2\Lambda$. 

We now repeat this process. Suppose that we can extend $(M^\infty,g^\infty)$ backwards in time to a complete solution $(M^\infty,g^\infty(t))$ which is defined for $t \in [\tau_m,0]$, and satisfies $\Lambda_m := \sup_{t \in [\tau_m,0]} \sup_{M^\infty} \text{\rm scal}_{g^\infty(t)} < \infty$. Let $\tau_{m+1} := \tau_m - \frac{1}{200\eta \Lambda_m}$. We claim that we can extend the solution $(M^\infty,g^\infty(t))$ backward to the interval $[\tau_{m+1},0]$, and $\Lambda_{m+1} := \sup_{t \in [\tau_{m+1},0]} \sup_{M^\infty} \text{\rm scal}_{g^\infty(t)} \leq 2\Lambda_m$.

Indeed, if this is not possible, then there exists a number $A>1$ with the property that $P(x_j,t_j,A Q_j^{-\frac{1}{2}},(\tau_m-\frac{1}{100\eta \Lambda_m}) \, Q_j^{-1})$ contains points modified by surgery for $j$ sufficiently large. Let $s_j \in [0,\frac{1}{100\eta\Lambda_m}]$ be the largest number such that $P(x_j,t_j,A Q_j^{-\frac{1}{2}},(\tau_m-s_j) Q_j^{-1})$ is free of surgeries. The pointwise curvature derivative estimate gives 
\[\sup_{P(x_j,t_j,A Q_j^{-\frac{1}{2}},(\tau_m-s_j) Q_j^{-1})} \text{\rm scal} \leq 2\Lambda_m \, Q_j\] 
if $j$ is sufficiently large. Let us choose $\alpha_m \in [0,\frac{n-1}{2})$ so that $\frac{\alpha_m}{K_{\text{\rm stnd}} \, (1-\frac{2\alpha_m}{n-1})} \geq 8\Lambda_m \, (\frac{1}{100\eta\Lambda_m}-\tau_m)$. If $(s_j-\tau_m) Q_j^{-1} \geq \alpha_m h_j^2$ for $j$ sufficiently large, then Proposition \ref{evolution.of.surgical.cap} together with the lower bound for the scalar curvature on the standard solution (cf. Theorem \ref{properties.of.standard.solution}) gives 
\begin{align*} 
\sup_{P(x_j,t_j,A Q_j^{-\frac{1}{2}},(\tau_m-s_j) Q_j^{-1})} \text{\rm scal} 
&\geq \frac{1}{2 \, K_{\text{\rm stnd}} \, (1-\frac{2\alpha_m}{n-1})} \, h_j^{-2} \\ 
&\geq \frac{\alpha_m}{2 \, K_{\text{\rm stnd}} \, (1-\frac{2\alpha_m}{n-1})} \, (s_j-\tau_m)^{-1} \, Q_j \\ 
&\geq 4\Lambda_m \, Q_j 
\end{align*}
for $j$ sufficiently large, which is impossible. Consequently, $(s_j-\tau_m) Q_j^{-1} \leq \alpha_m h_j^2$ for $j$ sufficiently large. Since $\hat{\delta}_j \to 0$, Proposition \ref{evolution.of.surgical.cap} implies that the parabolic neighborhood $P(x_j,t_j,A Q_j^{-\frac{1}{2}},(\tau_m-s_j) Q_j^{-1})$ is, after dilating by the factor $h_j$, arbitrarily close to a piece of the standard solution when $j$ is sufficiently large. Using Corollary \ref{canonical.neighborhood.theorem.for.standard.solution.2}, we conclude that $(x_j,t_j)$ lies on a $2\varepsilon$-neck or a $2\varepsilon$-cap when $j$ is sufficiently large. If $(x_j,t_j)$ lies on an $2\varepsilon$-neck, then this neck is actually a strong $2\varepsilon$-neck, since we are assuming that each $\hat{\delta}_j$-neck on which we perform surgery has a large backward parabolic neighborhood that is free of surgeries (cf. Definition \ref{definition.ricci.flow.with.surgery}). Moreover, Corollary \ref{canonical.neighborhood.theorem.for.standard.solution.2} implies that $|D \text{\rm scal}| \leq 2\eta \, \text{\rm scal}^{\frac{3}{2}}$ and $|\frac{\partial}{\partial t} \text{\rm scal}| \leq 2\eta \, \text{\rm scal}^2$ at the point $(x_j,t_j)$. Therefore, the point $(x_j,t_j)$ satisfies the conclusion of Theorem \ref{canonical.neighborhood.property.improved.accuracy}, and this contradicts statement (iii). Thus, we may extend the flow $(M^\infty,g^\infty(t))$ backward to the interval $[\tau_{m+1},0]$, where $\tau_{m+1} := \tau_m - \frac{1}{200\eta\Lambda_m}$, and we have $\Lambda_{m+1} := \sup_{t \in [\tau_{m+1},0]} \sup_{M^\infty} \text{\rm scal}_{g^\infty(t)} \leq 2\Lambda_m$.

Let $\tau^* = \lim_{m \to \infty} \tau_m \leq -\frac{1}{100\eta\Lambda}$. Using a standard diagonal sequence argument, we obtain a complete, smooth limit flow $(M^\infty,g^\infty(t))$ which is defined on the interval $(\tau^*,0]$ and which has bounded curvature for each $t \in (\tau^*,0]$. 

\textit{Step 6:} We next show that $\tau^* = -\infty$. Indeed, if $\tau^* > -\infty$, then $\lim_{m \to \infty} (\tau_m-\tau_{m+1}) = 0$, hence $\lim_{m \to \infty} \Lambda_m = \infty$. Arguing as in Step 6 in the proof of Theorem \ref{high.curvature.regions.1}, we can show that the limit flow $(M^\infty,g^\infty(t))$, $t \in (\tau^*,0]$ has bounded curvature. This contradicts the fact that $\lim_{m \to \infty} \Lambda_m = \infty$. Thus, $\tau^* = -\infty$. Hence, if we dilate the flow $(M,g(t))$ around the point $(x_j,t_j)$ by the factor $Q_j$, then (after passing to a subsequence) the rescaled flows converge to an ancient solution which is complete; has bounded curvature; is weakly PIC2; and satisfies $R - \theta \, \text{\rm scal} \, \text{\rm id} \owedge \text{\rm id} \in \text{\rm PIC}$. By Proposition \ref{noncollapsing.for.ricci.flows.with.surgery}, the limiting ancient solution is $\kappa$-noncollapsed for some $\kappa>0$ which depends only on the initial data. By Corollary \ref{structure.theorem.refined.version}, the point $(x_j,t_j)$ has a Canonical Neighborhood which is either a strong $2\varepsilon$-neck with center at $x_j$; or a $2\varepsilon$-cap; or a closed manifold diffeomorphic to $S^n/\Gamma$; or a quotient neck. Recall that we have discarded all connected components which are diffeomorphic to $S^n/\Gamma$ (see Definition \ref{definition.ricci.flow.with.surgery}). Hence, the Canonical Neighborhood of $(x_j,t_j)$ cannot be a closed manifold diffeomorphic to $S^n/\Gamma$. If the Canonical Neighborhood of $(x_j,t_j)$ is a quotient neck, then Theorem A.1 in \cite{Brendle4} implies that the underlying manifold contains a non-trivial incompressible $(n-1)$-dimensional space form, contrary to our assumption. Consequently, the point $(x_j,t_j)$ has a Canonical Neighborhood which is either a strong $2\varepsilon$-neck with center at $x_j$ or a $2\varepsilon$-cap. Finally, Corollary \ref{pointwise.derivative.estimate.with.universal.constant} implies that $|D \text{\rm scal}| \leq 2\eta \, \text{\rm scal}^{\frac{3}{2}}$ and $|\frac{\partial}{\partial t} \text{\rm scal}| \leq 2\eta \, \text{\rm scal}^2$ at the point $(x_j,t_j)$. In summary, we have shown that the point $(x_j,t_j)$ satisfies the conclusion of Theorem \ref{canonical.neighborhood.property.improved.accuracy}. This contradicts statement (iii). This completes the proof of Theorem \ref{canonical.neighborhood.property.improved.accuracy}.

\section{Global existence of surgically modified flows}

\label{global.existence}

As in the preceding sections, we fix a compact initial manifold $(M,g_0)$ of dimension $n \geq 12$ which has positive isotropic curvature and does not contain any non-trivial incompressible $(n-1)$-dimensional space forms. In this section, we will show that there exists a Ricci flow with surgery starting from $(M,g_0)$, which exists globally and becomes extinct in finite time. We begin by finalizing our choice of the surgery parameters. As above, we fix a small number $\varepsilon>0$. Having chosen $\varepsilon$, we choose numbers $\hat{r},\hat{\delta}$ such that the conclusion of Theorem \ref{canonical.neighborhood.property.improved.accuracy} holds. Having chosen $\varepsilon,\hat{r},\hat{\delta}$, we choose $h$ so that the following holds: 

\begin{proposition}[cf. G.~Perelman \cite{Perelman2}, Lemma 4.3]
\label{existence.of.fine.necks}
Given $\varepsilon,\hat{r},\hat{\delta}$, we can find a small number $h \in (0,\hat{\delta}\hat{r})$ with the following property. Suppose that we have a Ricci flow with surgery with parameters $\varepsilon,\hat{r},\hat{\delta},h$ which is defined on the time interval $[0,T)$ and goes singular at time $T$. Let $x$ be a point which lies in an $4\varepsilon$-horn in $(\Omega,g(T))$ and has curvature $\text{\rm scal}(x,T) = h^{-2}$. Then the parabolic neighborhood $P(x,T,\hat{\delta}^{-1} h,-\hat{\delta}^{-1} h^2)$ is free of surgeries. Moreover, $P(x,T,\hat{\delta}^{-1} h,-\hat{\delta}^{-1} h^2)$ is a $\hat{\delta}$-neck.
\end{proposition}

\textbf{Proof.} 
Suppose that the assertion is false. Then there exists a sequence of positive numbers $h_j \to 0$, a sequence of Ricci flows with surgery $\mathcal{M}^{(j)}$ and a sequence of points $x_j$ with the following properties:
\begin{itemize} 
\item[(i)] The flow $\mathcal{M}^{(j)}$ has surgery parameters $\varepsilon,\hat{r},h_j,\hat{\delta}$. It is defined on the time interval $[0,T_j)$ and goes singular as $t \to T_j$.
\item[(ii)] The point $(x_j,T_j)$ lies on an $4\varepsilon$-horn and $\text{\rm scal}(x_j,T_j) = h_j^{-2}$.
\item[(iii)] The parabolic neighborhood $P(x_j,T_j,\hat{\delta}^{-1} h_j,-\hat{\delta}^{-1} h_j^2)$ contains points modified by surgery, or it is not a $\hat{\delta}$-neck.
\end{itemize} 
By Definition \ref{definition.ricci.flow.with.surgery}, we have $|D \text{\rm scal}| \leq 4\eta \, \text{\rm scal}^{\frac{3}{2}}$ whenever $\text{\rm scal} \geq \hat{r}^{-2}$. Since $h_j \to 0$, it follows that $\inf_{x \in B_{g(T_j)}(x_j,Ah_j)} \text{\rm scal}(x,T_j) \geq (1+2\eta A)^{-2} \, h_j^{-2}$. In particular, if $j$ is sufficiently large (depending on $A$), then $\inf_{x \in B_{g(T_j)}(x_j,Ah_j)} \text{\rm scal}(x,T_j) \geq 10 \, (\hat{\delta}\hat{r})^{-2}$. 

We claim that for each $A>1$ there exists a constant $Q(A)$ (depending on $A$, but not on $j$) such that $\sup_{x \in B_{g(T_j)}(x_j,Ah_j)} \text{\rm scal}(x,T_j) \leq Q(A) \, h_j^{-2}$ if $j$ is sufficiently large. The proof of this statement is analogous to Claim 2 in Theorem 12.1 in \cite{Perelman1}. Indeed, if no such constant $Q(A)$ exists, then we can find a sequence of points $(y_j,T_j)$ (lying on the same horn as $(x_j,T_j)$) such that the blow-up limit around $(y_j,T_j)$ is a piece of non-flat metric cone. Since the flow $\mathcal{M}^{(j)}$ satisfies the Canonical Neighborhood Assumption with accuracy $4\varepsilon$, the point $(y_j,T_j)$ either lies on a strong $4\varepsilon$-neck or on a $4\varepsilon$-cap. The second case can easily be ruled out, so $(y_j,T_j)$ must lie on a strong $4\varepsilon$-neck. In particular, there exists a small parabolic neighborhood of $(y_j,T_j)$ which is free of surgeries. In view of Proposition \ref{curvature.pinching.for.surgically.modified.flows}, the blow-up limit around $(y_j,T_j)$ is weakly PIC2. Hence, Proposition \ref{splitting.a} implies that the limit cannot be a piece of a non-flat metric cone. This proves the claim.

In particular, if $j$ is sufficiently large (depending on $A$), then the distance of the point $x_j$ from either end of the horn is at least $Ah_j$.

We now continue with the proof of Proposition \ref{existence.of.fine.necks}. Let us fix a number $A>1$. Since the point $(x_j,T_j)$ lies on a $4\varepsilon$-horn, no point in $B_{g(T_j)}(x_j,Ah_j)$ can lie on a $4\varepsilon$-cap. Hence, the Canonical Neighborhood Assumption implies that every point in $B_{g(T_j)}(x_j,Ah_j)$ lies on a strong $4\varepsilon$-neck. Using Shi's estimate, we obtain bounds for all the covariant derivatives of the curvature tensor in $B_{g(T_j)}(x_j,\frac{1}{2} Ah_j)$. Note that these bounds may depend on $A$, but are independent of $j$. We now pass to the limit, sending $j \to \infty$ first and $A \to \infty$ second. In the limit, we obtain a complete manifold with two ends which, by Proposition \ref{curvature.pinching.for.surgically.modified.flows}, is uniformly PIC and weakly PIC2. By the Cheeger-Gromoll splitting theorem, the limit is isometric to a product $X \times \mathbb{R}$; moreover, the cross-section $X$ is compact and is nearly isometric to $S^{n-1}$. Since every point in $B_{g(T_j)}(x_j,Ah_j)$ lies on a strong $4\varepsilon$-neck, we conclude that, for each $A>1$, the parabolic neighborhood $P(x_j,T_j,Ah_j,-\frac{3h_j^2}{4})$ is free of surgeries if $j$ is sufficiently large (depending on $A$). After rescaling and passing to the limit, we obtain a solution to the Ricci flow which is defined on the time interval $[-\frac{1}{2},0]$ and which splits off a line. Now, if $j$ is sufficiently large (depending on $A$), then no point in $P(x_j,T_j,Ah_j,-\frac{h_j^2}{2})$ can lie on a $4\varepsilon$-cap. Hence, if $j$ is sufficiently large (depending on $A$), then every point in the parabolic neighborhood $P(x_j,T_j,Ah_j,-\frac{h_j^2}{2})$ lies on a strong $4\varepsilon$-neck. This allows us to extend the limit solution backward in time to the interval $[-1,0]$. Repeating this argument, we can extend the limit solution backwards in time, so that it is defined on $[-1,0]$, $[-\frac{3}{2},0]$, $[-2,0]$, etc. To summarize, we produce a limit solution which is ancient; uniformly PIC; weakly PIC2; and splits as a product of a line with a manifold diffeomorphic to $S^{n-1}$. By Theorem \ref{ancient.solutions.bounded.curvature.uniformly.pic1}, the limiting solution is a family of standard cylinders. Therefore, if $j$ is sufficiently large, then the parabolic neighborhood $P(x_j,T_j,\hat{\delta}^{-1} h_j,-\hat{\delta}^{-1} h_j^2)$ is free of surgeries, and $P(x_j,T_j,\hat{\delta}^{-1} h_j,-\hat{\delta}^{-1} h_j^2)$ is a $\hat{\delta}$-neck. This contradicts (iii). \\

We are now able to prove the main result of this section:

\begin{theorem}
\label{global.existence.of.flow.with.surgery}
Let us fix a small number $\varepsilon>0$. Let $\hat{r},\hat{\delta}$ be chosen as described at the beginning of this section, and let $h$ be chosen as in Proposition \ref{existence.of.fine.necks}. Then there exists a Ricci flow with surgery with parameters $\varepsilon,\hat{r},\hat{\delta},h$, which is defined on some finite time interval $[0,T)$ and becomes extinct as $t \to T$.
\end{theorem}

\textbf{Proof.} 
We evolve the initial metric $g_0$ by smooth Ricci flow until the flow becomes singular for the first time. It follows from Theorem \ref{canonical.neighborhood.property.improved.accuracy} and a standard continuity argument that the flow satisfies the Canonical Neighborhood Property with accuracy $2\varepsilon$ on all scales less than $2\hat{r}$, up until the first singular time. At the first singular time, we perform finitely many surgeries on $\hat{\delta}$-necks which have curvature level $h^{-2}$. The existence of such necks is ensured by Proposition \ref{existence.of.fine.necks}. After performing surgery, we restart the flow, and continue until the second singular time. Using Theorem \ref{canonical.neighborhood.property.improved.accuracy} and a standard continuity argument, we conclude that the flow with surgery satisfies the Canonical Neighborhood Property with accuracy $2\varepsilon$ on all scales less than $2\hat{r}$, up until the second singular time. Consequently, Proposition \ref{existence.of.fine.necks} ensures that, at the second singular time, we can again find $\hat{\delta}$-necks on which to perform surgery. After performing surgery, we continue the flow until the third singular time. Again, Theorem \ref{canonical.neighborhood.property.improved.accuracy} guarantees that the flow with surgery satisfies the Canonical Neighborhood Property with accuracy $2\varepsilon$ on all scales less than $2\hat{r}$, up until the third singular time. We can now perform surgery again, and repeat the process. 

Since each surgery reduces the volume by at least $c(n) \, h^n$, we have an upper bound for the number of surgeries. By Proposition \ref{finite.time.extinction}, the flow with surgery must become extinct by time $\frac{n}{2 \, \inf_{x \in M} \text{\rm scal}(x,0)}$ at the latest. This completes the proof of Theorem \ref{global.existence.of.flow.with.surgery}. \\

\begin{corollary}
\label{topological.classification}
The manifold $M$ is diffeomorphic to a connected sum of finitely many spaces, each of which is a quotient of $S^n$ or $S^{n-1} \times \mathbb{R}$ by standard isometries.
\end{corollary}

\textbf{Proof.} 
This follows by combining Theorem \ref{global.existence.of.flow.with.surgery} with Proposition \ref{change.of.topology}.

\appendix 

\section{Auxiliary results}

\label{aux}

In this section, we collect various technical results which are needed in Sections \ref{preliminary.pinching.estimate} -- \ref{second.family.of.cones}. We assume throughout that $n \geq 5$.

\begin{lemma}
\label{algebraic.fact}
Let $0 \leq \zeta \leq 1$ and $0 < \rho \leq 1$. Assume that $H$ is a symmetric bilinear form with the property that the largest eigenvalue of $H$ is bounded from above by $\frac{1}{2} \, \text{\rm tr}(H)$, and the sum of the two smallest eigenvalues of $H$ is bounded from below by $\frac{2(1-\zeta)}{n} \, \text{\rm tr}(H)$. Then 
\begin{align*} 
&\frac{n-2}{n} \, \text{\rm tr}(H) \, (H_{11}+H_{22}) - \rho \, ((\tracefreeH^2)_{11}+(\tracefreeH^2)_{22}) \\ 
&\geq \frac{2}{n^2} \, \Big ( (n-2)(1-\zeta) - 2\zeta^2\rho \, \frac{n^2-2n+2}{(n-2)^2} \Big ) \, \text{\rm tr}(H)^2 
\end{align*} 
for every pair of orthonormal vectors $\{e_1,e_2\}$, where $\tracefreeH$ denotes the tracefree part of $H$.
\end{lemma}

\textbf{Proof.} 
By scaling, we may assume that $\text{\rm tr}(H)=n$. Moreover, we may assume that $H$ is diagonal with diagonal entries $\lambda_1 \leq \hdots \leq \lambda_n$. By assumption, each eigenvalue is bounded from above by $\frac{n}{2}$, and the sum of any two distinct eigenvalues of $H$ is bounded from below by $2(1-\zeta)$. It suffices to prove that 
\[(n-2) \, (\lambda_i+\lambda_j) - \rho \, ((\lambda_i-1)^2+(\lambda_j-1)^2) \geq 2(n-2)(1-\zeta) - 4\zeta^2 \rho \, \frac{n^2-2n+2}{(n-2)^2}\] 
or, equivalently, 
\[\lambda_i \, \Big ( \frac{n-2}{\rho}+2-\lambda_i \Big ) + \lambda_j \, \Big ( \frac{n-2}{\rho}+2-\lambda_j \Big ) \geq 2(1-\zeta) \, \Big ( \frac{n-2}{\rho}+1+\zeta \Big ) - \frac{2n^2\zeta^2}{(n-2)^2}.\] 
for $i<j$. Note that $\frac{n-2}{\rho}+2 \geq n$. Moreover, $\lambda_n \leq \frac{n}{2}$ by assumption. We distinguish two cases:

\textit{Case 1:} Suppose first that $\lambda_i \geq 1-\zeta$. Since $i<j$, we have $1-\zeta \leq \lambda_i \leq \lambda_j \leq \frac{n}{2}$. Since the function $\lambda \mapsto \lambda \, \big ( \frac{n-2}{\rho}+2-\lambda \big )$ is monotone increasing on the interval $(-\infty,\frac{n}{2}]$, we obtain 
\begin{align*} 
&\lambda_i \, \Big ( \frac{n-2}{\rho}+2-\lambda_i \Big ) + \lambda_j \, \Big ( \frac{n-2}{\rho}+2-\lambda_j \Big ) \\ 
&\geq 2(1-\zeta) \, \Big ( \frac{n-2}{\rho}+1+\zeta \Big ), 
\end{align*} 
which implies the claim.

\textit{Case 2:} Suppose finally that $\lambda_i \leq 1-\zeta$. Since the sum of any two distinct eigenvalues is at least $2(1-\zeta)$, we obtain $2(1-\zeta)-\lambda_i \leq \lambda_j \leq \frac{n}{2}$. Since the function $\lambda \mapsto \lambda \, \big ( \frac{n-2}{\rho}+2-\lambda \big )$ is monotone increasing on the interval $(-\infty,\frac{n}{2}]$, we obtain   
\begin{align*} 
&\lambda_i \, \Big ( \frac{n-2}{\rho}+2-\lambda_i \Big ) + \lambda_j \, \Big ( \frac{n-2}{\rho}+2-\lambda_j \Big ) \\ 
&\geq \lambda_i \, \Big ( \frac{n-2}{\rho}+2-\lambda_i \Big ) + (2(1-\zeta)-\lambda_i) \, \Big ( \frac{n-2}{\rho}+2\zeta+\lambda_i \Big ) \\ 
&= 2(1-\zeta) \, \Big ( \frac{n-2}{\rho}+1+\zeta \Big ) - 2(1-\zeta-\lambda_i)^2. 
\end{align*} 
Since the sum of any two distinct eigenvalues is at least $2(1-\zeta)$, it follows that $\lambda_i \geq 1-\frac{2(n-1)\zeta}{n-2}$. Consequently, $0 \leq 1-\zeta-\lambda_i \leq \frac{n\zeta}{n-2}$. Thus, we conclude that 
\begin{align*} 
&\lambda_i \, \Big ( \frac{n-2}{\rho}+2-\lambda_i \Big ) + \lambda_j \, \Big ( \frac{n-2}{\rho}+2-\lambda_j \Big ) \\ 
&\geq 2(1-\zeta) \, \Big ( \frac{n-2}{\rho}+1+\zeta \Big ) - \frac{2n^2\zeta^2}{(n-2)^2}, 
\end{align*} 
as claimed. \\

\begin{lemma}
\label{wedge.product}
Let $H$ be a symmetric bilinear form. If $H$ is weakly two-positive, then $H \owedge H \in \text{\rm PIC}$.
\end{lemma}

\textbf{Proof.} 
Let $\zeta,\eta \in \mathbb{C}^n$ be linearly independent vectors satisfying $g(\zeta,\zeta)=g(\zeta,\eta)=g(\eta,\eta)=0$. We claim that $(H \owedge H)(\zeta,\eta,\bar{\zeta},\bar{\eta}) \geq 0$. We can find vectors $z,w \in \text{\rm span}\{\zeta,\eta\}$ such that $g(z,\bar{z})=g(w,\bar{w})=2$, $g(z,\bar{w})=0$, and $H(z,\bar{w}) = 0$. The identities $g(\zeta,\zeta)=g(\zeta,\eta)=g(\eta,\eta)=0$ give $g(z,z)=g(z,w)=g(w,w)=0$. Consequently, we may write $z=e_1+ie_2$ and $w=e_3+ie_4$ for some orthonormal four-frame $\{e_1,e_2,e_3,e_4\} \subset \mathbb{R}^n$. Using the identity $H(z,\bar{w}) = 0$, we obtain 
\[(H \owedge H)(z,w,\bar{z},\bar{w}) = 2 \, (H_{11}+H_{22})(H_{33}+H_{44}) \geq 0.\] 
Since $\text{\rm span}\{\zeta,\eta\} = \text{\rm span}\{z,w\}$, we conclude that $(H \owedge H)(\zeta,\eta,\bar{\zeta},\bar{\eta}) \geq 0$, as claimed. \\

\begin{lemma}
\label{estimate.for.largest.ricci.eigenvalue}
Suppose that $S \in \text{\rm PIC}$. Then the largest eigenvalue of the Ricci tensor of $S$ is bounded from above by $\frac{1}{2} \, \text{\rm scal}(S)$. 
\end{lemma}

\textbf{Proof.} 
Since $S$ has nonnegative isotropic curvature, we have 
\[\text{\rm scal}(S) - 2 \, \text{\rm Ric}(S)_{nn} = \sum_{k,l=1}^{n-1} S_{klkl} \geq 0,\] 
as claimed. \\

\begin{lemma}
\label{Ric.four.positive}
Assume $S \in \text{\rm PIC}$. Then $\text{\rm Ric}(S)_{11}+\text{\rm Ric}(S)_{22}+\text{\rm Ric}(S)_{33}+\text{\rm Ric}(S)_{44} \geq 0$ for every orthonormal four-frame $\{e_1,e_2,e_3,e_4\}$. 
\end{lemma}

\textbf{Proof.} 
The condition $S \in \text{\rm PIC}$ implies $\text{\rm Ric}(S)_{11}+\text{\rm Ric}(S)_{33} \geq 2S_{1313}$, $\text{\rm Ric}(S)_{11}+\text{\rm Ric}(S)_{44} \geq 2S_{1414}$, $\text{\rm Ric}(S)_{22}+\text{\rm Ric}(S)_{33} \geq 2S_{2323}$, $\text{\rm Ric}(S)_{22}+\text{\rm Ric}(S)_{44} \geq 2S_{2424}$. Taking the sum of all four inequalities yields
\begin{align*} 
&\text{\rm Ric}(S)_{11}+\text{\rm Ric}(S)_{22}+\text{\rm Ric}(S)_{33}+\text{\rm Ric}(S)_{44} \\ 
&\geq S_{1313}+S_{1414}+S_{2323}+S_{2424} \geq 0, 
\end{align*}
as claimed. \\

%

\begin{proposition} 
\label{minimum.isotropic.curvature}
Suppose that $R(t)$ is a solution of the Hamilton ODE $\frac{d}{dt} R = Q(R)$, and $\kappa(t)$ is a nonnegative function satisfying $\frac{d}{dt} \kappa(t) \leq 2(n-1)\kappa(t)^2$. Then the condition $R(t) - \frac{1}{2} \, \kappa(t) \, \text{\rm id} \owedge \text{\rm id} \in \text{\rm PIC}$ is preserved. 
\end{proposition}

\textbf{Proof.} 
Let $S(t) = R(t) - \frac{1}{2} \, \kappa(t) \, \text{\rm id} \owedge \text{\rm id}$. As in \cite{Brendle3}, we compute 
\[Q(S(t)) + 2\kappa(t) \, \text{\rm Ric}(S) \owedge \text{\rm id} + (n-1) \, \kappa(t)^2 \, \text{\rm id} \owedge \text{\rm id} = Q(R(t)).\] 
This gives 
\[\frac{d}{dt} S(t) = Q(S(t)) + 2\kappa(t) \, \text{\rm Ric}(S(t)) \owedge \text{\rm id} - \frac{1}{2} \, \Big ( \frac{d}{dt} \kappa(t) - 2(n-1)\kappa(t)^2 \Big ) \, \text{\rm id} \owedge \text{\rm id}.\] 
If $S \in \text{\rm PIC}$, then Proposition 7.5 in \cite{Brendle-book} implies $Q(S) \in T_S \text{\rm PIC}$. Moreover,if $S \in \text{\rm PIC}$, then the sum of the four smallest eigenvalues of $\text{\rm Ric}(S)$ is nonnegative by Lemma \ref{Ric.four.positive}, and consequently $\text{\rm Ric}(S) \owedge \text{\rm id} \in \text{\rm PIC}$. Consequently, the condition $S(t) \in \text{\rm PIC}$ is preserved. \\

\begin{lemma}
\label{pic.implies.pic1.with.error.term}
Assume that $S \in \text{\rm PIC}$. Then 
\begin{align*} 
&S_{1313}+\lambda^2 S_{1414}+S_{2323}+\lambda^2 S_{2424}-2\lambda S_{1234} \\ 
&+ \frac{1-\lambda^2}{n-4} \, (\text{\rm Ric}(S)_{11}+\text{\rm Ric}(S)_{22}-2S_{1212}) \geq 0 
\end{align*} 
for every orthonormal four-frame $\{e_1,e_2,e_3,e_4\}$ and every $\lambda \in [0,1]$.
\end{lemma}

\textbf{Proof.} 
Since $S \in \text{\rm PIC}$, we have 
\[S_{1313}+S_{1414}+S_{2323}+S_{2424}-2S_{1234} \geq 0\] 
for every orthonormal four-frame $\{e_1,e_2,e_3,e_4\}$. We now replace $e_4$ by $\lambda e_4 \pm \sqrt{1-\lambda^2} e_p$ for some fixed $p \in \{5,\hdots,n\}$, and take the arithmetic mean of the resulting inequalities. This gives 
\[S_{1313}+\lambda^2 S_{1414}+S_{2323}+\lambda^2 S_{2424}-2\lambda S_{1234} + (1-\lambda^2) (S_{1p1p}+S_{2p2p}) \geq 0\] 
for all $p \in \{5,\hdots,n\}$. In the next step, we take the mean value over all $p \in \{5,\hdots,n\}$. Using the inequality 
\[\sum_{p=5}^n (S_{1p1p}+S_{2p2p}) \leq \sum_{p=3}^n (S_{1p1p}+S_{2p2p}) = \text{\rm Ric}(S)_{11}+\text{\rm Ric}(S)_{22}-2S_{1212},\] 
we obtain 
\begin{align*} 
&S_{1313}+\lambda^2 S_{1414}+S_{2323}+\lambda^2 S_{2424}-2\lambda S_{1234} \\ 
&+ \frac{1-\lambda^2}{n-4} \, (\text{\rm Ric}(S)_{11}+\text{\rm Ric}(S)_{22}-2S_{1212}) \geq 0. 
\end{align*} 
This proves the assertion. \\

\begin{lemma}
\label{curvature.tensor.in.one.extra.dimension}
Let $S$ be an algebraic curvature tensor and let $H$ be a symmetric bilinear form on $\mathbb{R}^n$. Let us define an algebraic curvature tensor $T$ on $\mathbb{R}^{n+1}$ by 
\[T_{ijkl} = S_{ijkl}, \quad T_{ijk0} = 0, \quad T_{i0k0} = H_{ik}\] 
for $i,j,k,l \in \{1,\hdots,n\}$. If 
\[S_{1313}+\lambda^2 S_{1414}+S_{2323}+\lambda^2 S_{2424}-2\lambda S_{1234} + (1-\lambda^2) \, (H_{11}+H_{22}) \geq 0\] 
for every orthonormal four-frame $\{e_1,e_2,e_3,e_4\}$ and every $\lambda \in [0,1]$, then $T \in \text{\rm PIC}$.
\end{lemma} 

\textbf{Proof.} 
Let $\tilde{\zeta},\tilde{\eta} \in \mathbb{C}^{n+1}$ be two linearly independent vectors satisfying $g(\tilde{\zeta},\tilde{\zeta})=g(\tilde{\zeta},\tilde{\eta})=g(\tilde{\eta},\tilde{\eta})=0$. We claim that $T(\tilde{\zeta},\tilde{\eta},\bar{\tilde{\zeta}},\bar{\tilde{\eta}}) \geq 0$. We can find a non-zero vector $z \in \mathbb{C}^n$ such that $\tilde{z} = (z,0) \in \text{\rm span}\{\tilde{\zeta},\tilde{\eta}\}$. Moreover, we can find a non-zero vector $\tilde{w} \in \mathbb{C}^{n+1}$ such that $\tilde{w} \in \text{\rm span}\{\tilde{\zeta},\tilde{\eta}\}$ and $g(\tilde{z},\bar{\tilde{w}})=0$. We normalize $\tilde{z}$ and $\tilde{w}$ such that $g(\tilde{z},\bar{\tilde{z}})=g(\tilde{w},\bar{\tilde{w}})=2$. We can further arrange that $\tilde{w} = (w,ia)$, where $w \in \mathbb{C}^n$ and $a \in [0,\infty)$. The identities $g(\tilde{\zeta},\tilde{\zeta})=g(\tilde{\zeta},\tilde{\eta})=g(\tilde{\eta},\tilde{\eta})=0$ imply $g(\tilde{z},\tilde{z})=g(\tilde{z},\tilde{w})=g(\tilde{w},\tilde{w})=0$. Consequently, we may write $z=e_1+ie_2$, $w=e_3+i\lambda e_4$, and $a=\sqrt{1-\lambda^2}$ for some orthonormal four-frame $\{e_1,e_2,e_3,e_4\} \subset \mathbb{R}^n$ and some $\lambda \in [0,1]$. This implies 
\begin{align*} 
&T(\tilde{z},\tilde{w},\bar{\tilde{z}},\bar{\tilde{w}}) \\ 
&= S_{1313}+\lambda^2 S_{1414}+S_{2323}+\lambda^2 S_{2424}-2\lambda S_{1234} + (1-\lambda^2) \, (H_{11}+H_{22}) \geq 0. 
\end{align*}
Since $\text{\rm span}\{\tilde{\zeta},\tilde{\eta}\}=\text{\rm span}\{\tilde{z},\tilde{w}\}$, it follows that $T(\tilde{\zeta},\tilde{\eta},\bar{\tilde{\zeta}},\bar{\tilde{\eta}}) \geq 0$, as claimed. \\

\begin{proposition}
\label{interpolation.between.pic.and.pic1}
Let $S$ be an algebraic curvature tensor and let $H$ be a symmetric bilinear form on $\mathbb{R}^n$ such that 
\[Z := S_{1313}+\lambda^2 S_{1414}+S_{2323}+\lambda^2 S_{2424}-2\lambda S_{1234} + (1-\lambda^2) \, (H_{11}+H_{22}) \geq 0\] 
for every orthonormal four-frame $\{e_1,e_2,e_3,e_4\}$ and every $\lambda \in [0,1]$. Moreover, suppose that $Z=0$ for one particular orthonormal four-frame $\{e_1,e_2,e_3,e_4\}$ and one particular $\lambda \in [0,1)$. Then 
\begin{align*} 
&Q(S)_{1313}+\lambda^2 Q(S)_{1414}+Q(S)_{2323}+\lambda^2 Q(S)_{2424}-2\lambda Q(S)_{1234} \\ 
&+ (H \owedge H)_{1313}+\lambda^2 (H \owedge H)_{1414} \\ 
&+ (H \owedge H)_{2323}+\lambda^2 (H \owedge H)_{2424} \\ 
&- 2\lambda (H \owedge H)_{1234} \\ 
&+ 2(1-\lambda^2) \, ((S * H)_{11}+(S * H)_{22}) \\ 
&\geq (1+\lambda^2) (H_{11}+H_{22})^2 
\end{align*} 
for this particular four-frame $\{e_1,e_2,e_3,e_4\}$ and this particular $\lambda \in [0,1)$. Here, $S * H$ is defined by $(S * H)_{ik} := \sum_{p,q=1}^n S_{ipkq} H_{pq}$.
\end{proposition} 

\textbf{Proof.} 
We define a curvature tensor $T$ on $\mathbb{R}^{n+1}$ by 
\[T_{ijkl} = S_{ijkl}, \quad T_{ijk0} = 0, \quad T_{i0k0} = H_{ik}\] 
for $i,j,k,l \in \{1,\hdots,n\}$. By Lemma \ref{curvature.tensor.in.one.extra.dimension}, we know that $T \in \text{\rm PIC}$. Suppose that $Z=0$ for one particular orthonormal four-frame $\{e_1,e_2,e_3,e_4\} \subset \mathbb{R}^n$ and one particular $\lambda \in [0,1)$. Then $\{e_1,e_2,e_3,\lambda e_4 + \sqrt{1-\lambda^2} e_0\}$ is an orthonormal four-frame in $\mathbb{R}^{n+1}$ which has zero isotropic curvature for $T$. Using Proposition 7.4 in \cite{Brendle-book}, we obtain 
\begin{align*} 
0 &\leq \sum_{p,q=0}^n (T_{1p1q} + T_{2p2q}) \\ 
&\hspace{10mm} \cdot (T_{3p3q} + \lambda^2 T_{4p4q} + \lambda \sqrt{1-\lambda^2} T_{4p0q} + \lambda \sqrt{1-\lambda^2} T_{0p4q} + (1-\lambda^2) T_{0p0q}) \\ 
&- \sum_{p,q=0}^n T_{12pq} \, (\lambda T_{34pq} + \sqrt{1-\lambda^2} T_{30pq}) \\ 
&- \sum_{p,q=0}^n (T_{1p3q} + \lambda T_{2p4q} + \sqrt{1-\lambda^2} T_{2p0q}) \, (T_{3p1q} + \lambda T_{4p2q} + \sqrt{1-\lambda^2} T_{0p2q}) \\ 
&- \sum_{p,q=0}^n (\lambda T_{1p4q} + \sqrt{1-\lambda^2} T_{1p0q} - T_{2p3q}) \, (\lambda T_{4p1q} + \sqrt{1-\lambda^2} T_{0p1q} - T_{3p2q}) \\ 
&= \sum_{p,q=0}^n (T_{1p1q} + T_{2p2q}) \, (T_{3p3q} + \lambda^2 T_{4p4q} + (1-\lambda^2) T_{0p0q}) \\ 
&- \sum_{p,q=0}^n T_{12pq} \, \lambda T_{34pq} \\ 
&- \sum_{p,q=0}^n (T_{1p3q} + \lambda T_{2p4q}) \, (T_{3p1q} + \lambda T_{4p2q}) \\ 
&- \sum_{p,q=0}^n (\lambda T_{1p4q} - T_{2p3q}) \, (\lambda T_{4p1q} - T_{3p2q}). 
\end{align*} 
This implies 
\begin{align*} 
0 &\leq \sum_{p,q=1}^n (T_{1p1q} + T_{2p2q}) \, (T_{3p3q} + \lambda^2 T_{4p4q} + (1-\lambda^2) T_{0p0q}) \\ 
&- \sum_{p,q=1}^n T_{12pq} \, \lambda T_{34pq} \\ 
&- \sum_{p,q=1}^n (T_{1p3q} + \lambda T_{2p4q}) \, (T_{3p1q} + \lambda T_{4p2q}) \\ 
&- \sum_{p,q=1}^n (\lambda T_{1p4q} - T_{2p3q}) \, (\lambda T_{4p1q} - T_{3p2q}) \\ 
&+ (H_{11}+H_{22}) \, (H_{33}+\lambda^2 H_{44}) \\ 
&- (H_{13} + \lambda H_{24})^2 - (\lambda H_{14} - H_{23})^2 \\ 
&= \sum_{p,q=1}^n (S_{1p1q} + S_{2p2q}) \, (S_{3p3q} + \lambda^2 S_{4p4q} + (1-\lambda^2) H_{pq}) \\ 
&- \sum_{p,q=1}^n S_{12pq} \, \lambda S_{34pq} \\ 
&- \sum_{p,q=1}^n (S_{1p3q} + \lambda S_{2p4q}) \, (S_{3p1q} + \lambda S_{4p2q}) \\ 
&- \sum_{p,q=1}^n (\lambda S_{1p4q} - S_{2p3q}) \, (\lambda S_{4p1q} - S_{3p2q}) \\ 
&+ (H_{11}+H_{22}) \, (H_{33}+\lambda^2 H_{44}) \\ 
&- (H_{13} + \lambda H_{24})^2 - (\lambda H_{14} - H_{23})^2. 
\end{align*}
Since $\lambda \in [0,1)$, we have $Z = \frac{\partial Z}{\partial \lambda} = 0$. This gives $\lambda S_{1414} + \lambda S_{2424} - S_{1234} = \lambda (H_{11}+H_{22})$ and $S_{1313} + S_{2323} - \lambda S_{1234} = -(H_{11}+H_{22})$. Consequently,  
\begin{align*} 
(1+\lambda^2) (H_{11}+H_{22})^2 
&= \big [ (\lambda S_{1414} + S_{1423}) + (\lambda S_{2424} - S_{1324}) \big ]^2 \\ 
&+ \big [ (S_{1313} - \lambda S_{1324}) + (S_{2323} + \lambda S_{1423}) \big ]^2 \\ 
&\leq 2 (\lambda S_{1414} + S_{1423})^2 + 2 (\lambda S_{2424} - S_{1324})^2 \\ 
&+ 2 (S_{1313} - \lambda S_{1324})^2 + 2 (S_{2323} + \lambda S_{1423})^2 \\ 
&\leq \sum_{p,q=1}^n (S_{13pq} - \lambda S_{24pq})^2 + \sum_{p,q=1}^n (\lambda S_{14pq} + S_{23pq})^2. 
\end{align*}
Putting these facts together, we conclude that 
\begin{align*} 
&Q(S)_{1313}+\lambda^2 Q(S)_{1414}+Q(S)_{2323}+\lambda^2 Q(S)_{2424}-2\lambda Q(S)_{1234} \\ 
&= \sum_{p,q=1}^n (S_{13pq} - \lambda S_{24pq})^2 + \sum_{p,q=1}^n (\lambda S_{14pq} + S_{23pq})^2 \\ 
&+ 2 \sum_{p,q=1}^n (S_{1p1q} + S_{2p2q}) \, (S_{3p3q} + \lambda^2 S_{4p4q}) \\ 
&- 2 \sum_{p,q=1}^n S_{12pq} \, \lambda S_{34pq} \\ 
&- 2 \sum_{p,q=1}^n (S_{1p3q} + \lambda S_{2p4q}) \, (S_{3p1q} + \lambda S_{4p2q}) \\ 
&- 2 \sum_{p,q=1}^n (\lambda S_{1p4q} - S_{2p3q}) \, (\lambda S_{4p1q} - S_{3p2q}) \\ 
&\geq (1+\lambda^2) (H_{11}+H_{22})^2 - 2 \, (1-\lambda^2) \sum_{p,q=1}^n (S_{1p1q}+S_{2p2q}) \, H_{pq} \\ 
&- 2 \, (H_{11}+H_{22}) \, (H_{33}+\lambda^2 H_{44}) \\ 
&+ 2 \, (H_{13} + \lambda H_{24})^2 + 2 \, (\lambda H_{14} - H_{23})^2. 
\end{align*} 
This completes the proof. \\

Finally, we consider tensors $S$ and $U$ which do not necessarily satisfy the first Bianchi identity. The following result generalizes Proposition 7.4 in \cite{Brendle-book} (see also \cite{Wilking}): 

\begin{proposition}
\label{sharp}
Assume that $S$ and $U$ are tensors which satisfy all the algebraic properties of a curvature tensor, except for the first Bianchi identity. Moreover, suppose that $U \in \text{\rm PIC}$ and $S-U \geq 0$. Then $S^\# \in T_U \text{\rm PIC}$. 
\end{proposition}

\textbf{Proof.} 
Consider a complex two-form $\varphi$ with the property that $\varphi = (e_1+ie_2) \wedge (e_3+ie_4)$ for some orthonormal four-frame $\{e_1,e_2,e_3,e_4\}$ and $U(\varphi,\bar{\varphi})=0$. Since $U \in \text{\rm PIC}$, $U(e^{-t\sigma} \varphi e^{t\sigma},\overline{e^{-t\sigma} \varphi e^{t\sigma}}) \geq 0$ for all $\sigma \in so(n,\mathbb{C})$ and all $t \in \mathbb{R}$. Since the second derivative at $t=0$ is nonnegative, we obtain 
\[U([[\varphi,\sigma],\sigma],\bar{\varphi}) + U(\varphi,\overline{[[\varphi,\sigma],\sigma]}) + 2 \, U([\varphi,\sigma],\overline{[\varphi,\sigma]}) \geq 0\] 
for all $\sigma \in so(n,\mathbb{C})$ (cf. \cite{Wilking}). If we replace $\sigma$ by $i\sigma$ and add the resulting inequalities, we obtain $U([\varphi,\sigma],\overline{[\varphi,\sigma]}) \geq 0$ for all $\sigma \in so(n,\mathbb{C})$. Since $S-U \geq 0$, it follows that $S([\varphi,\sigma],\overline{[\varphi,\sigma]}) \geq 0$ for all $\sigma \in so(n,\mathbb{C})$. Let us define a linear transformation $L: so(n,\mathbb{C}) \to so(n,\mathbb{C})$ by $L\sigma = [\varphi,\sigma]$. The adjoint $L^*: so(n,\mathbb{C}) \to so(n,\mathbb{C})$ is given by $L^*\sigma = -[\bar{\varphi},\sigma]$. Let $P: so(n,\mathbb{C}) \to so(n,\mathbb{C})$ be a linear transformation with the property that $\sigma = LP\sigma$ whenever $\sigma$ lies in the image of $L$. Clearly, $L = LPL$, and consequently $L^* = L^*P^*L^*$. Since $L^*SL \geq 0$ and $LSL^* \geq 0$, we obtain 
\[\text{\rm tr}(SLSL^*) = \text{\rm tr}((L^*SL)P(LSL^*)P^*) \geq 0.\] 
From this, we deduce that $S^\#(\varphi,\bar{\varphi}) \geq 0$. Since this inequality holds for every two-form $\varphi = (e_1+ie_2) \wedge (e_3+ie_4)$ satisfying $U(\varphi,\bar{\varphi}) = 0$, we conclude that $S^\# \in T_U \text{\rm PIC}$. \\


\begin{thebibliography}{99} 
\bibitem{Bohm-Wilking} 
C.~B\"ohm and B.~Wilking, \textit{Manifolds with positive curvature operator are space forms,} Ann. of Math. 167, 1079--1097 (2008)

\bibitem{Brendle1}
S.~Brendle, \textit{A general convergence result for the Ricci flow in higher dimensions,} Duke Math. J. 145, 585--601 (2008)

\bibitem{Brendle2}
S.~Brendle, \textit{A generalization of Hamilton's differential Harnack inequality for the Ricci flow,} J. Diff. Geom. 82, 207--227 (2009)

\bibitem{Brendle3}
S.~Brendle, \textit{Einstein manifolds with nonnegative isotropic curvature are locally symmetric,} Duke Math. J. 151, 1--21 (2010)

\bibitem{Brendle-book}
S.~Brendle, \textit{Ricci Flow and the Sphere Theorem,} Graduate Studies in Mathematics, vol. 111, American Mathematical Society (2010)

\bibitem{Brendle4}
S.~Brendle, \textit{Ricci flow with surgery in higher dimensions,} Ann. of Math. 187, 263--299 (2018)

\bibitem{Brendle-Huisken-Sinestrari}
S.~Brendle, G.~Huisken, and C.~Sinestrari, \textit{Ancient solutions to the Ricci flow with pinched curvature,} Duke Math. J. 158, 537--551 (2011)

\bibitem{Brendle-Schoen} 
S.~Brendle and R.~Schoen, \textit{Manifolds with $1/4$-pinched curvature are space forms,} J. Amer. Math. Soc. 22, 287--307 (2009)

\bibitem{Cheeger-Gromoll1} 
J.~Cheeger and D.~Gromoll, \textit{The splitting theorem for manifolds of nonnegative Ricci curvature,} J. Diff. Geom. 6, 119--128 (1971)

\bibitem{Cheeger-Gromoll2}
J.~Cheeger and D.~Gromoll, \textit{On the structure of complete manifolds of nonnegative curvature,} Ann. of Math. 96, 413--443 (1972)

\bibitem{Chen-Zhu} 
B.~Chen and X.~Zhu, \textit{Ricci flow with surgery on four-manifolds with positive isotropic curvature,} J. Diff. Geom. 74, 177--264 (2006) 

\bibitem{Chen-Tang-Zhu}
B.~Chen, S.~Tang, and X.~Zhu, \textit{Complete classification of compact four-manifolds with positive isotropic curvature,} J. Diff. Geom. 91, 41--80 (2012)

\bibitem{Chow-Lu}
B.~Chow and P.~Lu, \textit{The maximum principle for systems of parabolic equations subject to an avoidance set,} Pacific J. Math. 214, 201--222 (2004)

\bibitem{Chow-et-al}
B.~Chow, S.-C.~Chu, D.~Glickenstein, C.~Guenther, J.~Isenberg, T.~Ivey, D.~Knopf, P.~Lu, F.~Luo, and L.~Ni, \textit{The Ricci Flow: Techniques and Applications, Part II: Analytic Aspects,} Mathematical Surveys and Monographs vol. 144, American Mathematical Society (2008)

\bibitem{Fraser}
A.M.~Fraser, \textit{Fundamental groups of manifolds with positive isotropic curvature,} Ann. of Math. 158, 345--354 (2003)

\bibitem{Gromov}
M.~Gromov, \textit{Positive curvature, macroscopic dimension, spectral gaps and higher signatures,} in: Functional analysis on the eve of the 21st century, vol. II (New Brunswick, NJ, 1993), pp.~1--213, Progr. Math. 132, Birkh\"auser, 1996

\bibitem{Hamilton1}
R.~Hamilton, \textit{Three-manifolds with positive Ricci curvature,} J. Diff. Geom. 17, 255--306 (1982)

\bibitem{Hamilton2} 
R.~Hamilton, \textit{Four-manifolds with positive curvature operator,} J. Diff. Geom. 24, 153--179 (1986)

\bibitem{Hamilton3} 
R.~Hamilton, \textit{The Harnack estimate for the Ricci flow,} J. Diff. Geom. 37, 225--243 (1993)

\bibitem{Hamilton4}
R.~Hamilton, \textit{The formation of singularities in the Ricci flow,} Surveys in Differential Geometry, vol. II, 7--136, International Press, Somerville MA (1995)

\bibitem{Hamilton5}
R.~Hamilton, \textit{Four-manifolds with positive isotropic curvature,} Comm. Anal. Geom. 5, 1--92 (1997)

\bibitem{Huisken} 
G.~Huisken, \textit{Ricci deformation of the metric on a Riemannian manifold,} J. Diff. Geom. 21, 47--62 (1985)

\bibitem{Huisken-Sinestrari1}
G.~Huisken and C.~Sinestrari, \textit{Convexity estimates for mean curvature flow and singularities of mean convex surfaces,} Acta Math. 183, 45--70 (1999)

\bibitem{Huisken-Sinestrari2}
G.~Huisken and C.~Sinestrari, \textit{Mean curvature flow with surgeries of two-convex hypersurfaces,} Invent. Math. 175, 137--221 (2009)

\bibitem{Ivey}
T.~Ivey, \textit{Ricci solitons on compact three-manifolds,} Differential Geom. Appl. 3, 301--307 (1993)

\bibitem{Margerin1}
C.~Margerin, \textit{Pointwise pinched manifolds are space forms,} Geometric measure theory and the calculus of variations (Arcata 1984), 307--328, Proc. Sympos. Pure Math. 44, Amer. Math. Soc., Providence RI (1986)

\bibitem{Margerin2}
C.~Margerin, \textit{A sharp characterization of the smooth $4$-sphere in curvature terms,} Comm. Anal. Geom. 6, 21--65 (1998)

\bibitem{Margerin3}
C.~Margerin, \textit{D\'eformations de structures Riemanniennes,} unpublished manuscript

\bibitem{Micallef-Moore}
M.~Micallef and J.D.~Moore, \textit{Minimal two-spheres and the topology of manifolds with positive curvature on totally isotropic two-planes,} Ann. of Math. 127, 199--227 (1988)

\bibitem{Micallef-Wang}
M.~Micallef and M.Y.~Wang, \textit{Metrics with nonnegative isotropic curvature,} Duke Math. J. 72, 649--672 (1993)

\bibitem{Morgan-Tian}
J.~Morgan and G.~Tian, \textit{Ricci flow and the Poincar\'e conjecture,} Amer. Math. Soc. 2007

\bibitem{Nishikawa}
S.~Nishikawa, \textit{Deformation of Riemannian metrics and manifolds with bounded curvature ratios,} Geometric measure theory and the calculus of variations (Arcata 1984), 343--352, Proc. Sympos. Pure Math. 44, Amer. Math. Soc., Providence RI (1986)

\bibitem{Perelman1} 
G.~Perelman, \textit{The entropy formula for the Ricci flow and its geometric applications,} arxiv:0211159

\bibitem{Perelman2}
G.~Perelman, \textit{Ricci flow with surgery on three-manifolds,} arxiv:0303109

\bibitem{Perelman3}
G.~Perelman, \textit{Finite extinction time for solutions to the Ricci flow on certain three-manifolds,} arxiv:0307245

\bibitem{Wilking}
B.~Wilking, \textit{A Lie algebraic approach to Ricci flow invariant curvature conditions and Harnack inequalities,} J. Reine Angew. Math. 679, 223--247 (2013)

\bibitem{Yokota}
T.~Yokota, \textit{Complete ancient solutions to the Ricci flow with pinched curvature,} Comm. Anal. Geom. 25, 485--506 (2017)
\end{thebibliography}
\end{document}